\documentclass[a4paper]{article}

\usepackage{theorem}
\newtheorem{lem}{Lemma}[section]
\newtheorem{corollary}[lem]{Corollary}
\newtheorem{remark}[lem]{Remark}

\usepackage{amssymb}

\def\GAP{\textsf{GAP}}
\def\ATLAS{\textsf{ATLAS}}

\def\Aut{{\rm Aut}}
\def\PGL{{\rm PGL}}
  
  \def\F{{\mathbb F}}
\def\tthdump#1{#1}
\tthdump{\def\URL#1#2{\texttt{#1}}}

\def\Fix{{\rm Fix}}
\def\total{{\sigma}}
\def\prop{{P}}
\def\calC{{\cal C}}
\def\sprbound{{\cal P}}
\def\sprtotal{{\cal S}}
\def\M{{\cal M}}
\tthdump{\def\MM{{\tilde{\M}}}}

\def\fpr{{\mu }}
\def\PSL{{\rm PSL}}
\def\SL{{\rm SL}}
\def\SO{{\rm SO}}
\def\SU{{\rm SU}}
\def\GU{{\rm GU}}
\def\PGO{{\rm PGO}}
\def\PSO{{\rm PSO}}
\def\PSU{{\rm PSU}}

\def\PSp{{\rm PSp}}
\def\Sp{{\rm Sp}}

\def\GL{{\rm GL}}
\def\Aut{{\rm Aut}}

\def\PGammaL{{\rm P\hbox{$\Gamma$}L}}
\def\GammaL{{\rm \hbox{$\Gamma$}L}}
\def\GO{{\rm GO}}
\def\POmega{{\rm P\hbox{$\Omega$}}}

\begin{document}

\tthdump{\title{{\GAP} Computations Concerning Probabilistic Generation of Finite Simple Groups}}

\author{\textsc{Thomas Breuer} \\[0.5cm]
\textit{Lehrstuhl D f\"ur Mathematik} \\
\textit{RWTH, 52056 Aachen, Germany}}

\date{March 28th, 2012}

\maketitle

\abstract{
This is a collection of examples showing how
the {\GAP} system~\cite{GAP}
can be used to compute information about the probabilistic generation of
finite almost simple groups.
It includes all examples that were needed for the computational results
in~\cite{BGK}.

The purpose of this writeup is twofold.
On the one hand, the computations are documented this way.
On the other hand, the {\GAP} code shown for the examples can be used as
test input for automatic checking of the data and the functions used.

A first version of this document, which was based on {\GAP}~4.4.10,
had been accessible in the web since April~2006
and is available in the arXiv (no. 0710.3267) since October~2007.
The differences to the current version are as follows.
\begin{itemize}
\item
  The format of the {\GAP} output was adjusted to the changed behaviour
  of {\GAP}~4.5.
  This affects mainly the way how {\GAP} records are printed.
\item
  Several computations are now easier because more character tables of
  almost simple groups and maximal subgroups of such groups are available
  in the {\GAP} Character Table Library.
  (The more involved computations from the original version have been kept
  in the file.)
\item
  The computation of all conjugacy classes of a subgroup of $\POmega^+(12,3)$
  has been replaced by the computation of the conjugacy classes of elements
  of prime order in this subgroup.
\item
  The irreducible element chosen in the simple group $\POmega^-(10,3)$ has
  order $61$ not $122$.
\end{itemize}}

\textwidth16cm
\oddsidemargin0pt

\parskip 1ex plus 0.5ex minus 0.5ex
\parindent0pt

\tableofcontents


\section{Overview}

The main purpose of this note is to document the {\GAP} computations
that were carried out in order to obtain the computational results
in~\cite{BGK}.
Table~\ref{thetable} lists the simple groups among these examples.
The first column gives the group names, the second and third columns contain
a plus sign $+$ or a minus sign $-$, depending on whether the quantities
$\total(G,s)$ and $\prop(G,s)$, respectively, are less than $1/3$.
The fourth column lists the orders of elements $s$ which either prove
the $+$ signs or cover most of the cases for proving these signs.
The fifth column lists the sections in this note where the example is
treated.
The rows of the table are ordered alphabetically w.r.t.~the group names.

In order to keep this note self-contained,
we first describe the theory needed, in Section~\ref{background}.
The translation of the relevant formulae into {\GAP} functions
can be found in Section~\ref{functions}.
Then Section~\ref{chartheor} describes the computations that only require
(ordinary) character tables in the
{\GAP} Character Table Library~\cite{CTblLib1.2}.
Computations using also the groups are shown in Section~\ref{hard}.
In each of the last two sections, the examples are ordered alphabetically
w.r.t.~the names of the simple groups.

\begin{table}
\caption{Computations needed in~\cite{BGK}}\label{thetable}
\begin{center}
\tthdump{\begin{tabular}[t]{cc}}
\begin{tabular}[t]{|l|c|c|r|r|}
   \hline
      $G$           & $\total < \frac{1}{3}$
                          & $\prop < \frac{1}{3}$
                                &  $|s|$ & see \\
   \tthdump{\hline\hline}
      $A_5$         & $-$ & $-$ &      5 & \ref{A5}  \\
      $A_6$         & $-$ & $-$ &      4 & \ref{A6}  \\
      $A_7$         & $-$ & $-$ &      7 & \ref{A7}  \\
      $A_8$         & $+$ &     &     15 & \ref{easyloop}, \ref{SL} \\
      $A_9$         & $+$ &     &      9 & \ref{easyloop}, \ref{Aodd} \\
      $A_{11}$      & $+$ &     &     11 & \ref{easyloop}, \ref{Aodd} \\
      $A_{13}$      & $+$ &     &     13 & \ref{easyloop}, \ref{Aodd} \\
      $A_{15}$      & $+$ &     &     15 & \ref{Aodd} \\
      $A_{17}$      & $+$ &     &     17 & \ref{Aodd} \\
      $A_{19}$      & $+$ &     &     19 & \ref{Aodd} \\
      $A_{21}$      & $+$ &     &     21 & \ref{Aodd} \\
      $A_{23}$      & $+$ &     &     23 & \ref{Aodd} \\
   \hline
      $L_3(2)$      & $+$ &     &      7 & \ref{easyloop}, \ref{easyloopaut}, \\
                    &     &     &        & \ref{SL}, \ref{L32} \\
      $L_3(3)$      & $+$ &     &     13 & \ref{easyloop}, \ref{easyloopaut}, \\
                    &     &     &        & \ref{SL} \\
      $L_3(4)$      & $+$ &     &      7 & \ref{easyloop}, \ref{easyloopaut} \\
      $L_4(3)$      & $+$ &     &     20 & \ref{easyloop}, \ref{SL} \\
      $L_4(4)$      & $+$ &     &     85 & \ref{SL} \\
      $L_6(2)$      & $+$ &     &     63 & \ref{SL} \\
      $L_6(3)$      & $+$ &     &    182 & \ref{SL} \\
      $L_6(4)$      & $+$ &     &    455 & \ref{SL} \\
      $L_6(5)$      & $+$ &     &   1953 & \ref{SL} \\
      $L_8(2)$      & $+$ &     &    255 & \ref{SL} \\
      $L_{10}(2)$   & $+$ &     &   1023 & \ref{SL} \\
   \hline
      $M_{11}$      & $-$ & $-$ &     11 & \ref{spreadM11} \\
      $M_{12}$      & $-$ & $+$ &     10 & \ref{sporaut}, \ref{spreadM12} \\
   \hline
\end{tabular}
\tthdump{&}
\begin{tabular}[t]{|l|c|c|r|r|}
   \hline
      $G$           & $\total < \frac{1}{3}$
                          & $\prop < \frac{1}{3}$
                                &  $|s|$ & see \\
   \tthdump{\hline\hline}
      $O^+_8(2)$    & $-$ & $-$ &     15 & \ref{O8p2} \\
      $O^+_8(3)$    & $-$ & $-$ &     20 & \ref{O8p3} \\
      $O^+_8(4)$    & $+$ &     &     65 & \ref{O8p4} \\
      $O^+_{10}(2)$ & $+$ &     &     45 & \ref{O10p2} \\
      $O^+_{12}(2)$ & $+$ &     &     85 & \ref{O12p2} \\
      $O^+_{12}(3)$ & $+$ &     &    205 & \ref{O12p3} \\
      $O^-_8(2)$    & $+$ &     &     17 & \ref{easyloop} \\
      $O^-_8(3)$    & $+$ &     &     41 & \ref{O8m3} \\
      $O^-_{10}(2)$ & $+$ &     &     33 & \ref{O10m2} \\
      $O^-_{10}(3)$ & $+$ &     &    122 & \ref{O10m3} \\
      $O^-_{12}(2)$ & $+$ &     &     65 & \ref{O12m2} \\
      $O^-_{14}(2)$ & $+$ &     &    129 & \ref{O14m2} \\
      $O_7(3)$      & $-$ & $-$ &     14 & \ref{O73} \\
   \hline
      $S_4(4)$      & $+$ &     &     17 & \ref{easyloop}, \ref{easyloopaut} \\
      $S_6(2)$      & $-$ & $-$ &      9 & \ref{S62} \\
      $S_6(3)$      & $+$ &     &     14 & \ref{easyloop}, \ref{easyloopaut} \\
      $S_6(4)$      & $+$ &     &     65 & \ref{S64} \\
      $S_8(2)$      & $-$ & $-$ &     17 & \ref{S82} \\
      $S_8(3)$      & $+$ &     &     41 & \ref{S83} \\
   \hline
      $U_3(3)$      & $+$ &     &      6 & \ref{easyloop}, \ref{easyloopaut} \\
      $U_3(5)$      & $+$ &     &     10 & \ref{easyloop}, \ref{easyloopaut} \\
      $U_4(2)$      & $-$ & $-$ &      9 & \ref{U42} \\
      $U_4(3)$      & $-$ & $+$ &      7 & \ref{U43} \\
      $U_4(4)$      & $+$ &     &     65 & \ref{U44} \\
      $U_5(2)$      & $+$ &     &     11 & \ref{easyloop} \\
      $U_6(2)$      & $+$ &     &     11 & \ref{U62} \\
      $U_6(3)$      & $+$ &     &    122 & \ref{U63} \\
      $U_8(2)$      & $+$ &     &    129 & \ref{U82} \\
   \hline
\end{tabular}
\tthdump{\end{tabular}}
\end{center}
\end{table}

Contrary to~\cite{BGK}, {\ATLAS} notation is used throughout this note,
because the identifiers used for character tables in the
{\GAP} Character Table Library follow mainly the {\ATLAS}~\cite{CCN85}.
For example, we write
$L_d(q)$ for $\PSL(d,q)$,
$S_d(q)$ for $\PSp(d,q)$,
$U_d(q)$ for $\PSU(d,q)$, and
$O^+_{2d}(q)$, $O^-_{2d}(q)$, $O_{2d+1}(q)$ for
$\POmega^+(2d,q)$, $\POmega^-(2d,q)$, $\POmega(2d+1,q)$, respectively.

Furthermore, in the case of classical groups,
the character tables of the (almost) \emph{simple} groups
are considered
not the tables of the matrix groups (which are in fact often not available
in the {\GAP} Character Table Library).
Consequently, also element orders and the description of maximal subgroups
refer to the (almost) simple groups not to the matrix groups.

This note contains also several examples that are not needed for the proofs
in~\cite{BGK}.
Besides several small simple groups $G$ whose character table
is contained in the {\GAP} Character Table Library
and for which enough information is available for computing $\total(G)$,
in Section~\ref{easyloop},
a few such examples appear in individual sections.
In the table of contents, the section headers of the latter kind of examples
are marked with an asterisk $(\ast)$.

The examples use the {\GAP} Character Table Library,
the {\GAP} Library of Tables of Marks,
and the {\GAP} interface~\cite{AtlasRep} to the
{\ATLAS} of Group Representations~\cite{AGR},
so we first load these three packages in the required versions.
The {\GAP} output was adjusted to the versions shown below;
in older versions, features necessary for the computations may be missing,
and it may happen that newer versions, the behaviour is different.

Also, we force the assertion level to zero;
this is the default in interactive {\GAP} sessions
but the level is automatically set to $1$ when a file is read with
\verb|ReadTest|.

\begin{verbatim}
    gap> CompareVersionNumbers( GAPInfo.Version, "4.5.0" );
    true
    gap> LoadPackage( "ctbllib", "1.2" );
    true
    gap> LoadPackage( "tomlib", "1.2" );
    true
    gap> LoadPackage( "atlasrep", "1.5" );
    true
    gap> SetAssertionLevel( 0 );
\end{verbatim}

Some of the computations in Section~\ref{hard}
require about $800$ MB of space (on $32$ bit machines).
Therefore we check whether {\GAP} was started with sufficient maximal
memory; the command line option for this is \verb|-o 800m|.

\begin{verbatim}
    gap> max:= GAPInfo.CommandLineOptions.o;;
    gap> IsSubset( max, "m" ) and Int( Filtered( max, IsDigitChar ) ) >= 800;
    true
\end{verbatim}

Several computations involve calls to the {\GAP} function \verb|Random|.
In order to make the results of individual examples reproducible,
independent of the rest of the computations,
we reset the relevant random number generators
whenever this is appropriate.
For that, we store the initial states in the variable \verb|staterandom|,
and provide a function for resetting the random number generators.
(The \verb|Random| calls in the {\GAP} library use the two random number
generators \verb|GlobalRandomSource| and \verb|GlobalMersenneTwister|.)

\begin{verbatim}
    gap> staterandom:= [ State( GlobalRandomSource ),
    >                    State( GlobalMersenneTwister ) ];;
    gap> ResetGlobalRandomNumberGenerators:= function()
    >     Reset( GlobalRandomSource, staterandom[1] );
    >     Reset( GlobalMersenneTwister, staterandom[2] );
    > end;;
\end{verbatim}

\section{Prerequisites}\label{background}

\subsection{Theoretical Background}

Let $G$ be a finite group, $S$ the socle of $G$,
and denote by $G^{\times}$ the set of nonidentity elements in $G$.
For $s, g \in G^{\times}$, let
$\prop( g, s ):=
   |\{ h \in G; S \not\subseteq \langle s^h, g \rangle \}| / |G|$,
the proportion of elements in the class $s^G$ which fail to generate
at least $S$ with $g$;
we set $\prop( G, s ):= \max\{ \prop( g, s ); g \in G^{\times} \}$.
We are interested in finding a class $s^G$ of elements in $S$
such that $\prop( G, s ) < 1/3$ holds.

First consider $g \in S$,
and let $\M(S,s)$ denote the set of those maximal subgroups of $S$
that contain $s$.
We have
\begin{eqnarray*}
   |\{ h \in S; S \not\subseteq \langle s^h, g \rangle \}|
      & = & |\{ h \in S; \langle s, h g h^{-1} \rangle \not= S \}| \\
      & \leq & \sum_{M \in \M(S,s)} |\{ h \in S; h g h^{-1} \in M \}|
\end{eqnarray*}
Since $h g h^{-1} \in M$ holds if and only if the coset $M h$ is fixed by $g$
under the permutation action of $S$ on the right cosets of $M$ in $S$,
we get that
$|\{ h \in S; h g h^{-1} \in M \}| = |C_S(g)| \cdot |g^S \cap M|
 = |M| \cdot 1_M^S(g)$,
where $1_M^S$ is the permutation character of this action,
of degree $|S|/|M|$.
Thus
\begin{eqnarray*}
   |\{ h \in S; \langle s, h g h^{-1} \rangle \not= S \}| / |S|
      & \leq & \sum_{M \in \M(S,s)} 1_M^S(g) / 1_M^S(1) .
\end{eqnarray*}
We abbreviate the right hand side of this inequality by $\total( g, s )$,
set $\total( S, s ):= \max\{ \total( g, s ); g \in S^{\times} \}$,
and choose a transversal $T$ of $S$ in $G$.
Then $\prop( g, s ) \leq |T|^{-1} \cdot \sum_{t \in T} \total( g^t, s )$
and thus $\prop( G, s ) \leq \total( S, s )$ holds.

If $S = G$ and if $\M(G,s)$ consists of a single maximal subgroup $M$ of $G$
then equality holds,
i.e., $\prop( g, s ) = \total( g, s ) = 1_M^S(g) / 1_M^S(1)$.

The quantity $1_M^S(g) / 1_M^S(1) = |g^S \cap M| / |g^S|$
is the proportion of fixed points of $g$ in the permutation action of $S$
on the right cosets of its subgroup $M$.
This is called the \emph{fixed point ratio} of $g$ w.~r.~t.~$S/M$,
and is denoted as $\fpr(g,S/M)$.

For a subgroup $M$ of $S$, the number $n$ of $S$-conjugates of $M$
containing $s$ is equal to $|M^S| \cdot |s^S \cap M| / |s^S|$.
To see this, consider the set $\{ (s^h, M^k); h, k \in S, s^h \in M^k \}$,
the cardinality of which can be counted either as $|M^S| \cdot |s^S \cap M|$
or as $|s^S| \cdot n$.
So we get $n = |M| \cdot 1_M^S(s) / |N_S(M)|$.

If $S$ is a finite \emph{nonabelian simple} group
then each maximal subgroup in $S$ is self-normalizing,
and we have $n = 1_M^S(s)$ if $M$ is maximal.
So we can replace the summation over $\M(S,s)$ by one over a set
$\MM(S,s)$ of representatives of conjugacy classes
of maximal subgroups of $S$,
and get that
\[
   \total( g, s ) = \sum_{M \in \MM(S,s)}
     \frac{1_M^S(s) \cdot 1_M^S(g)}{1_M^S(1)}.
\]
Furthermore, we have $|\M(S,s)| = \sum_{M \in \MM(S,s)} 1_M^S(s)$.

In the following, we will often deal with the quantities
$\total(S):= \min\{ \total( S, s ); s \in S^{\times} \}$ and
$\sprtotal(S):= \lceil 1 / \total(S) - 1 \rceil$.
These values can be computed easily from the primitive
permutation characters of $S$.

Analogously, we set
$\prop(S):= \min \{ \prop( S, s ); s \in S^{\times} \}$ and
$\sprbound(S):= \lceil 1 / \prop(S) - 1 \rceil$.
Clearly we have $\prop(S) \leq \total(S)$ and $\sprbound(S) \geq \sprtotal(S)$.

One interpretation of $\sprbound(S)$ is that if this value is at least $k$
then it follows that for any $g_1, g_2, \ldots, g_k \in S^{\times}$,
there is some $s \in S$ such that $S = \langle g_i, s \rangle$,
for $1 \leq i \leq k$.
In this case, $S$ is said to have \emph{spread} at least $k$.
(Note that the lower bound $\sprtotal(S)$ for $\sprbound(S)$ can be computed
from the list of primitive permutation characters of $S$.)

Moreover, $\sprbound(S) \geq k$ implies that the element $s$ can be chosen
uniformly from a fixed conjugacy class of $S$.
This is called \emph{uniform spread} at least $k$ in~\cite{BGK}.

It is proved in~\cite{GK} that all finite simple groups have uniform spread
at least $1$,
that is, for any element $x \in S^{\times}$,
there is an element $y$ in a prescribed class of $S$
such that $G = \langle x, y \rangle$ holds.
In~\cite[Corollary~1.3]{BGK},
it is shown that all finite simple groups have uniform spread at least $2$,
and the finite simple groups with (uniform) spread exactly $2$ are listed.

Concerning the spread, it should be mentioned that the methods used here and
in~\cite{BGK} are nonconstructive in the sense that they cannot be used for
finding an element $s$ that generates $G$ together with each of the $k$
prescribed elements $g_1, g_2, \ldots, g_k$.

Now consider $g \in G \setminus S$.
Since $\prop( g^k, s ) \geq \prop( g, s )$ for any positive integer $k$,
we can assume that $g$ has prime order $p$, say.
We set $H = \langle S, g \rangle \leq G$, with $[H:S] = p$,
choose a transversal $T$ of $H$ in $G$,
let $\M^{\prime}(H,s):= \M(H,s) \setminus \{ S \}$,
and let $\MM^{\prime}(H,s)$ denote a set of representatives of
$H$-conjugacy classes of these groups.
As above,
\begin{eqnarray*}
   |\{ h \in H; S \not\subseteq \langle s^h, g \rangle \}| / |H|
    & = & |\{ h \in H; \langle s^h, g \rangle \not= H \}| / |H| \\
    & \leq & \sum_{M \in \M^{\prime}(H,s)}
                     |\{ h \in H; h g h^{-1} \in M \}| / |H| \\
    & = & \sum_{M \in \M^{\prime}(H,s)} 1_M^H(g) / 1_M^H(1) \\
    & = & \sum_{M \in \MM^{\prime}(H,s)} 1_M^H(g) \cdot 1_M^H(s) / 1_M^H(1)
\end{eqnarray*}
(Note that no summand for $M = S$ occurs,
so each group in $\MM^{\prime}(H,s)$ is self-normalizing.)
We abbreviate the right hand side by $\total(H,g,s)$,
and set $\total^{\prime}( H, s ) =
  \max\{ \total(H,g,s); g \in H \setminus S, |g| = [H:S] \}$.
Then we get
$\prop( g, s ) \leq |T|^{-1} \cdot \sum_{t \in T} \total(H^t,g^t,s)$
and thus
\[
   \prop( G, s ) \leq \max\{ \prop( S, s ),
     \max\{ \total^{\prime}( H, s );
               S \leq H \leq G, [H:S] \mbox{\rm\ prime} \} \} .
\]

For convenience, we set
$\prop^{\prime}(G,s) = \max\{ \prop(g,s); g \in G \setminus S \}$.

\subsection{Computational Criteria}\label{criteria}

The following criteria will be used when we have to show
the existence or nonexistence of $x_1, x_2, \ldots, x_k$, and $s \in G$
with the property
$\langle x_i, s \rangle = G$ for $1 \leq i \leq k$.
Note that manipulating lists of integers (representing fixed or moved points)
is much more efficient than testing whether certain permutations generate
a given group.

\begin{lem}\label{existsgoodconjugate}
Let $G$ be a finite group, $s \in G^{\times}$,
and $X = \bigcup_{M \in \M(G,s)} G/M$.
For $x_1, x_2, \ldots, x_k \in G$,
the conjugate $s^{\prime}$ of $s$ satisfies
$\langle x_i, s^{\prime} \rangle = G$ for $1 \leq i \leq k$
if and only if
$\Fix_{X}(s^{\prime}) \cap \bigcup_{i=1}^k \Fix_{X}(x_i) = \emptyset$
holds.
\end{lem}

\tthdump{\begin{proof}}
If $s^g \in U \leq G$ for some $g \in G$ then
$\Fix_{X}(U) = \emptyset$ if and only if $U = G$ holds;
note that $\Fix_{X}(G) = \emptyset$,
and $\Fix_{X}(U) = \emptyset$ implies that $U \not\subseteq h^{-1} M h$
holds for all $h \in G$ and $M \in \M(G,s)$, thus $U = G$.

Applied to $U = \langle x_i, s^{\prime} \rangle$,
we get $\langle x_i, s^{\prime} \rangle = G$ if and only if
$\Fix_{X}(s^{\prime}) \cap \Fix_{X}(x_i) = \Fix_{X}(U) = \emptyset$.
\tthdump{\end{proof}}


\begin{corollary}\label{existsgoodconjugate1}
If $\M(G,s) = \{ M \}$ in the situation of Lemma~\ref{existsgoodconjugate}
then there is a conjugate $s^{\prime}$ of $s$ that satisfies
$\langle x_i, s^{\prime} \rangle = G$ for $1 \leq i \leq k$
if and only if
$\bigcup_{i=1}^k \Fix_{X}(x_i) \not= X$.
\end{corollary}

\begin{corollary}\label{disprovespread}
Let $G$ be a finite simple group
and let $X$  be a $G$-set such that each $g \in G$ fixes at least one
point in $X$ but that $\Fix_{X}(G) = \emptyset$ holds.
If $x_1, x_2, \ldots x_k$ are elements in $G$ such that
$\bigcup_{i=1}^k \Fix_{X}(x_i) = X$ holds
then for each $s \in G$ there is at least one $i$ with
$\langle x_i, s \rangle \not= G$.
\end{corollary}

\section{{\GAP} Functions for the Computations}\label{functions}

After the introduction of general utilities in Section~\ref{utils},
we distinguish two different tasks.
Section~\ref{ctfun} introduces functions that will be used
in the following to compute $\total(g,s)$ with character-theoretic methods.
Functions for computing $\prop(g,s)$ or an upper bound for this value
will be introduced in Section~\ref{groups}.

The {\GAP} functions shown in Section~\ref{functions}
are collected in the file \verb|tst/probgen.g| that is
distributed with the {\GAP} Character Table Library,
see~\URL{http://www.math.rwth-aachen.de/\~{}Thomas.Breuer/ctbllib}{http://www.math.rwth-aachen.de/~Thomas.Breuer/ctbllib}

The functions have been designed for the examples in the later sections,
they could be generalized and optimized for other examples.
It is not our aim to provide a package for this functionality.

\subsection{General Utilities}\label{utils}

Let \verb|list| be a dense list and \verb|prop| be a unary function that returns
\verb|true| or \verb|false| when applied to the entries of \verb|list|.
\verb|PositionsProperty| returns the set of positions in \verb|list| for which
\verb|true| is returned.

\begin{verbatim}
    gap> if not IsBound( PositionsProperty ) then
    >      PositionsProperty:= function( list, prop )
    >        return Filtered( [ 1 .. Length( list ) ], i -> prop( list[i] ) );
    >      end;
    >    fi;
\end{verbatim}

%

The following two functions implement loops over ordered triples
(and quadruples, respectively) in a Cartesian product.
A prescribed function \verb|prop| is subsequently applied to the triples
(quadruples),
and if the result of this call is \verb|true| then this triple (quadruple)
is returned immediately;
if none of the calls to \verb|prop| yields \verb|true| then \verb|fail| is returned.

\begin{verbatim}
    gap> BindGlobal( "TripleWithProperty", function( threelists, prop )
    >     local i, j, k, test;
    > 
    >     for i in threelists[1] do
    >       for j in threelists[2] do
    >         for k in threelists[3] do
    >           test:= [ i, j, k ];
    >           if prop( test ) then
    >               return test;
    >           fi;
    >         od;
    >       od;
    >     od;
    > 
    >     return fail;
    > end );
    
    gap> BindGlobal( "QuadrupleWithProperty", function( fourlists, prop )
    >     local i, j, k, l, test;
    > 
    >     for i in fourlists[1] do
    >       for j in fourlists[2] do
    >         for k in fourlists[3] do
    >           for l in fourlists[4] do
    >             test:= [ i, j, k, l ];
    >             if prop( test ) then
    >               return test;
    >             fi;
    >           od;
    >         od;
    >       od;
    >     od;
    > 
    >     return fail;
    > end );
\end{verbatim}

Of course one could do better by considering \emph{un}ordered $n$-tuples
when several of the argument lists are equal,
and in practice, backtrack searches would often allow one to prune parts
of the search tree in early stages.
However, the above loops are not time critical in the examples presented
here, so the possible improvements are not worth the effort for our
purposes.

The function \verb|PrintFormattedArray| prints the matrix \verb|array| in a columnwise
formatted way.
(The only diference to the {\GAP} library function \verb|PrintArray| is that
\verb|PrintFormattedArray| chooses each column width according to the entries
only in this column not w.r.t.~the whole matrix.)

\begin{verbatim}
    gap> BindGlobal( "PrintFormattedArray", function( array )
    >      local colwidths, n, row;
    >      array:= List( array, row -> List( row, String ) );
    >      colwidths:= List( TransposedMat( array ),
    >                        col -> Maximum( List( col, Length ) ) );
    >      n:= Length( array[1] );
    >      for row in List( array, row -> List( [ 1 .. n ],
    >                   i -> FormattedString( row[i], colwidths[i] ) ) ) do
    >        Print( "  ", JoinStringsWithSeparator( row, " " ), "\n" );
    >      od;
    > end );
\end{verbatim}

Finally, \verb|CleanWorkspace| is a utility for reducing the space needed.
This is achieved by unbinding those user variables
that are not write protected and are not mentioned in the list
\verb|NeededVariables| of variable names that are bound now,
and by flushing the caches of tables of marks and character tables.

\begin{verbatim}
    gap> BindGlobal( "NeededVariables", NamesUserGVars() );
    gap> BindGlobal( "CleanWorkspace", function()
    >       local name, record;
    > 
    >       for name in Difference( NamesUserGVars(), NeededVariables ) do
    >        if not IsReadOnlyGlobal( name ) then
    >          UnbindGlobal( name );
    >        fi;
    >      od;
    >      for record in [ LIBTOMKNOWN, LIBTABLE ] do
    >        for name in RecNames( record.LOADSTATUS ) do
    >          Unbind( record.LOADSTATUS.( name ) );
    >          Unbind( record.( name ) );
    >        od;
    >      od;
    > end );
\end{verbatim}

The function \verb|PossiblePermutationCharacters| takes two ordinary character
tables \verb|sub| and \verb|tbl|,
computes the possible class fusions from \verb|sub| to \verb|tbl|,
then induces the trivial character of \verb|sub| to \verb|tbl|, w.r.t.~these fusions,
and returns the set of these class functions.
(So if \verb|sub| and \verb|tbl| are the character tables of groups $H$ and $G$,
respectively, where $H$ is a subgroup of $G$,
then the result contains the permutation character $1_H^G$.)

Note that the columns of the character tables
in the {\GAP} Character Table Library
are not explicitly associated with particular conjugacy classes of the
corresponding groups,
so from the character tables,
we can compute only \emph{possible} class fusions,
i.e., maps between the columns of two tables that satisfy certain
necessary conditions, see the section about the function
\verb|PossibleClassFusions| in the {\GAP} Reference Manual for details.
There is no problem if the permutation character is uniquely determined
by the character tables, in all other cases we give ad hoc arguments
for resolving the ambiguities.

\begin{verbatim}
    gap> BindGlobal( "PossiblePermutationCharacters", function( sub, tbl )
    >     local fus, triv;
    > 
    >     fus:= PossibleClassFusions( sub, tbl );
    >     if fus = fail then
    >       return fail;
    >     fi;
    >     triv:= [ TrivialCharacter( sub ) ];
    > 
    >     return Set( List( fus, map -> Induced( sub, tbl, triv, map )[1] ) );
    > end );
\end{verbatim}

\subsection{Character-Theoretic Computations}\label{ctfun}

We want to use the {\GAP} libraries of character tables
and of tables of marks, and proceed in three steps.

First we extract the primitive permutation characters from the library
information if this is available;
for that, we write the function \verb|PrimitivePermutationCharacters|.
Then the result can be used as the input for the function \verb|ApproxP|,
which computes the values $\total( g, s )$.
Finally, the functions \verb|ProbGenInfoSimple| and \verb|ProbGenInfoAlmostSimple|
compute $\sprtotal( G )$.

For a group $G$ whose character table $T$ is contained in the {\GAP}
character table library, the complete set of primitive permutation
characters can be easily computed if the character tables of all maximal
subgroups and their class fusions into $T$ are known
(in this case, we check whether the attribute \verb|Maxes| of $T$ is bound)
or if the table of marks of $G$ and the class fusion from $T$ into this
table of marks are known
(in this case, we check whether the attribute \verb|FusionToTom| of $T$ is bound).
If the attribute \verb|UnderlyingGroup| of $T$ is bound then this group
can be used to compute the primitive permutation characters.
The latter happens if $T$ was computed from the group object in {\GAP};
for tables in the {\GAP} character table library,
this is not the case by default.

The {\GAP} function \verb|PrimitivePermutationCharacters| tries to compute
the primitive permutation characters of a group using this information;
it returns the required list of characters if this can be computed this way,
otherwise \verb|fail| is returned.
(For convenience, we use the {\GAP} mechanism of \emph{attributes}
in order to store the permutation characters in the character table object
once they have been computed.)

\begin{verbatim}
    gap> DeclareAttribute( "PrimitivePermutationCharacters", IsCharacterTable );
    gap> InstallMethod( PrimitivePermutationCharacters,
    >     [ IsCharacterTable ],
    >     function( tbl )
    >     local maxes, tom, G;
    > 
    >     if HasMaxes( tbl ) then
    >       maxes:= List( Maxes( tbl ), CharacterTable );
    >       if ForAll( maxes, s -> GetFusionMap( s, tbl ) <> fail ) then
    >         return List( maxes, subtbl -> TrivialCharacter( subtbl )^tbl );
    >       fi;
    >     elif HasFusionToTom( tbl ) then
    >       tom:= TableOfMarks( tbl );
    >       maxes:= MaximalSubgroupsTom( tom );
    >       return PermCharsTom( tbl, tom ){ maxes[1] };
    >     elif HasUnderlyingGroup( tbl ) then
    >       G:= UnderlyingGroup( tbl );
    >       return List( MaximalSubgroupClassReps( G ),
    >                    M -> TrivialCharacter( M )^tbl );
    >     fi;
    > 
    >     return fail;
    > end );
\end{verbatim}

The function \verb|ApproxP| takes a list \verb|primitives| of primitive permutation
characters of a group $G$, say,
and the position \verb|spos| of the class $s^G$ in the character table of $G$.

Assume that the elements in \verb|primitives| have the form $1_M^G$,
for suitable maximal subgroups $M$ of $G$,
and let $\MM$ be the set of these groups $M$.
\verb|ApproxP| returns the class function $\psi$ of $G$ that is defined by
$\psi(1) = 0$ and
\[
   \psi(g) = \sum_{M \in \MM}
                 \frac{1_M^G(s) \cdot 1_M^G(g)}{1_M^G(1)}
\]
otherwise.

If \verb|primitives| contains all those primitive permutation characters $1_M^G$
of $G$ (with multiplicity according to the number of conjugacy classes
of these maximal subgroups) that do not vanish at $s$,
and if all these $M$ are self-normalizing in $G$
--this holds for example if $G$ is a finite simple group--
then $\psi(g) = \total( g, s )$ holds.

\begin{verbatim}
    gap> BindGlobal( "ApproxP", function( primitives, spos )
    >     local sum;
    > 
    >     sum:= ShallowCopy( Sum( List( primitives,
    >                                   pi -> pi[ spos ] * pi / pi[1] ) ) );
    >     sum[1]:= 0;
    > 
    >     return sum;
    > end );
\end{verbatim}

Note that for computations with permutation characters,
it would make the functions more complicated (and not more efficient)
if we would consider only elements $g$ of prime order,
and only one representative of Galois conjugate classes.

The next functions needed in this context compute $\total(S)$ and
$\sprtotal( S )$, for a simple group $S$,
and $\total^{\prime}(G,s)$ for an almost simple group $G$ with socle $S$,
respectively.

\verb|ProbGenInfoSimple| takes the character table \verb|tbl| of $S$ as its argument.
If the full list of primitive permutation characters of $S$ cannot be
computed with \verb|PrimitivePermutationCharacters| then the function returns
\verb|fail|.
Otherwise \verb|ProbGenInfoSimple| returns a list containing
the identifier of the table,
the value $\total(S)$,
the integer $\sprtotal( S )$,
a list of {\ATLAS} names of representatives of Galois families of those
classes of elements $s$ for which $\total(S) = \total( S, s )$ holds,
and the list of the corresponding cardinalities $|\M(S,s)|$.

\begin{verbatim}
    gap> BindGlobal( "ProbGenInfoSimple", function( tbl )
    >     local prim, max, min, bound, s;
    >     prim:= PrimitivePermutationCharacters( tbl );
    >     if prim = fail then
    >       return fail;
    >     fi;
    >     max:= List( [ 1 .. NrConjugacyClasses( tbl ) ],
    >                 i -> Maximum( ApproxP( prim, i ) ) );
    >     min:= Minimum( max );
    >     bound:= Inverse( min );
    >     if IsInt( bound ) then
    >       bound:= bound - 1;
    >     else
    >       bound:= Int( bound );
    >     fi;
    >     s:= PositionsProperty( max, x -> x = min );
    >     s:= List( Set( List( s, i -> ClassOrbit( tbl, i ) ) ), i -> i[1] );
    >     return [ Identifier( tbl ),
    >              min,
    >              bound,
    >              AtlasClassNames( tbl ){ s },
    >              Sum( List( prim, pi -> pi{ s } ) ) ];
    > end );
\end{verbatim}

\verb|ProbGenInfoAlmostSimple| takes the character tables \verb|tblS| and \verb|tblG|
of $S$ and $G$, and a list \verb|sposS| of class positions (w.r.t.~\verb|tblS|)
as its arguments.
It is assumed that $S$ is simple and has prime index in $G$.
If \verb|PrimitivePermutationCharacters| can compute the full list
of primitive permutation characters of $G$ then the function returns
a list containing
the identifier of \verb|tblG|,
the maximum $m$ of $\total^{\prime}( G, s )$,
for $s$ in the classes described by \verb|sposS|,
a list of {\ATLAS} names (in $G$) of the classes of elements $s$
for which this maximum is attained,
and the list of the corresponding cardinalities $|\M^{\prime}(G,s)|$.
When \verb|PrimitivePermutationCharacters| returns \verb|fail|,
also \verb|ProbGenInfoAlmostSimple| returns \verb|fail|.

\begin{verbatim}
    gap> BindGlobal( "ProbGenInfoAlmostSimple", function( tblS, tblG, sposS )
    >     local p, fus, inv, prim, sposG, outer, approx, l, max, min,
    >           s, cards, i, names;
    > 
    >     p:= Size( tblG ) / Size( tblS );
    >     if not IsPrimeInt( p )
    >        or Length( ClassPositionsOfNormalSubgroups( tblG ) ) <> 3 then
    >       return fail;
    >     fi;
    >     fus:= GetFusionMap( tblS, tblG );
    >     if fus = fail then
    >       return fail;
    >     fi;
    >     inv:= InverseMap( fus );
    >     prim:= PrimitivePermutationCharacters( tblG );
    >     if prim = fail then
    >       return fail;
    >     fi;
    >     sposG:= Set( fus{ sposS } );
    >     outer:= Difference( PositionsProperty(
    >                 OrdersClassRepresentatives( tblG ), IsPrimeInt ), fus );
    >     approx:= List( sposG, i -> ApproxP( prim, i ){ outer } );
    >     if IsEmpty( outer ) then
    >       max:= List( approx, x -> 0 );
    >     else
    >       max:= List( approx, Maximum );
    >     fi;
    >     min:= Minimum( max);
    >     s:= sposG{ PositionsProperty( max, x -> x = min ) };
    >     cards:= List( prim, pi -> pi{ s } );
    >     for i in [ 1 .. Length( prim ) ] do
    >       # Omit the character that is induced from the simple group.
    >       if ForAll( prim[i], x -> x = 0 or x = prim[i][1] ) then
    >         cards[i]:= 0;
    >       fi;
    >     od;
    >     names:= AtlasClassNames( tblG ){ s };
    >     Perform( names, ConvertToStringRep );
    > 
    >     return [ Identifier( tblG ),
    >              min,
    >              names,
    >              Sum( cards ) ];
    > end );
\end{verbatim}

The next function computes $\total(G,s)$ from
the character table \verb|tbl| of a simple or almost simple group $G$,
the name \verb|sname| of the class of $s$ in this table,
the list \verb|maxes| of the character tables of all subgroups $M$
with $M \in \M(G,s)$,
and the list \verb|numpermchars| of the numbers of possible permutation characters
induced from \verb|maxes|.
If the string \verb|"outer"| is given as an optional argument then $G$ is assumed
to be an automorphic extension of a simple group $S$, with $[G:S]$ a prime,
and $\total^{\prime}(G,s)$ is returned.
In both situations,
the result is \verb|fail| if the numbers of possible permutation characters
induced from \verb|maxes| do not coincide with the numbers prescribed in
\verb|numpermchars|.

\begin{verbatim}
    gap> BindGlobal( "SigmaFromMaxes", function( arg )
    >     local t, sname, maxes, numpermchars, prim, spos, outer;
    > 
    >     t:= arg[1];
    >     sname:= arg[2];
    >     maxes:= arg[3];
    >     numpermchars:= arg[4];
    >     prim:= List( maxes, s -> PossiblePermutationCharacters( s, t ) );
    >     spos:= Position( AtlasClassNames( t ), sname );
    >     if ForAny( [ 1 .. Length( maxes ) ],
    >                i -> Length( prim[i] ) <> numpermchars[i] ) then
    >       return fail;
    >     elif Length( arg ) = 5 and arg[5] = "outer" then
    >       outer:= Difference(
    >           PositionsProperty( OrdersClassRepresentatives( t ), IsPrimeInt ),
    >           ClassPositionsOfDerivedSubgroup( t ) );
    >       return Maximum( ApproxP( Concatenation( prim ), spos ){ outer } );
    >     else
    >       return Maximum( ApproxP( Concatenation( prim ), spos ) );
    >     fi;
    > end );
\end{verbatim}

The following function allows us to extract information about $\M(G,s)$
from the character table \verb|tbl| of $G$ and a list \verb|snames| of class positions
of $s$.
If \verb|Maxes( tbl )| is stored then the names of the character tables of the
subgroups in $\M(G,s)$ and the number of conjugates are printed,
otherwise \verb|fail| is printed.

\begin{verbatim}
    gap> BindGlobal( "DisplayProbGenMaxesInfo", function( tbl, snames )
    >     local mx, prim, i, spos, nonz, indent, j;
    > 
    >     if not HasMaxes( tbl ) then
    >       Print( Identifier( tbl ), ": fail\n" );
    >       return;
    >     fi;
    > 
    >     mx:= List( Maxes( tbl ), CharacterTable );
    >     prim:= List( mx, s -> TrivialCharacter( s )^tbl );
    >     Assert( 1, SortedList( prim ) =
    >                SortedList( PrimitivePermutationCharacters( tbl ) ) );
    >     for i in [ 1 .. Length( prim ) ] do
    >       # Deal with the case that the subgroup is normal.
    >       if ForAll( prim[i], x -> x = 0 or x = prim[i][1] ) then
    >         prim[i]:= prim[i] / prim[i][1];
    >       fi;
    >     od;
    > 
    >     spos:= List( snames,
    >                  nam -> Position( AtlasClassNames( tbl ), nam ) );
    >     nonz:= List( spos, x -> PositionsProperty( prim, pi -> pi[x] <> 0 ) );
    >     for i in [ 1 .. Length( spos ) ] do
    >       Print( Identifier( tbl ), ", ", snames[i], ": " );
    >       indent:= ListWithIdenticalEntries(
    >           Length( Identifier( tbl ) ) + Length( snames[i] ) + 4, ' ' );
    >       if not IsEmpty( nonz[i] ) then
    >         Print( Identifier( mx[ nonz[i][1] ] ), "  (",
    >                prim[ nonz[i][1] ][ spos[i] ], ")\n" );
    >         for j in [ 2 .. Length( nonz[i] ) ] do
    >           Print( indent, Identifier( mx[ nonz[i][j] ] ), "  (",
    >                prim[ nonz[i][j] ][ spos[i] ], ")\n" );
    >         od;
    >       else
    >         Print( "\n" );
    >       fi;
    >     od;
    > end );
\end{verbatim}

\subsection{Computations with Groups}\label{groups}

Here, the task is to compute $\prop(g,s)$ or $\prop(G,s)$ using explicit
computations with $G$,
where the character-theoretic bounds are not sufficient.

We start with small utilities that make the examples shorter.

For a finite solvable group \verb|G|,
the function \verb|PcConjugacyClassReps| returns a list of representatives of
the conjugacy classes of \verb|G|,
which are computed using a polycyclic presentation for \verb|G|.

\begin{verbatim}
    gap> BindGlobal( "PcConjugacyClassReps", function( G )
    >      local iso;
    > 
    >      iso:= IsomorphismPcGroup( G );
    >      return List( ConjugacyClasses( Image( iso ) ),
    >               c -> PreImagesRepresentative( iso, Representative( c ) ) );
    > end );
\end{verbatim}

For a finite group \verb|G|, a list \verb|primes| of prime integers,
and a normal subgroup \verb|N| of \verb|G|,
the function \verb|ClassesOfPrimeOrder| returns a list of those conjugacy classes
of \verb|G|
that are not contained in \verb|N| and whose elements' orders occur in \verb|primes|.

For each prime $p$ in \verb|primes|, first class representatives of order $p$
in a Sylow $p$ subgroup of \verb|G| are computed,
then the representatives in \verb|N| are discarded,
and then representatives w.~r.~t.~conjugacy in \verb|G| are computed.

(Note that this approach may be inappropriate
for example if a large elementary abelian Sylow $p$ subgroup occurs,
and if the conjugacy tests in \verb|G| are expensive, see Section~\ref{O8p4}.)

\begin{verbatim}
    gap> BindGlobal( "ClassesOfPrimeOrder", function( G, primes, N )
    >      local ccl, p, syl, reps;
    > 
    >      ccl:= [];
    >      for p in primes do
    >        syl:= SylowSubgroup( G, p );
    >        reps:= Filtered( PcConjugacyClassReps( syl ),
    >                   r -> Order( r ) = p and not r in N );
    >        Append( ccl, DuplicateFreeList( List( reps,
    >                                          r -> ConjugacyClass( G, r ) ) ) );
    >      od;
    > 
    >      return ccl;
    > end );
\end{verbatim}

The function \verb|IsGeneratorsOfTransPermGroup| takes a
\emph{transitive} permutation group \verb|G| and a list \verb|list| of elements in \verb|G|,
and returns \verb|true| if the elements in \verb|list| generate \verb|G|,
and \verb|false| otherwise.
The main point is that the return value \verb|true| requires the group
generated by \verb|list| to be transitive, and the check for transitivity
is much cheaper than the test whether this group is equal to \verb|G|.

\begin{verbatim}
    gap> BindGlobal( "IsGeneratorsOfTransPermGroup", function( G, list )
    >     local S;
    > 
    >     if not IsTransitive( G ) then
    >       Error( "<G> must be transitive on its moved points" );
    >     fi;
    >     S:= SubgroupNC( G, list );
    > 
    >     return IsTransitive( S, MovedPoints( G ) ) and Size( S ) = Size( G );
    > end );
\end{verbatim}

\verb|RatioOfNongenerationTransPermGroup| takes a \emph{transitive} permutation
group \verb|G| and two elements \verb|g| and \verb|s| of \verb|G|,
and returns the proportion $\prop(g,s)$.
(The function tests the (non)generation only for representatives of
$C_G(g)$-$C_G(s)$-double cosets.
Note that for $c_1 \in C_G(g)$, $c_2 \in C_G(s)$,
and a representative $r \in G$,
we have $\langle g^{c_1 r c_2}, s \rangle = \langle g^r, s \rangle^{c_2}$.)

\begin{verbatim}
    gap> BindGlobal( "RatioOfNongenerationTransPermGroup", function( G, g, s )
    >     local nongen, pair;
    > 
    >     if not IsTransitive( G ) then
    >       Error( "<G> must be transitive on its moved points" );
    >     fi;
    >     nongen:= 0;
    >     for pair in DoubleCosetRepsAndSizes( G, Centralizer( G, g ),
    >                     Centralizer( G, s ) ) do
    >       if not IsGeneratorsOfTransPermGroup( G, [ s, g^pair[1] ] ) then
    >         nongen:= nongen + pair[2];
    >       fi;
    >     od;
    > 
    >     return nongen / Size( G );
    > end );
\end{verbatim}

Let $G$ be a group, and let \verb|groups| be a list $[ G_1, G_2, \ldots, G_n ]$
of permutation groups such that $G_i$ describes the action of $G$ on a set
$\Omega_i$, say.
Moreover, we require that for $1 \leq i, j \leq n$,
mapping the \verb|GeneratorsOfGroup| list of $G_i$ to that of $G_j$
defines an isomorphism.
\verb|DiagonalProductOfPermGroups| takes \verb|groups| as its argument,
and returns the action of $G$ on the disjoint union of
$\Omega_1, \Omega_2, \ldots, \Omega_n$.

\begin{verbatim}
    gap> BindGlobal( "DiagonalProductOfPermGroups", function( groups )
    >     local prodgens, deg, i, gens, D, pi;
    > 
    >     prodgens:= GeneratorsOfGroup( groups[1] );
    >     deg:= NrMovedPoints( prodgens );
    >     for i in [ 2 .. Length( groups ) ] do
    >       gens:= GeneratorsOfGroup( groups[i] );
    >       D:= MovedPoints( gens );
    >       pi:= MappingPermListList( D, [ deg+1 .. deg+Length( D ) ] );
    >       deg:= deg + Length( D );
    >       prodgens:= List( [ 1 .. Length( prodgens ) ],
    >                        i -> prodgens[i] * gens[i]^pi );
    >     od;
    > 
    >     return Group( prodgens );
    > end );
\end{verbatim}


The following two functions are used to reduce checks of generation
to class representatives of maximal order.
Note that if $\langle s, g \rangle$ is a proper subgroup of $G$ then
also $\langle s^k, g \rangle$ is a proper subgroup of $G$,
so we need not check powers $s^k$ different from $s$ in this situation.

For an ordinary character table \verb|tbl|,
the function \verb|RepresentativesMaximallyCyclicSubgroups|
returns a list of class positions, containing one class of generators
for each class of maximally cyclic subgroups.

\begin{verbatim}
    gap> BindGlobal( "RepresentativesMaximallyCyclicSubgroups", function( tbl )
    >     local n, result, orders, p, pmap, i, j;
    > 
    >     # Initialize.
    >     n:= NrConjugacyClasses( tbl );
    >     result:= BlistList( [ 1 .. n ], [ 1 .. n ] );
    > 
    >     # Omit powers of smaller order.
    >     orders:= OrdersClassRepresentatives( tbl );
    >     for p in Set( Factors( Size( tbl ) ) ) do
    >       pmap:= PowerMap( tbl, p );
    >       for i in [ 1 .. n ] do
    >         if orders[ pmap[i] ] < orders[i] then
    >           result[ pmap[i] ]:= false;
    >         fi;
    >       od;
    >     od;
    > 
    >     # Omit Galois conjugates.
    >     for i in [ 1 .. n ] do
    >       if result[i] then
    >         for j in ClassOrbit( tbl, i ) do
    >           if i <> j then
    >             result[j]:= false;
    >           fi;
    >         od;
    >       fi;
    >     od;
    > 
    >     # Return the result.
    >     return ListBlist( [ 1 .. n ], result );
    > end );
\end{verbatim}

Let \verb|G| be a finite group, \verb|tbl| be the ordinary character table of \verb|G|,
and \verb|cols| be a list of class positions in \verb|tbl|,
for example the list returned by \verb|RepresentativesMaximallyCyclicSubgroups|.
The function \verb|ClassesPerhapsCorrespondingToTableColumns| returns the sublist
of those conjugacy classes of \verb|G| for which the corresponding column in \verb|tbl|
can be contained in \verb|cols|, according to element order and class size.

\begin{verbatim}
    gap> BindGlobal( "ClassesPerhapsCorrespondingToTableColumns",
    >    function( G, tbl, cols )
    >     local orders, classes, invariants;
    > 
    >     orders:= OrdersClassRepresentatives( tbl );
    >     classes:= SizesConjugacyClasses( tbl );
    >     invariants:= List( cols, i -> [ orders[i], classes[i] ] );
    > 
    >     return Filtered( ConjugacyClasses( G ),
    >         c -> [ Order( Representative( c ) ), Size(c) ] in invariants );
    > end );
\end{verbatim}

The next function computes,
for a finite group $G$ and subgroups $M_1, M_2, \ldots, M_n$ of $G$,
an upper bound for
$\max \{ \sum_{i=1}^n \fpr(g,G/M_i); g \in G \setminus Z(G) \}$.
So if the $M_i$ are the groups in $\M(G,s)$, for some $s \in G^{\times}$,
then we get an upper bound for $\total(G,s)$.

The idea is that for $M \leq G$ and $g \in G$ of order $p$, we have
\[
   \fpr(g,G/M) = |g^G \cap M| / |g^G|
               \leq \sum_{h \in C} |h^M| / |g^G|
               =    \sum_{h \in C} |h^M| \cdot |C_G(g)| / |G| ,
\]
where $C$ is a set of class representatives $h \in M$ of all those classes
that satisfy $|h| = p$ and $|C_G(h)| = |C_G(g)|$,
and in the case that $G$ is a permutation group additionally that
$h$ and $g$ move the same number of points.
(Note that it is enough to consider elements of \emph{prime} order.)

For computing the maximum of the rightmost term in this inequality,
for $g \in G \setminus Z(G)$,
we need not determine the $G$-conjugacy of class representatives in $M$.
Of course we pay the price that the result may be larger than the
leftmost term.
However,
if the maximal sum is in fact taken only over a single class representative,
we are sure that equality holds.
Thus we return a list of length two, containing the maximum of the
right hand side of the above inequality and a Boolean value indicating
whether this is equal to $\max \{ \fpr(g,G/M); g \in G \setminus Z(G) \}$
or just an upper bound.

The arguments for \verb|UpperBoundFixedPointRatios| are the group \verb|G|,
a list \verb|maxesclasses| such that the $i$-th entry is a list of conjugacy
classes of $M_i$, which covers all classes of prime element order in $M_i$,
and either \verb|true| or \verb|false|, where \verb|true| means that the \emph{exact} value
of $\total(G,s)$ is computed, not just an upper bound;
this can be much more expensive because of the conjugacy tests in $G$
that may be necessary.
(We try to reduce the number of conjugacy tests in this case,
the second half of the code is not completely straightforward.
The special treatment of conjugacy checks for elements with the same sets
of fixed points is essential in the computation of $\total^{\prime}(G,s)$
for $G = \PGL(6,4)$;
the critical input line is \verb|ApproxPForOuterClassesInGL( 6, 4 )|,
see Section~\ref{SLaut}.
Currently the standard {\GAP} conjugacy test for
an element of order three and its inverse in $G \setminus G^{\prime}$
requires hours of CPU time, whereas the check for existence of a conjugating
element in the stabilizer of the common set of fixed points of the two
elements is almost free of charge.)

\verb|UpperBoundFixedPointRatios| can be used to compute $\total^{\prime}(G,s)$
in the case that $G$ is an automorphic extension of a simple group $S$,
with $[G:S] = p$ a prime;
if $\M^{\prime}(G,s) = \{ M_1, M_2, \ldots, M_n \}$ then
the $i$-th entry of \verb|maxesclasses| must contain only
the classes of element order $p$ in $M_i \setminus (M_i \cap S)$.

\begin{verbatim}
    gap> BindGlobal( "UpperBoundFixedPointRatios",
    >    function( G, maxesclasses, truetest )
    >     local myIsConjugate, invs, info, c, r, o, inv, pos, sums, max, maxpos,
    >           maxlen, reps, split, i, found, j;
    > 
    >     myIsConjugate:= function( G, x, y )
    >       local movx, movy;
    > 
    >       movx:= MovedPoints( x );
    >       movy:= MovedPoints( y );
    >       if movx = movy then
    >         G:= Stabilizer( G, movx, OnSets );
    >       fi;
    >       return IsConjugate( G, x, y );
    >     end;
    > 
    >     invs:= [];
    >     info:= [];
    > 
    >     # First distribute the classes according to invariants.
    >     for c in Concatenation( maxesclasses ) do
    >       r:= Representative( c );
    >       o:= Order( r );
    >       # Take only prime order representatives.
    >       if IsPrimeInt( o ) then
    >         inv:= [ o, Size( Centralizer( G, r ) ) ];
    >         # Omit classes that are central in `G'.
    >         if inv[2] <> Size( G ) then
    >           if IsPerm( r ) then
    >             Add( inv, NrMovedPoints( r ) );
    >           fi;
    >           pos:= First( [ 1 .. Length( invs ) ], i -> inv = invs[i] );
    >           if pos = fail then
    >             # This class is not `G'-conjugate to any of the previous ones.
    >             Add( invs, inv );
    >             Add( info, [ [ r, Size( c ) * inv[2] ] ] );
    >           else
    >             # This class may be conjugate to an earlier one.
    >             Add( info[ pos ], [ r, Size( c ) * inv[2] ] );
    >           fi;
    >         fi;
    >       fi;
    >     od;
    > 
    >     if info = [] then
    >       return [ 0, true ];
    >     fi;
    > 
    >     repeat
    >       # Compute the contributions of the classes with the same invariants.
    >       sums:= List( info, x -> Sum( List( x, y -> y[2] ) ) );
    >       max:= Maximum( sums );
    >       maxpos:= Filtered( [ 1 .. Length( info ) ], i -> sums[i] = max );
    >       maxlen:= List( maxpos, i -> Length( info[i] ) );
    > 
    >       # Split the sets with the same invariants if necessary
    >       # and if we want to compute the exact value.
    >       if truetest and not 1 in maxlen then
    >         # Make one conjugacy test.
    >         pos:= Position( maxlen, Minimum( maxlen ) );
    >         reps:= info[ maxpos[ pos ] ];
    >         if myIsConjugate( G, reps[1][1], reps[2][1] ) then
    >           # Fuse the two classes.
    >           reps[1][2]:= reps[1][2] + reps[2][2];
    >           reps[2]:= reps[ Length( reps ) ];
    >           Unbind( reps[ Length( reps ) ] );
    >         else
    >           # Split the list. This may require additional conjugacy tests.
    >           Unbind( info[ maxpos[ pos ] ] );
    >           split:= [ reps[1], reps[2] ];
    >           for i in [ 3 .. Length( reps ) ] do
    >             found:= false;
    >             for j in split do
    >               if myIsConjugate( G, reps[i][1], j[1] ) then
    >                 j[2]:= reps[i][2] + j[2];
    >                 found:= true;
    >                 break;
    >               fi;
    >             od;
    >             if not found then
    >               Add( split, reps[i] );
    >             fi;
    >           od;
    > 
    >           info:= Compacted( Concatenation( info,
    >                                            List( split, x -> [ x ] ) ) );
    >         fi;
    >       fi;
    >     until 1 in maxlen or not truetest;
    > 
    >     return [ max / Size( G ), 1 in maxlen ];
    > end );
\end{verbatim}


Suppose that $C_1, C_2, C_3$ are conjugacy classes in $G$,
and that we have to prove,
for each $(x_1, x_2, x_3) \in C_1 \times C_2 \times C_3$,
the existence of an element $s$ in a prescribed class $C$ of $G$ such that
$\langle x_1, s \rangle = \langle x_2, s \rangle = \langle x_2, s \rangle = G$
holds.

We have to check only representatives under the conjugation action of $G$
on $C_1 \times C_2 \times C_3$.
For each representative, we try a prescribed number of random elements in $C$.
If this is successful then we are done.
The following two functions implement this idea.

For a group $G$ and a list $[ g_1, g_2, \ldots, g_n ]$ of elements in $G$,
\verb|OrbitRepresentativesProductOfClasses| returns a list
$R(G, g_1, g_2, \ldots, g_n)$ of representatives of $G$-orbits
on the Cartesian product $g_1^G \times g_2^G \times \cdots \times g_n^G$.

The idea behind this function is to choose $R(G, g_1) = \{ ( g_1 ) \}$
in the case $n = 1$,
and, for $n > 1$,
\begin{eqnarray*}
   R(G, g_1, g_2, \ldots, g_n) & = & \{ (h_1, h_2, \ldots, h_n) \mid
          (h_1, h_2, \ldots, h_{n-1}) \in R(G, g_1, g_2, \ldots, g_{n-1}), \\
          & & h_n = g_n^d, \mbox{\rm\ for\ } d \in D \} ,
\end{eqnarray*}
where $D$ is a set of representatives of double cosets
$C_G(g_n) \setminus G / \cap_{i=1}^{n-1} C_G(h_i)$.


\begin{verbatim}
    gap> BindGlobal( "OrbitRepresentativesProductOfClasses",
    >    function( G, classreps )
    >     local cents, n, orbreps;
    > 
    >     cents:= List( classreps, x -> Centralizer( G, x ) );
    >     n:= Length( classreps );
    > 
    >     orbreps:= function( reps, intersect, pos )
    >       if pos > n then
    >         return [ reps ];
    >       fi;
    >       return Concatenation( List(
    >           DoubleCosetRepsAndSizes( G, cents[ pos ], intersect ),
    >             r -> orbreps( Concatenation( reps, [ classreps[ pos ]^r[1] ] ),
    >                  Intersection( intersect, cents[ pos ]^r[1] ), pos+1 ) ) );
    >     end;
    > 
    >     return orbreps( [ classreps[1] ], cents[1], 2 );
    > end );
\end{verbatim}


The function \verb|RandomCheckUniformSpread| takes
a transitive permutation group $G$,
a list of class representatives $g_i \in G$, an element $s \in G$,
and a positive integer $N$.
The return value is \verb|true| if for each representative of $G$-orbits
on the product of the classes $g_i^G$,
a good conjugate of $s$ is found in at most $N$ random tests.

\begin{verbatim}
    gap> BindGlobal( "RandomCheckUniformSpread", function( G, classreps, s, try )
    >     local elms, found, i, conj;
    > 
    >     if not IsTransitive( G, MovedPoints( G ) ) then
    >       Error( "<G> must be transitive on its moved points" );
    >     fi;
    > 
    >     # Compute orbit representatives of G on the direct product,
    >     # and try to find a good conjugate of s for each representative.
    >     for elms in OrbitRepresentativesProductOfClasses( G, classreps ) do
    >       found:= false;
    >       for i in [ 1 .. try ] do
    >         conj:= s^Random( G );
    >         if ForAll( elms,
    >               x -> IsGeneratorsOfTransPermGroup( G, [ x, conj ] ) ) then
    >           found:= true;
    >           break;
    >         fi;
    >       od;
    >       if not found then
    >         return elms;
    >       fi;
    >     od;
    > 
    >     return true;
    > end );
\end{verbatim}

Of course this approach is not suitable for \emph{dis}proving the existence
of $s$, but it is much cheaper than an exhaustive search in the class $C$.
(Typically, $|C|$ is large whereas the $|C_i|$ are small.)


%

The following function can be used to verify that a given $n$-tuple
$(x_1, x_2, \ldots, x_n)$ of elements in a group $G$ has the property
that for all elements $g \in G$, at least one $x_i$ satisfies
$\langle x_i, g \rangle$.
The arguments are a transitive permutation group $G$,
a list of class representatives in $G$, and the $n$-tuple in question.
The return value is a conjugate $g$ of the given representatives
that has the property if such an element exists,
and \verb|fail| otherwise.

\begin{verbatim}
    gap> BindGlobal( "CommonGeneratorWithGivenElements",
    >    function( G, classreps, tuple )
    >     local inter, rep, repcen, pair;
    > 
    >     if not IsTransitive( G, MovedPoints( G ) ) then
    >       Error( "<G> must be transitive on its moved points" );
    >     fi;
    > 
    >     inter:= Intersection( List( tuple, x -> Centralizer( G, x ) ) );
    >     for rep in classreps do
    >       repcen:= Centralizer( G, rep );
    >       for pair in DoubleCosetRepsAndSizes( G, repcen, inter ) do
    >         if ForAll( tuple,
    >            x -> IsGeneratorsOfTransPermGroup( G, [ x, rep^pair[1] ] ) ) then
    >           return rep;
    >         fi;
    >       od;
    >     od;
    > 
    >     return fail;
    > end );
\end{verbatim}

\section{Character-Theoretic Computations}\label{chartheor}

In this section, we apply the functions introduced in Section~\ref{ctfun}
to the character tables of simple groups
that are available in the {\GAP} Character Table Library.

Our first examples are the sporadic simple groups, in Section~\ref{spor},
then their automorphism groups are considered in Section~\ref{sporaut}.

Then we consider those other simple groups for which {\GAP} provides
enough information for automatically computing an upper bound on
$\total(G,s)$ --see Section~\ref{easyloop}--
and their automorphic extensions --see Section~\ref{easyloopaut}.

After that, individual groups are considered.

\subsection{Sporadic Simple Groups}\label{spor}

The {\GAP} Character Table Library contains the tables of maximal subgroups
of all sporadic simple groups except $B$ and $M$,
so all primitive permutation characters can be computed via the function
\verb|PrimitivePermutationCharacters| for $24$ of the $26$ sporadic simple groups.

\begin{verbatim}
    gap> sporinfo:= [];;
    gap> spornames:= AllCharacterTableNames( IsSporadicSimple, true,
    >                                        IsDuplicateTable, false );;
    gap> for tbl in List( spornames, CharacterTable ) do
    >      info:= ProbGenInfoSimple( tbl );
    >      if info <> fail then
    >        Add( sporinfo, info );
    >      fi;
    >    od;
\end{verbatim}

We show the result as a formatted table.

\begin{verbatim}
    gap> PrintFormattedArray( sporinfo );
       Co1    421/1545600         3671        [ "35A" ]    [ 4 ]
       Co2          1/270          269        [ "23A" ]    [ 1 ]
       Co3        64/6325           98        [ "21A" ]    [ 4 ]
       F3+ 1/269631216855 269631216854        [ "29A" ]    [ 1 ]
      Fi22         43/585           13        [ "16A" ]    [ 7 ]
      Fi23   2651/2416635          911        [ "23A" ]    [ 2 ]
        HN        4/34375         8593        [ "19A" ]    [ 1 ]
        HS        64/1155           18        [ "15A" ]    [ 2 ]
        He          3/595          198        [ "14C" ]    [ 3 ]
        J1           1/77           76        [ "19A" ]    [ 1 ]
        J2           5/28            5        [ "10C" ]    [ 3 ]
        J3          2/153           76        [ "19A" ]    [ 2 ]
        J4   1/1647124116   1647124115        [ "29A" ]    [ 1 ]
        Ly     1/35049375     35049374        [ "37A" ]    [ 1 ]
       M11            1/3            2        [ "11A" ]    [ 1 ]
       M12            1/3            2        [ "10A" ]    [ 3 ]
       M22           1/21           20        [ "11A" ]    [ 1 ]
       M23         1/8064         8063        [ "23A" ]    [ 1 ]
       M24       108/1265           11        [ "21A" ]    [ 2 ]
       McL      317/22275           70 [ "15A", "30A" ] [ 3, 3 ]
        ON       10/30723         3072        [ "31A" ]    [ 2 ]
        Ru         1/2880         2879        [ "29A" ]    [ 1 ]
       Suz       141/5720           40        [ "14A" ]    [ 3 ]
        Th       2/267995       133997 [ "27A", "27B" ] [ 2, 2 ]
\end{verbatim}

We see that in all these cases, $\total(G) < 1/2$ and thus
$\sprbound( G ) \geq 2$,
and all sporadic simple groups $G$ except $G = M_{11}$ and $G = M_{12}$
satisfy $\total(G) < 1/3$.
See~\ref{spreadM11} and~\ref{spreadM12} for a proof that also these
two groups have uniform spread at least three.

The structures and multiplicities of the maximal subgroups containing $s$
are as follows.

\begin{verbatim}
    gap> for entry in sporinfo do
    >      DisplayProbGenMaxesInfo( CharacterTable( entry[1] ), entry[4] );
    > od;
    Co1, 35A: (A5xJ2):2  (1)
              (A6xU3(3)):2  (2)
              (A7xL2(7)):2  (1)
    Co2, 23A: M23  (1)
    Co3, 21A: U3(5).3.2  (2)
              L3(4).D12  (1)
              s3xpsl(2,8).3  (1)
    F3+, 29A: 29:14  (1)
    Fi22, 16A: 2^10:m22  (1)
               (2x2^(1+8)):U4(2):2  (1)
               2F4(2)'  (4)
               2^(5+8):(S3xA6)  (1)
    Fi23, 23A: 2..11.m23  (1)
               L2(23)  (1)
    HN, 19A: U3(8).3_1  (1)
    HS, 15A: A8.2  (1)
             5:4xa5  (1)
    He, 14C: 2^1+6.psl(3,2)  (1)
             7^2:2psl(2,7)  (1)
             7^(1+2):(S3x3)  (1)
    J1, 19A: 19:6  (1)
    J2, 10C: 2^1+4b:a5  (1)
             a5xd10  (1)
             5^2:D12  (1)
    J3, 19A: L2(19)  (1)
             J3M3  (1)
    J4, 29A: frob  (1)
    Ly, 37A: 37:18  (1)
    M11, 11A: L2(11)  (1)
    M12, 10A: A6.2^2  (1)
              M12M4  (1)
              2xS5  (1)
    M22, 11A: L2(11)  (1)
    M23, 23A: 23:11  (1)
    M24, 21A: L3(4).3.2_2  (1)
              2^6:(psl(3,2)xs3)  (1)
    McL, 15A: 3^(1+4):2S5  (1)
              2.A8  (1)
              5^(1+2):3:8  (1)
    McL, 30A: 3^(1+4):2S5  (1)
              2.A8  (1)
              5^(1+2):3:8  (1)
    ON, 31A: L2(31)  (1)
             ONM8  (1)
    Ru, 29A: L2(29)  (1)
    Suz, 14A: J2.2  (2)
              (a4xpsl(3,4)):2  (1)
    Th, 27A: ThN3B  (1)
             ThM7  (1)
    Th, 27B: ThN3B  (1)
             ThM7  (1)
\end{verbatim}

For the remaining two sporadic simple groups, $B$ and $M$,
we choose suitable elements $s$.
If $G = B$ and $s \in G$ is of order $47$ then, by~\cite{Wil99},
$\M(G,s) = \{ 47:23 \}$.

\begin{verbatim}
    gap> SigmaFromMaxes( CharacterTable( "B" ), "47A",
    >        [ CharacterTable( "47:23" ) ], [ 1 ] );
    1/174702778623598780219392000000
\end{verbatim}

If $G = M$ and $s \in G$ is of order $59$ then, by~\cite{HW04},
$\M(G,s) = \{ L_2(59) \}$.
In this case, the permutation character is not uniquely determined by the
character tables, but all possibilities lead to the same value for
$\total(G)$.

\begin{verbatim}
    gap> t:= CharacterTable( "M" );;
    gap> s:= CharacterTable( "L2(59)" );;
    gap> pi:= PossiblePermutationCharacters( s, t );;
    gap> Length( pi );
    5
    gap> spos:= Position( OrdersClassRepresentatives( t ), 59 );
    152
    gap> Set( List( pi, x -> Maximum( ApproxP( [ x ], spos ) ) ) );
    [ 1/3385007637938037777290625 ]
\end{verbatim}

Essentially the same approach is taken in~\cite{GM01}.
However, there $s$ is restricted to classes of prime order.
Thus the results in the above table are better for $J_2$, $HS$, $M_{24}$,
$McL$, $He$, $Suz$, $Co_3$, $Fi_{22}$, $Ly$, $Th$, $Co_1$, and $J_4$.
Besides that, the value $10\,999$ claimed in~\cite{GM01}
for $\sprtotal( HN )$ is not correct.


\subsection{Automorphism Groups of Sporadic Simple Groups}\label{sporaut}

Next we consider the automorphism groups of the sporadic simple groups.
There are exactly $12$ cases where nontrivial outer automorphisms exist,
and then the simple group $S$ has index $2$ in its automorphism group $G$.

\begin{verbatim}
    gap> sporautnames:= AllCharacterTableNames( IsSporadicSimple, true,
    >                       IsDuplicateTable, false,
    >                       OfThose, AutomorphismGroup );;
    gap> sporautnames:= Difference( sporautnames, spornames );
    [ "F3+.2", "Fi22.2", "HN.2", "HS.2", "He.2", "J2.2", "J3.2", "M12.2",
      "M22.2", "McL.2", "ON.2", "Suz.2" ]
\end{verbatim}

First we compute the values $\total^{\prime}(G,s)$,
for the same $s \in S$ that were chosen for the simple group $S$
in Section~\ref{spor}.

For six of the groups $G$ in question,
the character tables of all maximal subgroups are available in the {\GAP}
Character Table Library.
In these cases, the values $\total^{\prime}( G, s )$
can be computed using \verb|ProbGenInfoAlmostSimple|.

\emph{(The above statement can meanwhile be replaced by the statement that
the character tables of all maximal subgroups are available for all twelve
groups.
We show the table results for all these groups but keep the individual
computations from the original computations.)}

\begin{verbatim}
    gap> sporautinfo:= [];;
    gap> fails:= [];;
    gap> for name in sporautnames do
    >      tbl:= CharacterTable( name{ [ 1 .. Position( name, '.' ) - 1 ] } );
    >      tblG:= CharacterTable( name );
    >      info:= ProbGenInfoSimple( tbl );
    >      info:= ProbGenInfoAlmostSimple( tbl, tblG,
    >          List( info[4], x -> Position( AtlasClassNames( tbl ), x ) ) );
    >      if info = fail then
    >        Add( fails, name );
    >      else
    >        Add( sporautinfo, info );
    >      fi;
    >    od;
    gap> PrintFormattedArray( sporautinfo );
       F3+.2         0         [ "29AB" ]    [ 1 ]
      Fi22.2  251/3861         [ "16AB" ]    [ 7 ]
        HN.2    1/6875         [ "19AB" ]    [ 1 ]
        HS.2    36/275          [ "15A" ]    [ 2 ]
        He.2   37/9520         [ "14CD" ]    [ 3 ]
        J2.2      1/15         [ "10CD" ]    [ 3 ]
        J3.2    1/1080         [ "19AB" ]    [ 1 ]
       M12.2      4/99          [ "10A" ]    [ 1 ]
       M22.2      1/21         [ "11AB" ]    [ 1 ]
       McL.2      1/63 [ "15AB", "30AB" ] [ 3, 3 ]
        ON.2   1/84672         [ "31AB" ]    [ 1 ]
       Suz.2 661/46332          [ "14A" ]    [ 3 ]
\end{verbatim}

Note that for $S = McL$, the bound $\total^{\prime}(G,s)$ for $G = S.2$
(in the second column) is worse than the bound for the simple group $S$.

The structures and multiplicities of the maximal subgroups containing $s$
are as follows.

\begin{verbatim}
    gap> for entry in sporautinfo do
    >      DisplayProbGenMaxesInfo( CharacterTable( entry[1] ), entry[3] );
    > od;
    F3+.2, 29AB: F3+  (1)
                 frob  (1)
    Fi22.2, 16AB: Fi22  (1)
                  Fi22.2M4  (1)
                  (2x2^(1+8)):(U4(2):2x2)  (1)
                  2F4(2)'.2  (4)
                  2^(5+8):(S3xS6)  (1)
    HN.2, 19AB: HN  (1)
                U3(8).6  (1)
    HS.2, 15A: HS  (1)
               S8x2  (1)
               5:4xS5  (1)
    He.2, 14CD: He  (1)
                2^(1+6)_+.L3(2).2  (1)
                7^2:2.L2(7).2  (1)
                7^(1+2):(S3x6)  (1)
    J2.2, 10CD: J2  (1)
                2^(1+4).S5  (1)
                (A5xD10).2  (1)
                5^2:(4xS3)  (1)
    J3.2, 19AB: J3  (1)
                19:18  (1)
    M12.2, 10A: M12  (1)
                (2^2xA5):2  (1)
    M22.2, 11AB: M22  (1)
                 L2(11).2  (1)
    McL.2, 15AB: McL  (1)
                 3^(1+4):4S5  (1)
                 Isoclinic(2.A8.2)  (1)
                 5^(1+2):(24:2)  (1)
    McL.2, 30AB: McL  (1)
                 3^(1+4):4S5  (1)
                 Isoclinic(2.A8.2)  (1)
                 5^(1+2):(24:2)  (1)
    ON.2, 31AB: ON  (1)
                31:30  (1)
    Suz.2, 14A: Suz  (1)
                J2.2x2  (2)
                (A4xL3(4):2_3):2  (1)
\end{verbatim}

Note that the maximal subgroups $L_2(19)$ of $J_3$ do not extend
to $J_3.2$ and that a class of maximal subgroups of the type $19:18$ appears
in $J_3.2$ whose intersection with $J_3$ is not maximal in $J_3$.
Similarly, the maximal subgroups $A_6.2^2$ of $M_{12}$ do not extend
to $M_{12}.2$.

For the other six groups, we use individual computations.

In the case $S = Fi_{24}^{\prime}$, the unique maximal subgroup $29:14$
that contains an element $s$ of order $29$ extends to a group of the type
$29:28$ in $Fi_{24}$,
which is a nonsplit extension of $29:14$.

\begin{verbatim}
    gap> SigmaFromMaxes( CharacterTable( "Fi24'.2" ), "29AB",
    >        [ CharacterTable( "29:28" ) ], [ 1 ], "outer" );
    0
\end{verbatim}

In the case $S = Fi_{22}$, there are four classes of maximal subgroups
that contain $s$ of order $16$.
They extend to $G = Fi_{22}.2$,
and none of the \emph{novelties} in $G$ (i.~e., subgroups of $G$ that are
maximal in $G$ but whose intersections with $S$ are not maximal in $S$)
contains $s$, cf.~\cite[p.~163]{CCN85}.

\begin{verbatim}
    gap> 16 in OrdersClassRepresentatives( CharacterTable( "U4(2).2" ) );
    false
    gap> 16 in OrdersClassRepresentatives( CharacterTable( "G2(3).2" ) );
    false
\end{verbatim}

The character tables of three of the four extensions are available in the
{\GAP} Character Table Library.
The permutation character on the cosets of the fourth extension can be
obtained as the extension of the permutation character of $S$ on the
cosets of its maximal subgroup of the type $2^{5+8}:(S_3 \times A_6)$.

\begin{verbatim}
    gap> t2:= CharacterTable( "Fi22.2" );;
    gap> prim:= List( [ "Fi22.2M4", "(2x2^(1+8)):(U4(2):2x2)", "2F4(2)" ],
    >        n -> PossiblePermutationCharacters( CharacterTable( n ), t2 ) );;
    gap> t:= CharacterTable( "Fi22" );;
    gap> pi:= PossiblePermutationCharacters(
    >             CharacterTable( "2^(5+8):(S3xA6)" ), t );
    [ Character( CharacterTable( "Fi22" ), [ 3648645, 56133, 10629, 2245, 567, 
          729, 405, 81, 549, 165, 133, 37, 69, 20, 27, 81, 9, 39, 81, 19, 1, 13, 
          33, 13, 1, 0, 13, 13, 5, 1, 0, 0, 0, 8, 4, 0, 0, 9, 3, 15, 3, 1, 1, 1, 
          1, 3, 3, 1, 0, 0, 0, 2, 1, 1, 0, 0, 0, 0, 0, 0, 0, 0, 1, 1, 2 ] ) ]
    gap> torso:= CompositionMaps( pi[1], InverseMap( GetFusionMap( t, t2 ) ) );
    [ 3648645, 56133, 10629, 2245, 567, 729, 405, 81, 549, 165, 133, 37, 69, 20, 
      27, 81, 9, 39, 81, 19, 1, 13, 33, 13, 1, 0, 13, 13, 5, 1, 0, 0, 0, 8, 4, 0, 
      9, 3, 15, 3, 1, 1, 1, 3, 3, 1, 0, 0, 2, 1, 0, 0, 0, 0, 0, 0, 1, 1, 2 ]
    gap> ext:= PermChars( t2, rec( torso:= torso ) );;
    gap> Add( prim, ext );
    gap> prim:= Concatenation( prim );;  Length( prim );
    4
    gap> spos:= Position( OrdersClassRepresentatives( t2 ), 16 );;
    gap> List( prim, x -> x[ spos ] );
    [ 1, 1, 4, 1 ]
    gap> sigma:= ApproxP( prim, spos );;
    gap> Maximum( sigma{ Difference( PositionsProperty(
    >                        OrdersClassRepresentatives( t2 ), IsPrimeInt ),
    >                        ClassPositionsOfDerivedSubgroup( t2 ) ) } );
    251/3861
\end{verbatim}

In the case $S = HN$, the unique maximal subgroup $U_3(8).3$
that contains the fixed element $s$ of order $19$
extends to a group of the type $U_3(8).6$ in $HN.2$.

\begin{verbatim}
    gap> SigmaFromMaxes( CharacterTable( "HN.2" ), "19AB",
    >        [ CharacterTable( "U3(8).6" ) ], [ 1 ], "outer" );
    1/6875
\end{verbatim}

In the case $S = HS$, there are two classes of maximal subgroups
that contain $s$ of order $15$.
They extend to $G = HS.2$,
and none of the novelties in $G$ contains $s$ (cf.~\cite[p.~80]{CCN85}).

\begin{verbatim}
    gap> SigmaFromMaxes( CharacterTable( "HS.2" ), "15A",
    >      [ CharacterTable( "S8x2" ),
    >        CharacterTable( "5:4" ) * CharacterTable( "A5.2" ) ], [ 1, 1 ],
    >      "outer" );
    36/275
\end{verbatim}

In the case $S = He$, there are three classes of maximal subgroups
that contain $s$ in the class {\tt 14C}.
They extend to $G = He.2$,
and none of the novelties in $G$ contains $s$ (cf.~\cite[p.~104]{CCN85}).
We compute the extensions of the corresponding primitive permutation
characters of $S$.

\begin{verbatim}
    gap> t:= CharacterTable( "He" );;
    gap> t2:= CharacterTable( "He.2" );;
    gap> prim:= PrimitivePermutationCharacters( t );;
    gap> spos:= Position( AtlasClassNames( t ), "14C" );;
    gap> prim:= Filtered( prim, x -> x[ spos ] <> 0 );;
    gap> map:= InverseMap( GetFusionMap( t, t2 ) );;
    gap> torso:= List( prim, pi -> CompositionMaps( pi, map ) );
    [ [ 187425, 945, 449, 0, 21, 21, 25, 25, 0, 0, 5, 0, 0, 7, 1, 0, 0, 1, 0, 1, 
          0, 0, 0, 0, 0, 0 ], 
      [ 244800, 0, 64, 0, 84, 0, 0, 16, 0, 0, 4, 24, 45, 3, 4, 0, 0, 0, 0, 1, 0, 
          0, 0, 0, 0, 0 ], 
      [ 652800, 0, 512, 120, 72, 0, 0, 0, 0, 0, 8, 8, 22, 1, 0, 0, 0, 0, 0, 1, 0, 
          0, 1, 1, 2, 0 ] ]
    gap> ext:= List( torso, x -> PermChars( t2, rec( torso:= x ) ) );
    [ [ Character( CharacterTable( "He.2" ), [ 187425, 945, 449, 0, 21, 21, 25, 
              25, 0, 0, 5, 0, 0, 7, 1, 0, 0, 1, 0, 1, 0, 0, 0, 0, 0, 0, 315, 15, 
              0, 0, 3, 7, 7, 3, 0, 0, 0, 1, 1, 0, 1, 1, 0, 0, 0 ] ) ], 
      [ Character( CharacterTable( "He.2" ), [ 244800, 0, 64, 0, 84, 0, 0, 16, 0, 
              0, 4, 24, 45, 3, 4, 0, 0, 0, 0, 1, 0, 0, 0, 0, 0, 0, 360, 0, 0, 0, 
              6, 0, 0, 0, 0, 0, 3, 2, 2, 0, 0, 0, 0, 0, 0 ] ) ], 
      [ Character( CharacterTable( "He.2" ), [ 652800, 0, 512, 120, 72, 0, 0, 0, 
              0, 0, 8, 8, 22, 1, 0, 0, 0, 0, 0, 1, 0, 0, 1, 1, 2, 0, 480, 0, 120, 
              0, 12, 0, 0, 0, 0, 0, 4, 0, 0, 0, 0, 0, 0, 1, 1 ] ) ] ]
    gap> spos:= Position( AtlasClassNames( t2 ), "14CD" );;
    gap> sigma:= ApproxP( Concatenation( ext ), spos );;
    gap> Maximum( sigma{ Difference( PositionsProperty(
    >                        OrdersClassRepresentatives( t2 ), IsPrimeInt ),
    >                        ClassPositionsOfDerivedSubgroup( t2 ) ) } );
    37/9520
\end{verbatim}

In the case $S = O'N$, the two classes of maximal subgroups of the type
$L_2(31)$ do not extend to $G = O'N.2$, and a class of novelties of the
structure $31:30$ appears (see~\cite[p.~132]{CCN85}).

\begin{verbatim}
    gap> SigmaFromMaxes( CharacterTable( "ON.2" ), "31AB",
    >        [ CharacterTable( "P:Q", [ 31, 30 ] ) ], [ 1 ], "outer" );
    1/84672
\end{verbatim}

Now we consider also $\total(G,\hat{s})$,
for suitable $\hat{s} \in G \setminus S$;
this yields lower bounds for the spread of the nonsimple groups $G$.
(These results are shown in the last two columns of~\cite[Table~9]{BGK}.)

As above, we use the known character tables of the maximal subgroups
in order to compute the optimal choice for $\hat{s} \in G \setminus S$.
(We may use the function \verb|ProbGenInfoSimple| although the groups are not
simple; all we need is that the relevant maximal subgroups are
self-normalizing.)

\begin{verbatim}
    gap> sporautinfo2:= [];;
    gap> for name in List( sporautinfo, x -> x[1] ) do
    >      Add( sporautinfo2, ProbGenInfoSimple( CharacterTable( name ) ) );
    >    od;
    gap> PrintFormattedArray( sporautinfo2 );
       F3+.2    19/5684  299        [ "42E" ]   [ 10 ] 
      Fi22.2 1165/20592   17        [ "24G" ]    [ 3 ] 
        HN.2     1/1425 1424        [ "24B" ]    [ 4 ] 
        HS.2     21/550   26        [ "20C" ]    [ 4 ] 
        He.2    33/4165  126        [ "24A" ]    [ 2 ] 
        J2.2       1/15   14        [ "14A" ]    [ 1 ] 
        J3.2   77/10260  133        [ "34A" ]    [ 1 ] 
       M12.2    113/495    4        [ "12B" ]    [ 3 ] 
       M22.2       8/33    4        [ "10A" ]    [ 4 ] 
       McL.2      1/135  134        [ "22A" ]    [ 1 ] 
        ON.2  61/109368 1792 [ "22A", "38A" ] [ 1, 1 ] 
       Suz.2      1/351  350        [ "28A" ]    [ 1 ]
    gap> for entry in sporautinfo2 do
    >      DisplayProbGenMaxesInfo( CharacterTable( entry[1] ), entry[4] );
    > od;
    F3+.2, 42E: 2^12.M24  (2)
                2^2.U6(2):S3x2  (1)
                2^(3+12).(L3(2)xS6)  (2)
                (S3xS3xG2(3)):2  (1)
                S6xL2(8):3  (1)
                7:6xS7  (1)
                7^(1+2)_+:(6xS3).2  (2)
    Fi22.2, 24G: Fi22.2M4  (1)
                 2^(5+8):(S3xS6)  (1)
                 3^5:(2xU4(2).2)  (1)
    HN.2, 24B: 2^(1+8)_+.(A5xA5).2^2  (1)
               5^2.5.5^2.4S5  (2)
               HN.2M13  (1)
    HS.2, 20C: (2xA6.2^2).2  (1)
               HS.2N5  (2)
               5:4xS5  (1)
    He.2, 24A: 2^(1+6)_+.L3(2).2  (1)
               S4xL3(2).2  (1)
    J2.2, 14A: L3(2).2x2  (1)
    J3.2, 34A: L2(17)x2  (1)
    M12.2, 12B: L2(11).2  (1)
                D8.(S4x2)  (1)
                3^(1+2):D8  (1)
    M22.2, 10A: M22.2M4  (1)
                A6.2^2  (1)
                L2(11).2  (2)
    McL.2, 22A: 2xM11  (1)
    ON.2, 22A: J1x2  (1)
    ON.2, 38A: J1x2  (1)
    Suz.2, 28A: (A4xL3(4):2_3):2  (1)
\end{verbatim}

In the other six cases,
we do not have the complete lists of primitive permutation characters,
so we choose a suitable element $\hat{s}$ for each group.
It is sufficient to prescribe $|\hat{s}|$, as follows.

\begin{verbatim}
    gap> sporautchoices:= [
    >        [ "Fi22",  "Fi22.2",  42 ],
    >        [ "Fi24'", "Fi24'.2", 46 ],
    >        [ "He",    "He.2",    42 ],
    >        [ "HN",    "HN.2",    44 ],
    >        [ "HS",    "HS.2",    30 ],
    >        [ "ON",    "ON.2",    38 ], ];;
\end{verbatim}

First we list the maximal subgroups of the corresponding simple groups
that contain the square of $\hat{s}$.

\begin{verbatim}
    gap> for triple in sporautchoices do
    >      tbl:= CharacterTable( triple[1] );
    >      tbl2:= CharacterTable( triple[2] );
    >      spos2:= PowerMap( tbl2, 2,
    >          Position( OrdersClassRepresentatives( tbl2 ), triple[3] ) );
    >      spos:= Position( GetFusionMap( tbl, tbl2 ), spos2 );
    >      DisplayProbGenMaxesInfo( tbl, AtlasClassNames( tbl ){ [ spos ] } );
    >    od;
    Fi22, 21A: O8+(2).3.2  (1)
               S3xU4(3).2_2  (1)
               A10.2  (1)
               A10.2  (1)
    F3+, 23A: Fi23  (1)
              F3+M7  (1)
    He, 21B: 3.A7.2  (1)
             7^(1+2):(S3x3)  (1)
             7:3xpsl(3,2)  (2)
    HN, 22A: 2.HS.2  (1)
    HS, 15A: A8.2  (1)
             5:4xa5  (1)
    ON, 19B: L3(7).2  (1)
             ONM2  (1)
             J1  (1)
\end{verbatim}

According to~\cite{CCN85}, exactly
the following maximal subgroups of the simple group $S$ in the above list
do {\bf not} extend to $\Aut(S)$:
The two $S_{10}$ type subgroups of $Fi_{22}$
and the two $L_3(7).2$ type subgroups of $O'N$.

Furthermore, the following maximal subgroups of $\Aut(S)$ with the property
that the intersection with $S$ is not maximal in $S$ have to be considered
whether they contain $s^{\prime}$:
$G_2(3).2$ and $3^5:(2 \times U_4(2).2)$ in $Fi_{22}.2$.
(Note that the order of the $7^{1+2}_+:(3 \times D_{16})$ type subgroup in
$O'N.2$ is obviously not divisible by $19$.)

\begin{verbatim}
    gap> 42 in OrdersClassRepresentatives( CharacterTable( "G2(3).2" ) );
    false
    gap> Size( CharacterTable( "U4(2)" ) ) mod 7 = 0;
    false
\end{verbatim}

So we take the extensions of the above maximal subgroups,
as described in~\cite{CCN85}.

\begin{verbatim}
    gap> SigmaFromMaxes( CharacterTable( "Fi22.2" ), "42A",
    >     [ CharacterTable( "O8+(2).3.2" ) * CharacterTable( "Cyclic", 2 ),
    >       CharacterTable( "S3" ) * CharacterTable( "U4(3).(2^2)_{122}" ) ],
    >     [ 1, 1 ] );
    163/1170
    gap> SigmaFromMaxes( CharacterTable( "Fi24'.2" ), "46A",
    >      [ CharacterTable( "Fi23" ) * CharacterTable( "Cyclic", 2 ),
    >        CharacterTable( "2^12.M24" ) ],
    >      [ 1, 1 ] );
    566/5481
    gap> SigmaFromMaxes( CharacterTable( "He.2" ), "42A",
    >      [ CharacterTable( "3.A7.2" ) * CharacterTable( "Cyclic", 2 ),
    >        CharacterTable( "7^(1+2):(S3x6)" ),
    >        CharacterTable( "7:6" ) * CharacterTable( "L3(2)" ) ],
    >      [ 1, 1, 1 ] );
    1/119
    gap> SigmaFromMaxes( CharacterTable( "HN.2" ), "44A",
    >      [ CharacterTable( "4.HS.2" ) ],
    >      [ 1 ] );
    997/192375
    gap> SigmaFromMaxes( CharacterTable( "HS.2" ), "30A",
    >      [ CharacterTable( "S8" ) * CharacterTable( "C2" ),
    >        CharacterTable( "5:4" ) * CharacterTable( "S5" ) ],
    >      [ 1, 1 ] );
    36/275
    gap> SigmaFromMaxes( CharacterTable( "ON.2" ), "38A",
    >      [ CharacterTable( "J1" ) * CharacterTable( "C2" ) ],
    >      [ 1 ] );
    61/109368
\end{verbatim}

\subsection{Other Simple Groups -- Easy Cases}\label{easyloop}

We are interested in simple groups $G$ for which \verb|ProbGenInfoSimple| does
not guarantee $\sprtotal(G) \geq 3$.
So we examine the remaining tables of simple groups in the {\GAP}
Character Table Library,
and distinguish the following three cases:
Either \verb|ProbGenInfoSimple| yields the lower bound at least three,
or a smaller bound,
or the computation of a lower bound fails because not enough information
is available to compute the primitive permutation characters.

\begin{verbatim}
    gap> names:= AllCharacterTableNames( IsSimple, true,
    >                                    IsDuplicateTable, false );;
    gap> names:= Difference( names, spornames );;
    gap> fails:= [];;
    gap> lessthan3:= [];;
    gap> atleast3:= [];;
    gap> for name in names do
    >      tbl:= CharacterTable( name );
    >      info:= ProbGenInfoSimple( tbl );
    >      if info = fail then
    >        Add( fails, name );
    >      elif info[3] < 3 then
    >        Add( lessthan3, info );
    >      else
    >        Add( atleast3, info );
    >      fi;
    >    od;
\end{verbatim}

For the following simple groups,
(currently) not enough information is available
in the {\GAP} Character Table Library and in the {\GAP} Library of Tables
of Marks,
for computing a lower bound for $\total(G)$.
Some of these groups will be dealt with in later sections,
and for the other groups, the bounds derived with theoretical arguments
in~\cite{BGK} are sufficient, so we need no {\GAP} computations for them.

\begin{verbatim}
    gap> fails;
    [ "2E6(2)", "2F4(8)", "A14", "A15", "A16", "A17", "A18", "E6(2)", "F4(2)", 
      "L4(4)", "L4(5)", "L4(9)", "L5(3)", "L6(2)", "L8(2)", "O10+(2)", "O10-(2)",
      "O10-(3)", "O7(5)", "O8-(3)", "O9(3)", "R(27)", "S10(2)", "S12(2)",
      "S4(7)", "S4(8)", "S4(9)", "S6(4)", "S6(5)", "S8(3)", "U4(4)", "U4(5)",
      "U5(3)", "U5(4)", "U7(2)" ]
\end{verbatim}

The following simple groups appear in~\cite[Table~1--6]{BGK}.
More detailed computations can be found in the sections~\ref{A5},
\ref{A6}, \ref{A7}, \ref{O8p2}, \ref{O8p3}, \ref{S62}, \ref{U42}, \ref{U43}.

\begin{verbatim}
    gap> PrintFormattedArray( lessthan3 );
          A5      1/3 2                [ "5A" ]       [ 1 ]
          A6      2/3 1                [ "5A" ]       [ 2 ]
          A7      2/5 2                [ "7A" ]       [ 2 ]
       O7(3)  199/351 1               [ "14A" ]       [ 3 ]
      O8+(2)  334/315 0 [ "15A", "15B", "15C" ] [ 7, 7, 7 ]
      O8+(3) 863/1820 2 [ "20A", "20B", "20C" ] [ 8, 8, 8 ]
       S6(2)      4/7 1                [ "9A" ]       [ 4 ]
       S8(2)     8/15 1               [ "17A" ]       [ 3 ]
       U4(2)    21/40 1               [ "12A" ]       [ 2 ]
       U4(3)   53/135 2                [ "7A" ]       [ 7 ]
\end{verbatim}

For the following simple groups $G$, the inequality $\total(G) < 1/3$
follows from the loop above.
The columns show the name of $G$, the values $\total(G)$ and $\sprtotal(G)$,
the class names of $s$ for which these values are attained,
and $|\M(G,s)|$.

(The entry for the group $L_7(2)$ would not fit on one screen line,
so we show it independently.)

\begin{verbatim}
    gap> PrintFormattedArray( Filtered( atleast3, l -> l[1] <> "L7(2)" ) );
      2F4(2)' 118/1755   14                           [ "16A" ]             [ 2 ]
       3D4(2)   1/5292 5291                           [ "13A" ]             [ 1 ]
          A10     3/10    3                           [ "21A" ]             [ 1 ]
          A11    2/105   52                           [ "11A" ]             [ 2 ]
          A12      2/9    4                           [ "35A" ]             [ 1 ]
          A13   4/1155  288                           [ "13A" ]             [ 5 ]
           A8     3/14    4                           [ "15A" ]             [ 1 ]
           A9     9/35    3                      [ "9A", "9B" ]          [ 4, 4 ]
        G2(3)      1/7    6                           [ "13A" ]             [ 3 ]
        G2(4)     1/21   20                           [ "13A" ]             [ 2 ]
        G2(5)     1/31   30                     [ "7A", "21A" ]         [ 10, 1 ]
      L2(101)    1/101  100                    [ "51A", "17A" ]          [ 1, 1 ]
      L2(103)  53/5253   99             [ "52A", "26A", "13A" ]       [ 1, 1, 1 ]
      L2(107)  55/5671  103 [ "54A", "27A", "18A", "9A", "6A" ] [ 1, 1, 1, 1, 1 ]
      L2(109)    1/109  108                    [ "55A", "11A" ]          [ 1, 1 ]
       L2(11)     7/55    7                            [ "6A" ]             [ 1 ]
      L2(113)    1/113  112                    [ "57A", "19A" ]          [ 1, 1 ]
      L2(121)    1/121  120                           [ "61A" ]             [ 1 ]
      L2(125)    1/125  124        [ "63A", "21A", "9A", "7A" ]    [ 1, 1, 1, 1 ]
       L2(13)     1/13   12                            [ "7A" ]             [ 1 ]
       L2(16)     1/15   14                           [ "17A" ]             [ 1 ]
       L2(17)     1/17   16                            [ "9A" ]             [ 1 ]
       L2(19)   11/171   15                           [ "10A" ]             [ 1 ]
       L2(23)   13/253   19                     [ "6A", "12A" ]          [ 1, 1 ]
       L2(25)     1/25   24                           [ "13A" ]             [ 1 ]
       L2(27)    5/117   23                     [ "7A", "14A" ]          [ 1, 1 ]
       L2(29)     1/29   28                           [ "15A" ]             [ 1 ]
       L2(31)   17/465   27                     [ "8A", "16A" ]          [ 1, 1 ]
       L2(32)     1/31   30              [ "3A", "11A", "33A" ]       [ 1, 1, 1 ]
       L2(37)     1/37   36                           [ "19A" ]             [ 1 ]
       L2(41)     1/41   40                     [ "21A", "7A" ]          [ 1, 1 ]
       L2(43)   23/903   39                    [ "22A", "11A" ]          [ 1, 1 ]
       L2(47)  25/1081   43        [ "24A", "12A", "8A", "6A" ]    [ 1, 1, 1, 1 ]
       L2(49)     1/49   48                           [ "25A" ]             [ 1 ]
       L2(53)     1/53   52                     [ "27A", "9A" ]          [ 1, 1 ]
       L2(59)  31/1711   55       [ "30A", "15A", "10A", "6A" ]    [ 1, 1, 1, 1 ]
       L2(61)     1/61   60                           [ "31A" ]             [ 1 ]
       L2(64)     1/63   62                    [ "65A", "13A" ]          [ 1, 1 ]
       L2(67)  35/2211   63                    [ "34A", "17A" ]          [ 1, 1 ]
       L2(71)  37/2485   67 [ "36A", "18A", "12A", "9A", "6A" ] [ 1, 1, 1, 1, 1 ]
       L2(73)     1/73   72                           [ "37A" ]             [ 1 ]
       L2(79)  41/3081   75       [ "40A", "20A", "10A", "8A" ]    [ 1, 1, 1, 1 ]
        L2(8)      1/7    6                      [ "3A", "9A" ]          [ 1, 1 ]
       L2(81)     1/81   80                           [ "41A" ]             [ 1 ]
       L2(83)  43/3403   79 [ "42A", "21A", "14A", "7A", "6A" ] [ 1, 1, 1, 1, 1 ]
       L2(89)     1/89   88              [ "45A", "15A", "9A" ]       [ 1, 1, 1 ]
       L2(97)     1/97   96                     [ "49A", "7A" ]          [ 1, 1 ]
       L3(11)   1/6655 6654                   [ "19A", "133A" ]          [ 1, 1 ]
        L3(2)      1/4    3                            [ "7A" ]             [ 1 ]
        L3(3)     1/24   23                           [ "13A" ]             [ 1 ]
        L3(4)      1/5    4                            [ "7A" ]             [ 3 ]
        L3(5)    1/250  249                           [ "31A" ]             [ 1 ]
        L3(7)   1/1372 1371                           [ "19A" ]             [ 1 ]
        L3(8)   1/1792 1791                           [ "73A" ]             [ 1 ]
        L3(9)   1/2880 2879                           [ "91A" ]             [ 1 ]
        L4(3)  53/1053   19                           [ "20A" ]             [ 1 ]
        L5(2)   1/5376 5375                           [ "31A" ]             [ 1 ]
       O8-(2)     1/63   62                           [ "17A" ]             [ 1 ]
        S4(4)     4/15    3                           [ "17A" ]             [ 2 ]
        S4(5)      1/5    4                           [ "13A" ]             [ 1 ]
        S6(3)    1/117  116                           [ "14A" ]             [ 2 ]
       Sz(32)   1/1271 1270                     [ "5A", "25A" ]          [ 1, 1 ]
        Sz(8)     1/91   90                            [ "5A" ]             [ 1 ]
       U3(11)   1/6655 6654                           [ "37A" ]             [ 1 ]
        U3(3)    16/63    3                     [ "6A", "12A" ]          [ 2, 2 ]
        U3(4)    1/160  159                           [ "13A" ]             [ 1 ]
        U3(5)   46/525   11                           [ "10A" ]             [ 2 ]
        U3(7)   1/1372 1371                           [ "43A" ]             [ 1 ]
        U3(8)   1/1792 1791                           [ "19A" ]             [ 1 ]
        U3(9)   1/3600 3599                           [ "73A" ]             [ 1 ]
        U5(2)     1/54   53                           [ "11A" ]             [ 1 ]
        U6(2)     5/21    4                           [ "11A" ]             [ 4 ]
    gap> First( atleast3, l -> l[1] = "L7(2)" );
    [ "L7(2)", 1/4388290560, 4388290559, [ "127A" ], [ 1 ] ]
\end{verbatim}

It should be mentioned that~\cite{BW1} states the following lower bounds
for the uniform spread of the groups $L_2(q)$.
\[
    \begin{array}{ll}
        q-2 & \mbox{\rm if $4 \leq q$ is even,} \\
        q-1 & \mbox{\rm if $11 \leq q \equiv 1 \pmod{4}$,} \\
        q-4 & \mbox{\rm if $11 \leq q \equiv -1 \pmod{4}$.}
    \end{array}
\]
These bounds appear in the third column of the above table.
Furthermore, \cite{BW1} states that the (uniform) spread
of alternating groups of even degree at least $8$ is exactly $4$.

For the sake of completeness, Table~\ref{maxtable} gives an overview of the
sets $\M(G,s)$ for those cases in the above list that are needed in~\cite{BGK}
but that do not require a further discussion here.
The structure of the maximal subgroups and the order of $s$ in the table
refer to the matrix groups not to the simple groups.
The number of the subgroups has been shown above,
the structure follows from~\cite{CCN85}.

\begin{table}
\caption{Maximal subgroups}\label{maxtable}
\[
   \begin{array}{|l|l|r|r|} \hline
        G           & \M(G,s) &  |s| & \mbox{\rm see~\cite{CCN85}} \\
      \hline\hline
        \SL(3,4) = 3.L_3(4)
                    & 3 \times L_3(2), 3 \times L_3(2), 3 \times L_3(2)
                            &   21 &  p.~23 \\
      \hline
        \Omega^-(8,2) = O^-_8(2)
                    & \Omega^-(4,4).2 = L_2(16).2
                            &   17 &  p.~89 \\
      \hline
        \Sp(4,4) = S_4(4)
                    & \Omega^-(4,4).2 = L_2(16).2, \Sp(2,16).2 = L_2(16).2
                            &   17 &  p.~44 \\
        \Sp(6,3) = 2.S_6(3)
                    & (4 \times U_3(3)).2, \Sp(2,17).3 = 2.L_2(27).3
                            &   28 & p.~113 \\
      \hline
        \SU(3,3) = U_3(3)
                    & 3^{1+2}_+:8, \GU(2,3) = 4.S_4
                            &    6 &  p.~14 \\
        \SU(3,5) = 3.U_3(5)
                    & 3 \times 5^{1+2}_+:8, \GU(2,5) = 3 \times 2S_5
                            &   30 &  p.~34 \\
        \SU(5,2) = U_5(2)
                    & L_2(11)
                            &   11 &  p.~73 \\
      \hline
   \end{array}
\]
\end{table}

\subsection{Automorphism Groups of other Simple Groups -- Easy Cases}%
\label{easyloopaut}

We deal with automorphic extensions of those simple groups that are
listed in Table~\ref{thetable} and that have been treated successfully in
Section~\ref{easyloop}.

For the following groups, \verb|ProbGenInfoAlmostSimple| can be used
because {\GAP} can compute their primitive permutation characters.

\begin{verbatim}
    gap> list:= [
    >   [ "A5", "A5.2" ],
    >   [ "A6", "A6.2_1" ],
    >   [ "A6", "A6.2_2" ],
    >   [ "A6", "A6.2_3" ],
    >   [ "A7", "A7.2" ],
    >   [ "A8", "A8.2" ],
    >   [ "A9", "A9.2" ],
    >   [ "A11", "A11.2" ],
    >   [ "L3(2)", "L3(2).2" ],
    >   [ "L3(3)", "L3(3).2" ],
    >   [ "L3(4)", "L3(4).2_1" ],
    >   [ "L3(4)", "L3(4).2_2" ],
    >   [ "L3(4)", "L3(4).2_3" ],
    >   [ "L3(4)", "L3(4).3" ],
    >   [ "S4(4)", "S4(4).2" ],
    >   [ "U3(3)", "U3(3).2" ],
    >   [ "U3(5)", "U3(5).2" ],
    >   [ "U3(5)", "U3(5).3" ],
    >   [ "U4(2)", "U4(2).2" ],
    >   [ "U4(3)", "U4(3).2_1" ],
    >   [ "U4(3)", "U4(3).2_3" ],
    > ];;
    gap> autinfo:= [];;
    gap> fails:= [];;
    gap> for pair in list do
    >      tbl:= CharacterTable( pair[1] );
    >      tblG:= CharacterTable( pair[2] );
    >      info:= ProbGenInfoSimple( tbl );
    >      spos:= List( info[4], x -> Position( AtlasClassNames( tbl ), x ) );
    >      Add( autinfo, ProbGenInfoAlmostSimple( tbl, tblG, spos ) );
    >    od;
    gap> PrintFormattedArray( autinfo );
           A5.2      0        [ "5AB" ]    [ 1 ]
         A6.2_1    2/3        [ "5AB" ]    [ 2 ]
         A6.2_2    1/6         [ "5A" ]    [ 1 ]
         A6.2_3      0        [ "5AB" ]    [ 1 ]
           A7.2   1/15        [ "7AB" ]    [ 1 ]
           A8.2  13/28       [ "15AB" ]    [ 1 ]
           A9.2    1/4        [ "9AB" ]    [ 1 ]
          A11.2  1/945       [ "11AB" ]    [ 1 ]
        L3(2).2    1/4        [ "7AB" ]    [ 1 ]
        L3(3).2   1/18       [ "13AB" ]    [ 1 ]
      L3(4).2_1   3/10        [ "7AB" ]    [ 3 ]
      L3(4).2_2  11/60         [ "7A" ]    [ 1 ]
      L3(4).2_3   1/12        [ "7AB" ]    [ 1 ]
        L3(4).3   1/64         [ "7A" ]    [ 1 ]
        S4(4).2      0       [ "17AB" ]    [ 2 ]
        U3(3).2    2/7 [ "6A", "12AB" ] [ 2, 2 ]
        U3(5).2   2/21        [ "10A" ]    [ 2 ]
        U3(5).3 46/525        [ "10A" ]    [ 2 ]
        U4(2).2  16/45       [ "12AB" ]    [ 2 ]
      U4(3).2_1 76/135         [ "7A" ]    [ 3 ]
      U4(3).2_3 31/162        [ "7AB" ]    [ 3 ]
\end{verbatim}

We see that from this list,
the two groups $A_6.2_1 = S_6$ and $U_4(3).2_1$ require further computations
(see Sections~\ref{A6} and~\ref{U43}, respectively)
because the bound in the second column is larger than $1/2$.

Also $U_4(2)$ is not done by the above,
because in~\cite[Table~4]{BGK}, an element $s$ of order $9$ is chosen
for the simple group, see Section~\ref{U42}.

Finally, we deal with automorphic extensions of the groups $L_4(3)$,
$O_8^-(2)$, $S_6(3)$, and $U_5(2)$.

For $S = L_4(3)$ and $s \in S$ of order $20$,
we have $\M(S,s) = \{ (4 \times A_6):2 \}$,
the subgroup has index $2\,106$, see~\cite[p.~69]{CCN85}.

\begin{verbatim}
    gap> t:= CharacterTable( "L4(3)" );;
    gap> prim:= PrimitivePermutationCharacters( t );;
    gap> spos:= Position( AtlasClassNames( t ), "20A" );;
    gap> prim:= Filtered( prim, x -> x[ spos ] <> 0 );
    [ Character( CharacterTable( "L4(3)" ), [ 2106, 106, 42, 0, 27, 27, 0, 46, 6, 
          6, 1, 7, 7, 0, 3, 3, 0, 0, 0, 1, 1, 1, 0, 0, 0, 0, 0, 1, 1 ] ) ]
\end{verbatim}

For the three automorphic extensions of the structure $G = S.2$,
we compute the extensions of the permutation character,
and the bounds $\total^{\prime}(G,s)$.

\begin{verbatim}
    gap> for name in [ "L4(3).2_1", "L4(3).2_2", "L4(3).2_3" ] do
    >      t2:= CharacterTable( name );
    >      map:= InverseMap( GetFusionMap( t, t2 ) );
    >      torso:= List( prim, pi -> CompositionMaps( pi, map ) );
    >      ext:= Concatenation( List( torso,
    >                              x -> PermChars( t2, rec( torso:= x ) ) ) );
    >      sigma:= ApproxP( ext, Position( OrdersClassRepresentatives( t2 ), 20 ) );
    >      max:= Maximum( sigma{ Difference( PositionsProperty(
    >                           OrdersClassRepresentatives( t2 ), IsPrimeInt ),
    >                           ClassPositionsOfDerivedSubgroup( t2 ) ) } );
    >      Print( name, ":\n", ext, "\n", max, "\n" );
    > od;
    L4(3).2_1:
    [ Character( CharacterTable( "L4(3).2_1" ), [ 2106, 106, 42, 0, 27, 0, 46, 6, 
          6, 1, 7, 0, 3, 0, 0, 1, 1, 0, 0, 0, 0, 0, 1, 1, 0, 4, 0, 0, 6, 6, 6, 6, 
          2, 0, 0, 0, 0, 0, 0, 1, 1, 1, 1 ] ) ]
    0
    L4(3).2_2:
    [ Character( CharacterTable( "L4(3).2_2" ), 
        [ 2106, 106, 42, 0, 27, 27, 0, 46, 6, 6, 1, 7, 7, 0, 3, 3, 0, 0, 0, 1, 1, 
          1, 0, 0, 0, 1, 306, 306, 42, 6, 10, 10, 0, 0, 15, 15, 3, 3, 3, 3, 0, 0, 
          1, 1, 0, 1, 1, 0, 0 ] ) ]
    17/117
    L4(3).2_3:
    [ Character( CharacterTable( "L4(3).2_3" ), [ 2106, 106, 42, 0, 27, 0, 46, 6, 
          6, 1, 7, 0, 3, 0, 0, 1, 1, 0, 0, 0, 1, 36, 0, 0, 6, 6, 2, 2, 2, 1, 1, 
          0, 0, 0 ] ) ]
    2/117
\end{verbatim}

For $S = O_8^-(2)$ and $s \in S$ of order $17$,
we have $\M(S,s) = \{ L_2(16).2 \}$,
the subgroup extends to $L_2(16).4$ in $S.2$, see~\cite[p.~89]{CCN85}.
This is a non-split extension, so $\total^{\prime}(S.2,s) = 0$ holds.

\begin{verbatim}
    gap> SigmaFromMaxes( CharacterTable( "O8-(2).2" ), "17AB",
    >        [ CharacterTable( "L2(16).4" ) ], [ 1 ], "outer" );
    0
\end{verbatim}

For $S = S_6(3)$ and $s \in S$ irreducible of order $14$,
we have $\M(S,s) = \{ (2 \times U_3(3)).2, L_2(27).3 \}$.
In $G = S.2$, the subgroups extend to $(4 \times U_3(3)).2$ and $L_2(27).6$,
respectively, see~\cite[p.~113]{CCN85}.
In order to show that $\total^{\prime}(G,s) = 7/3240$ holds,
we compute the primitive permutation characters of $S$
(cf.~Section~\ref{easyloop}) and the unique extensions to $G$
of those which are nonzero on $s$.

\begin{verbatim}
    gap> t:= CharacterTable( "S6(3)" );;
    gap> t2:= CharacterTable( "S6(3).2" );;
    gap> prim:= PrimitivePermutationCharacters( t );;
    gap> spos:= Position( AtlasClassNames( t ), "14A" );;
    gap> prim:= Filtered( prim, x -> x[ spos ] <> 0 );;
    gap> map:= InverseMap( GetFusionMap( t, t2 ) );;
    gap> torso:= List( prim, pi -> CompositionMaps( pi, map ) );;
    gap> ext:= List( torso, pi -> PermChars( t2, rec( torso:= pi ) ) );
    [ [ Character( CharacterTable( "S6(3).2" ), [ 155520, 0, 288, 0, 0, 0, 216, 
              54, 0, 0, 0, 0, 0, 0, 0, 0, 0, 0, 0, 0, 0, 6, 1, 0, 0, 0, 0, 0, 6, 
              0, 0, 0, 0, 0, 0, 0, 0, 0, 0, 1, 1, 1, 0, 0, 0, 0, 0, 0, 0, 0, 144, 
              288, 0, 0, 0, 6, 0, 0, 0, 0, 0, 0, 0, 0, 6, 0, 0, 0, 0, 0, 1, 1, 1, 
              1, 0, 0 ] ) ],
      [ Character( CharacterTable( "S6(3).2" ), [ 189540, 1620, 568, 0, 486, 0, 
              0, 27, 540, 84, 24, 0, 0, 0, 0, 0, 54, 0, 0, 10, 0, 7, 1, 6, 6, 0, 
              0, 0, 0, 0, 0, 18, 0, 0, 0, 6, 12, 0, 0, 0, 0, 1, 0, 0, 0, 0, 0, 0, 
              0, 0, 234, 64, 30, 8, 0, 3, 90, 6, 0, 4, 10, 6, 0, 2, 1, 0, 0, 0, 
              0, 0, 0, 0, 1, 1, 0, 0 ] ) ] ]
    gap> spos:= Position( AtlasClassNames( t2 ), "14A" );;
    gap> sigma:= ApproxP( Concatenation( ext ), spos );;
    gap> Maximum( sigma{ Difference(
    >      PositionsProperty( OrdersClassRepresentatives( t2 ), IsPrimeInt ),
    >      ClassPositionsOfDerivedSubgroup( t2 ) ) } );
    7/3240
\end{verbatim}

For $S = U_5(2)$ and $s \in S$ of order $11$,
we have $\M(S,s) = \{ L_2(11) \}$,
the subgroup extends to $L_2(11).2$ in $S.2$, see~\cite[p.~73]{CCN85}.

\begin{verbatim}
    gap> SigmaFromMaxes( CharacterTable( "U5(2).2" ), "11AB",
    >        [ CharacterTable( "L2(11).2" ) ], [ 1 ], "outer" );
    1/288
\end{verbatim}

Here we clean the workspace for the first time.
This may save more than $100$ megabytes, due to the fact that the caches
for tables of marks and character tables are flushed.

\begin{verbatim}
    gap> CleanWorkspace();
\end{verbatim}

\subsection{$O_8^-(3)$}\label{O8m3}

We show that $S = O_8^-(3) = \Omega^-(8,3)$ satisfies the following.
\begin{enumerate}
\item[(a)]
    For $s \in S$ of order $41$,
    $\M(S,s)$ consists of one group of the type
    $L_2(81).2_1 = \Omega^-(4,9).2$.
\item[(b)]
    $\total(S,s) = 1/567$.
\end{enumerate}

The only maximal subgroups of $S$ containing elements of order $41$
have the type $L_2(81).2_1$,
and there is one class of these subgroups, see~\cite[p.~141]{CCN85}.

\begin{verbatim}
    gap> SigmaFromMaxes( CharacterTable( "O8-(3)" ), "41A",
    >    [ CharacterTable( "L2(81).2_1" ) ], [ 1 ] );
    1/567
\end{verbatim}

\subsection{$O_{10}^+(2)$}\label{O10p2}

We show that $S = O_{10}^+(2) = \Omega^+(10,2)$ satisfies the following.
\begin{enumerate}
\item[(a)]
    For $s \in S$ of order $45$,
    $\M(S,s)$ consists of one group of the type
    $(A_5 \times U_4(2)).2 = (\Omega^-(4,2) \times \Omega^-(6,2)).2$.
\item[(b)]
    $\total(S,s) = 43/4\,216$.
\item[(c)]
    For $s$ as in (a),
    the maximal subgroup in (a) extends to $S_5 \times U_4(2).2$
    in $G = \Aut(S) = S.2$,
    and $\total^{\prime}(G,s) = 23/248$.
\end{enumerate}

The only maximal subgroups of $S$ containing elements of order $45$
are one class of groups
$H = (A_5 \times U_4(2)):2$, see~\cite[p.~146]{CCN85}.
(Note that none of the groups $S_8(2)$, $O_8^+(2)$, $L_5(2)$, $O_8^-(2)$,
and $A_8$ contains elements of order $45$.)
$H$ extends to subgroups of the type $H.2 = S_5 \times U_4(2):2$ in $G$,
so we can compute $1_H^S = (1_{H.2}^G)_S$.

\begin{verbatim}
    gap> ForAny( [ "S8(2)", "O8+(2)", "L5(2)", "O8-(2)", "A8" ],
    >            x -> 45 in OrdersClassRepresentatives( CharacterTable( x ) ) );
    false
    gap> t:= CharacterTable( "O10+(2)" );;
    gap> t2:= CharacterTable( "O10+(2).2" );;
    gap> s2:= CharacterTable( "A5.2" ) * CharacterTable( "U4(2).2" );
    CharacterTable( "A5.2xU4(2).2" )
    gap> pi:= PossiblePermutationCharacters( s2, t2 );;
    gap> spos:= Position( OrdersClassRepresentatives( t2 ), 45 );;
    gap> approx:= ApproxP( pi, spos );;
    gap> Maximum( approx{ ClassPositionsOfDerivedSubgroup( t2 ) } );
    43/4216
\end{verbatim}

Statement~(c) follows from considering the outer classes
of prime element order.

\begin{verbatim}
    gap> Maximum( approx{ Difference(
    >      PositionsProperty( OrdersClassRepresentatives( t2 ), IsPrimeInt ),
    >      ClassPositionsOfDerivedSubgroup( t2 ) ) } );
    23/248
\end{verbatim}

Alternatively, we can use \verb|SigmaFromMaxes|.

\begin{verbatim}
    gap> SigmaFromMaxes( t2, "45AB", [ s2 ], [ 1 ], "outer" );
    23/248
\end{verbatim}

\subsection{$O_{10}^-(2)$}\label{O10m2}

We show that $S = O_{10}^-(2) = \Omega^-(10,2)$ satisfies the following.
\begin{enumerate}
\item[(a)]
    For $s \in S$ of order $33$,
    $\M(S,s)$ consists of one group of the type
    $3 \times U_5(2) = \GU(5,2)$.
\item[(b)]
    $\total(S,s) = 1/119$.
\item[(c)]
    For $s$ as in (a),
    the maximal subgroup in (a) extends to $(3 \times U_5(2)).2$
    in $G$,
    and $\total^{\prime}(G,s) = 1/595$.
\end{enumerate}

The only maximal subgroups of $S$ containing elements of order $11$
have the types $A_{12}$ and $3 \times U_5(2)$, see~\cite[p.~147]{CCN85}.
So $3 \times U_5(2)$ is the unique class of subgroups containing elements
of order $33$.
This shows statement~(a),
and statement~(b) follows using \verb|SigmaFromMaxes|.

\begin{verbatim}
    gap> SigmaFromMaxes( CharacterTable( "O10-(2)" ), "33A",
    >    [ CharacterTable( "Cyclic", 3 ) * CharacterTable( "U5(2)" ) ], [ 1 ] );
    1/119
\end{verbatim}

The structure of the maximal subgroup of $G$ follows
from~\cite[p.~147]{CCN85}.
We create its character table with a generic construction
that is based on the fact that the outer automorphism acts nontrivially on
the two direct factors; this determines the character table uniquely.
(See~\cite{Auto} for details.)

\begin{verbatim}
    gap> tblG:= CharacterTable( "U5(2)" );;
    gap> tblMG:= CharacterTable( "Cyclic", 3 ) * tblG;;
    gap> tblGA:= CharacterTable( "U5(2).2" );;
    gap> acts:= PossibleActionsForTypeMGA( tblMG, tblG, tblGA );;
    gap> poss:= Concatenation( List( acts, pi ->
    >            PossibleCharacterTablesOfTypeMGA( tblMG, tblG, tblGA, pi,
    >                "(3xU5(2)).2" ) ) );
    [ rec( MGfusMGA := [ 1, 2, 3, 4, 4, 5, 5, 6, 7, 8, 9, 10, 11, 12, 12, 13, 13,
              14, 14, 15, 15, 16, 17, 17, 18, 18, 19, 20, 21, 21, 22, 22, 23, 23,
              24, 24, 25, 25, 26, 27, 27, 28, 28, 29, 29, 30, 30, 31, 32, 33, 34,
              35, 36, 37, 38, 39, 40, 41, 42, 43, 44, 45, 46, 47, 48, 49, 50, 51,
              52, 53, 54, 55, 56, 57, 58, 59, 60, 61, 62, 63, 64, 65, 66, 67, 68,
              69, 70, 71, 72, 73, 74, 75, 76, 77, 31, 32, 33, 35, 34, 37, 36, 38,
              39, 40, 41, 42, 43, 45, 44, 47, 46, 49, 48, 51, 50, 52, 54, 53, 56,
              55, 57, 58, 60, 59, 62, 61, 64, 63, 66, 65, 68, 67, 69, 71, 70, 73,
              72, 75, 74, 77, 76 ], table := CharacterTable( "(3xU5(2)).2" ) ) ]
\end{verbatim}

Now statement~(c) follows using \verb|SigmaFromMaxes|.

\begin{verbatim}
    gap> SigmaFromMaxes( CharacterTable( "O10-(2).2" ), "33AB",
    >        [ poss[1].table ], [ 1 ], "outer" );
    1/595
\end{verbatim}

\subsection{$O_{12}^+(2)$}\label{O12p2}

We show that $S = O_{12}^+(2) = \Omega^+(12,2)$ satisfies the following.
\begin{enumerate}
\item[(a)]
    For $s \in S$ of the type $4^- \perp 8^-$
    (i.~e., $s$ decomposes the natural $12$-dimensional module for
    $\GO^+_{12}(2) = S.2$ into an orthogonal sum of two irreducible
    modules of the dimensions $4$ and $8$, respectively) and of order $85$,
    $\M(S,s)$ consists of one group of the type
    $G_8 = (\Omega^-(4,2) \times \Omega^-(8,2)).2$ and
    two groups of the type $L_4(4).2^2 = \Omega^+(6,4).2^2$
    that are conjugate in $G = \Aut(S) = S.2 = \SO^+(12,2)$
    but \emph{not} conjugate in $S$.
\item[(b)]
    $\total(S,s) = 7\,675/1\,031\,184$.
\item[(c)]
    $\total^{\prime}(G,s) = 73/1\,008$.
\end{enumerate}

The element $s$ is a ppd$(12,2;8)$-element in the sense of~\cite{GPPS},
so the maximal subgroups of $S$ that contain $s$ are among the nine
cases (2.1)--(2.9) listed in this paper;
in the notation of this paper,
we have $q = 2$, $d = 12$, $e = 8$, and $r = 17$.
Case~(2.1) does not occur for orthogonal groups and $q = 2$,
according to~\cite{KlL90};
case~(2.2) contributes a unique maximal subgroup,
the stabilizer $G_8$ of the orthogonal decomposition;
the cases~(2.3), (2.4)~(a), (2.5), and (2.6)~(a) do not occur
because $r \not= e+1$ in our situation;
case~(2.4)~(b) describes extension field type subgroups that are contained
in $\GammaL(6,4)$,
which yields the candidates $\GU(6,2).2 \cong 3.U_6(2).S_3$
--but $3.U_6(2).3$ does not contain elements of order $85$--
and $\Omega^+(6,4).2^2 \cong L_4(4).2^2$
(two classes by~\cite[Prop.~4.3.14]{KlL90});
the cases~(2.6)~(b)--(c) and~(2.8) do not occur because they require
$d \leq 8$;
case~(2.7) does not occur because~\cite[Table~5]{GPPS} contains no entry for
$r = 2e+1 = 17$;
finally, case~(2.9) does not occur because it requires $e \in \{ d-1, d \}$
in the case $r = 2e+1$.

So we need the permutation characters of the actions on the cosets of
$L_4(4).2^2$ (two classes) and $G_8$.
According to~\cite[Prop.~4.1.6]{KlL90},
$G_8$ has the structure $(\Omega^-(4,2) \times \Omega^-(8,2)).2$.

Currently the {\GAP} Character Table Library does not contain the
character table of $S$, but the table of $G$ is available,
and we work with this table.

The two classes of $L_4(4).2^2$ type subgroups in $S$ are fused in $G$.
This can be seen from the fact that inducing the trivial character of
a subgroup $H_1 = L_4(4).2^2$ of $S$ to $G$ yields a character $\psi$
whose values are not all even;
note that if $H_1$ would extend in $G$ to a subgroup of twice the size
of $H_1$ then $\psi$ would be induced from a degree two character
of this subgroup whose values are all even,
and induction preserves this property.

\begin{verbatim}
    gap> CharacterTable( "O12+(2)" );
    fail
    gap> t:= CharacterTable( "O12+(2).2" );;
    gap> h1:= CharacterTable( "L4(4).2^2" );;
    gap> psi:= PossiblePermutationCharacters( h1, t );;
    gap> Length( psi );
    1
    gap> ForAny( psi[1], IsOddInt );
    true
\end{verbatim}

The fixed element $s$ of order $85$ is contained in a member of each
of the two conjugacy classes of the type $L_4(4).2^2$ in $S$,
since $S$ contains only one class of subgroups of the order $85$;
note that the order of the Sylow $17$ centralizer (in both $S$ and $G$)
is not divisible by $25$.

\begin{verbatim}
    gap> SizesCentralizers( t ){ PositionsProperty(
    >        OrdersClassRepresentatives( t ), x -> x = 17 ) } / 25;
    [ 408/5, 408/5 ]
\end{verbatim}

This implies that the restriction of $\psi$ to $S$
is the sum $\psi_S = \pi_1 + \pi_2$, say,
of the first two interesting permutation characters of $S$.

The subgroup $G_8$ of $S$ extends to a group of the structure
$H_2 = \Omega^-(4,2).2 \times \Omega^-(8,2).2$ in $G$,
inducing the trivial characters of $H_2$ to $G$ yields
a permutation character $\varphi$ of $G$ whose restriction to $S$
is the third permutation character $\varphi_S = \pi_3$, say.

\begin{verbatim}
    gap> h2:= CharacterTable( "S5" ) * CharacterTable( "O8-(2).2" );;
    gap> phi:= PossiblePermutationCharacters( h2, t );;
    gap> Length( phi );
    1
\end{verbatim}

We have $\pi_1(1) = \pi_2(1)$ and $\pi_1(s) = \pi_2(s)$,
the latter again because $S$ contains only one class of subgroups
of order $85$.

Now statement~(a) follows from the fact that $\pi_i(s) = 1$ holds
for $1 \leq i \leq 3$.

\begin{verbatim}
    gap> prim:= Concatenation( psi, phi );;
    gap> spos:= Position( OrdersClassRepresentatives( t ), 85 );
    213
    gap> List( prim, x -> x[ spos ] );
    [ 2, 1 ]
\end{verbatim}

For statement~(b), we compute $\total(S,s)$.
Note that we have to consider only classes inside $S = G^{\prime}$,
and that
\[
   \total( g, s )
     = \sum_{i=1}^3 \frac{\pi_i(s) \cdot \pi_i(g)}{\pi_i(1)}
     = \frac{\psi(s) \cdot \psi(g)}{\psi(1)}
       + \frac{\varphi(s) \cdot \varphi(g)}{\varphi(1)}
\]
holds for $g \in S^{\times}$,
so the characters $\psi$ and $\varphi$ are sufficient.

\begin{verbatim}
    gap> approx:= ApproxP( prim, spos );;
    gap> Maximum( approx{ ClassPositionsOfDerivedSubgroup( t ) } );
    7675/1031184
\end{verbatim}

Statement~(c) follows from considering the outer involution classes.
Note that by~\cite[Remark after Proposition~5.14]{BGK},
only the subgroup $H_2$ need to be considered,
no novelties appear.

\begin{verbatim}
    gap> Maximum( approx{ Difference(
    >      PositionsProperty( OrdersClassRepresentatives( t ), IsPrimeInt ),
    >      ClassPositionsOfDerivedSubgroup( t ) ) } );
    73/1008
\end{verbatim}

\subsection{$O_{12}^-(2)$}\label{O12m2}

We show that $S = O_{12}^-(2) = \Omega^-(12,2)$ satisfies the following.
\begin{enumerate}
\item[(a)]
    For $s \in S$ irreducible of order $2^6+1 = 65$,
    $\M(S,s)$ consists of two groups of the types
    $U_4(4).2 = \Omega^-(6,4).2$ and $L_2(64).3 = \Omega^-(4,8).3$,
    respectively.
\item[(b)]
    $\total(S,s) = 1/1\,023$.
\item[(c)]
    $\total^{\prime}(\Aut(S),s) = 1/347\,820$.
\end{enumerate}

By~\cite{Be00}, $\M(S,s)$ consists of extension field subgroups,
which have the structures $U_4(4).2$ and $L_2(64).3$, respectively,
and by~\cite[Prop.~4.3.16]{KlL90},
there is just one class of each of these types.

Currently the {\GAP} Character Table Library does not contain the
character table of $S$, but the table of $G = \Aut(S) = O_{12}^-(2).2$
is available.
So we compute the permutation characters $\pi_1, \pi_2$ of the extensions
of the groups in $\M(S,s)$ to $G$
--these maximal subgroups have the structures $U_4(4).4$ and $L_2(64).6$,
respectively--
and compute the fixed point ratios of the restrictions to $S$.

\begin{verbatim}
    gap> t:= CharacterTable( "O12-(2)" );
    fail
    gap> t:= CharacterTable( "O12-(2).2" );;
    gap> s1:= CharacterTable( "U4(4).4" );;
    gap> pi1:= PossiblePermutationCharacters( s1, t );;
    gap> s2:= CharacterTable( "L2(64).6" );;
    gap> pi2:= PossiblePermutationCharacters( s2, t );;
    gap> prim:= Concatenation( pi1, pi2 );;  Length( prim );
    2
\end{verbatim}

Now statement~(a) follows from the fact that $\pi_1(s) = \pi_2(s) = 1$ holds.

\begin{verbatim}
    gap> spos:= Position( OrdersClassRepresentatives( t ), 65 );;
    gap> List( prim, x -> x[ spos ] );
    [ 1, 1 ]
\end{verbatim}

For statement~(b), we compute $\total(S,s)$;
note that we have to consider only classes inside $S = G^{\prime}$.

\begin{verbatim}
    gap> approx:= ApproxP( prim, spos );;
    gap> Maximum( approx{ ClassPositionsOfDerivedSubgroup( t ) } );
    1/1023
\end{verbatim}

Statement~(c) follows from the values on the outer involution classes.

\begin{verbatim}
    gap> Maximum( approx{ Difference(
    >      PositionsProperty( OrdersClassRepresentatives( t ), IsPrimeInt ),
    >      ClassPositionsOfDerivedSubgroup( t ) ) } );
    1/347820
\end{verbatim}

\subsection{$S_6(4)$}\label{S64}

We show that $S = S_6(4) = \Sp(6,4)$ satisfies the following.
\begin{enumerate}
\item[(a)]
    For $s \in S$ irreducible of order $65$,
    $\M(S,s)$ consists of two groups of the types
    $U_4(4).2 = \Omega^-(6,4).2$ and $L_2(64).3 = \Sp(2,64).3$,
    respectively.
\item[(b)]
    $\total(S,s) = 16/63$.
\item[(c)]
    $\total^{\prime}(\Aut(S),s) = 0$.
\end{enumerate}

By~\cite{Be00}, the element $s$ is contained in maximal
subgroups of the given types,
and by~\cite[Prop.~4.3.10, 4.8.6]{KlL90},
there is exactly one class of these subgroups.

The character tables of these two subgroups are currently not contained
in the {\GAP} Character Table Library.
We compute the permutation character induced from the first subgroup
as the  unique character of the right degree that is combinatorially
possible (cf.~\cite{BP98copy}).

\begin{verbatim}
    gap> t:= CharacterTable( "S6(4)" );;
    gap> degree:= Size( t ) / ( 2 * Size( CharacterTable( "U4(4)" ) ) );;
    gap> pi1:= PermChars( t, rec( torso:= [ degree ] ) );;
    gap> Length( pi1 );
    1
\end{verbatim}

The index of the second subgroup is too large for this simpleminded
approach;
therefore, we first restrict the set of possible irreducible constituents
of the permutation character to those of $1_H^G$,
where $H$ is the derived subgroup of $L_2(64).3$,
for which the character table is available.

\begin{verbatim}
    gap> CharacterTable( "L2(64).3" );  CharacterTable( "U4(4).2" );
    fail
    fail
    gap> s:= CharacterTable( "L2(64)" );;
    gap> subpi:= PossiblePermutationCharacters( s, t );;
    gap> Length( subpi );
    1
    gap> scp:= MatScalarProducts( t, Irr( t ), subpi );;
    gap> nonzero:= PositionsProperty( scp[1], x -> x <> 0 );
    [ 1, 11, 13, 14, 17, 18, 32, 33, 56, 58, 59, 73, 74, 77, 78, 79, 80, 93, 95, 
      96, 103, 116, 117, 119, 120 ]
    gap> const:= RationalizedMat( Irr( t ){ nonzero } );;
    gap> degree:= Size( t ) / ( 3 * Size( s ) );
    5222400
    gap> pi2:= PermChars( t, rec( torso:= [ degree ], chars:= const ) );;
    gap> Length( pi2 );
    1
    gap> prim:= Concatenation( pi1, pi2 );;
\end{verbatim}

Now statement~(a) follows from the fact that $\pi_1(s) = \pi_2(s) = 1$ holds.

\begin{verbatim}
    gap> spos:= Position( OrdersClassRepresentatives( t ), 65 );;
    gap> List( prim, x -> x[ spos ] );
    [ 1, 1 ]
\end{verbatim}

For statement~(b), we compute $\total(G,s)$.

\begin{verbatim}
    gap> Maximum( ApproxP( prim, spos ) );
    16/63
\end{verbatim}

In order to prove statement~(c), we have to consider only the extensions
of the above permutation characters of $S$ to $\Aut(S) \cong S.2$
(cf.~\cite[Section~2.2]{BGK}).

\begin{verbatim}
    gap> t2:= CharacterTable( "S6(4).2" );;
    gap> tfust2:= GetFusionMap( t, t2 );;
    gap> cand:= List( prim, x -> CompositionMaps( x, InverseMap( tfust2 ) ) );;
    gap> ext:= List( cand, pi -> PermChars( t2, rec( torso:= pi ) ) );
    [ [ Character( CharacterTable( "S6(4).2" ), [ 2016, 512, 96, 128, 32, 120, 0, 
              6, 16, 40, 24, 0, 8, 136, 1, 6, 6, 1, 32, 0, 8, 6, 2, 0, 2, 0, 0, 
              4, 0, 16, 32, 1, 8, 2, 6, 2, 1, 2, 4, 0, 0, 1, 6, 0, 1, 10, 0, 1, 
              1, 0, 10, 10, 4, 0, 1, 0, 2, 0, 2, 1, 2, 2, 1, 1, 0, 0, 0, 1, 1, 1, 
              1, 0, 0, 0, 0, 0, 32, 0, 0, 8, 0, 0, 0, 0, 8, 8, 0, 0, 0, 0, 8, 0, 
              0, 0, 2, 2, 0, 2, 2, 0, 2, 2, 2, 0, 0 ] ) ], 
      [ Character( CharacterTable( "S6(4).2" ), [ 5222400, 0, 0, 0, 1280, 0, 960, 
              120, 0, 0, 0, 0, 0, 0, 1600, 0, 0, 0, 0, 0, 0, 0, 0, 0, 8, 1, 0, 0, 
              15, 0, 0, 0, 0, 0, 0, 0, 0, 0, 0, 0, 0, 1, 0, 0, 0, 0, 0, 0, 10, 0, 
              0, 0, 0, 0, 0, 1, 0, 0, 0, 0, 0, 0, 0, 0, 1, 1, 1, 1, 1, 1, 1, 0, 
              0, 0, 0, 960, 0, 0, 0, 16, 0, 24, 12, 0, 0, 0, 0, 0, 0, 0, 0, 0, 0, 
              0, 0, 4, 1, 0, 0, 3, 0, 0, 0, 0, 0 ] ) ] ]
    gap> spos2:= Position( OrdersClassRepresentatives( t2 ), 65 );;
    gap> sigma:= ApproxP( Concatenation( ext ), spos2 );;
    gap> Maximum( approx{ Difference(
    >      PositionsProperty( OrdersClassRepresentatives( t2 ), IsPrimeInt ),
    >      ClassPositionsOfDerivedSubgroup( t2 ) ) } );
    0
\end{verbatim}

For the simple group, we can \emph{alternatively} consider a reducible
element $s: 2 \perp 4$ of order $85$,
which is a multiple of the primitive prime divisor $r = 17$ of $4^4-1$.
So we have $e = 4$, $d = 6$, and $q = 4$, in the terminology of~\cite{GPPS}.
Then $\M(S,s)$ consists of two groups, of the types
$\Omega^+(6,4).2 \cong L_4(4).2_2$ and $\Sp(2,4) \times \Sp(4,4)$.
This can be shown by checking~\cite[Ex.~2.1--2.9]{GPPS}.
Ex.~2.1 yields the candidates $\Omega^{\pm}(6,4).2$,
but only $\Omega^+(6,4).2$ contains elements of order $85$.
Ex.~2.2 yields the stabilizer of a two-dimensional subspace,
which has the structure $\Sp(2,4) \times \Sp(4,4)$, by~\cite{KlL90}.
All other cases except Ex.~2.4~(b) are excluded by the fact that $r = 4e+1$,
and Ex.~2.4~(b) does not apply because $d/\gcd(d,e)$ is odd.

\begin{verbatim}
    gap> SigmaFromMaxes( CharacterTable( "S6(4)" ), "85A",
    >    [ CharacterTable( "L4(4).2_2" ),
    >      CharacterTable( "A5" ) * CharacterTable( "S4(4)" ) ], [ 1, 1 ] );
    142/455
\end{verbatim}

This bound is not as good as the one obtained from the irreducible
element of order $65$ used above.

\begin{verbatim}
    gap> 16/63 < 142/455;
    true
\end{verbatim}

\subsection{$\ast$~$S_6(5)$}\label{S65}

We show that $S = S_6(5) = \PSp(6,5)$ satisfies the following.
\begin{enumerate}
\item[(a)]
    For $s \in S$ of the type $2 \perp 4$
    (i.~e., the preimage of $s$ in $\Sp(6,5) = 2.G$ decomposes
    the natural $6$-dimensional module for
    $\Sp(6,5)$ into an orthogonal sum of two irreducible modules
    of the dimensions $2$ and $4$, respectively) and of order $78$,
    $\M(S,s)$ consists of one group of the type
    $G_2 = 2.(\PSp(2,5) \times \PSp(4,5))$.
\item[(b)]
    $\total(S,s) = 9/217$.
\end{enumerate}

The order of $s$ is a multiple of the primitive prime divisor $r = 13$
of $5^4-1$,
so we have $e = 4$, $d = 6$, and $q = 5$, in the terminology of~\cite{GPPS}.
We check~\cite[Ex.~2.1--2.9]{GPPS}.
Ex.~2.1 does not apply because the classes $\calC_5$ and $\calC_8$
are empty by~\cite[Table~3.5.C]{KlL90},
Ex.~2.2 yields exactly the stabilizer $G_2$ of a $2$-dimensional subspace,
Ex.~2.4~(b) does not apply because $d/\gcd(d,e)$ is odd,
and all other cases are excluded by the fact that $r = 3e+1$.

The group $G_2$ has the structure $2.(\PSp(2,5) \times \PSp(4,5))$,
which is a central product of $\Sp(2,5) \cong 2.A_5$ and $\Sp(4,5) = 2.S_4(5)$
(see~\cite[Prop.~4.1.3]{KlL90}).
The character table of $G_2$ can be derived from that of the direct product
of $2.A_5$ and $2.S_4(5)$,
by factoring out the diagonal central subgroup of order two.

\begin{verbatim}
    gap> t:= CharacterTable( "S6(5)" );;
    gap> s1:= CharacterTable( "2.A5" );;
    gap> s2:= CharacterTable( "2.S4(5)" );;
    gap> dp:= s1 * s2;
    CharacterTable( "2.A5x2.S4(5)" )
    gap> c:= Difference( ClassPositionsOfCentre( dp ), Union(
    >                        GetFusionMap( s1, dp ), GetFusionMap( s2, dp ) ) );
    [ 62 ]
    gap> s:= dp / c;
    CharacterTable( "2.A5x2.S4(5)/[ 1, 62 ]" )
\end{verbatim}

Now we compute $\total(S,s)$.

\begin{verbatim}
    gap> SigmaFromMaxes( t, "78A", [ s ], [ 1 ] );
    9/217
\end{verbatim}

\subsection{$S_8(3)$}\label{S83}

We show that $S = S_8(3) = \PSp(8,3)$ satisfies the following.
\begin{enumerate}
\item[(a)]
    For $s \in S$ irreducible of order $41$,
    $\M(S,s)$ consists of one group $M$ of the type
    $S_4(9).2 = \PSp(4,9).2$.
\item[(b)]
    $\total(S,s) = 1/546$.
\item[(c)]
    The preimage of $s$ in the matrix group $2.S_8(3) = \Sp(8,3)$
    can be chosen of order $82$,
    and the preimage of $M$ is $2.S_4(9).2 = \Sp(4,9).2$.
\end{enumerate}

By~\cite{Be00}, the only maximal subgroups of $S$ that contain irreducible
elements of order $(3^4+1)/2 = 41$ are of extension field type,
and by~\cite[Prop.~4.3.10]{KlL90}, these groups have the structure $S_4(9).2$
and there is exactly one class of these groups.

The group $U = S_4(9)$ has three nontrivial outer automorphisms,
the character table of the subgroup $U.2$ in question has the identifier
\verb!"S4(9).2_1"!,
which follows from the fact that the extensions of $U$ by the other two
outer automorphisms do not admit a class fusion into $S$.

\begin{verbatim}
    gap> t:= CharacterTable( "S8(3)" );;
    gap> pi:= List( [ "S4(9).2_1", "S4(9).2_2", "S4(9).2_3" ],
    >               name -> PossiblePermutationCharacters(
    >                           CharacterTable( name ), t ) );;
    gap> List( pi, Length );
    [ 1, 0, 0 ]
\end{verbatim}

Now statement~(a) follows from the fact that $(1_{U.2})^S(s) = 1$ holds.

\begin{verbatim}
    gap> spos:= Position( OrdersClassRepresentatives( t ), 41 );;
    gap> pi[1][1][ spos ];
    1
\end{verbatim}

Now we compute $\total(S,s)$ in order to show statement~(b).

\begin{verbatim}
    gap> Maximum( ApproxP( pi[1], spos ) );
    1/546
\end{verbatim}

Statement~(c) is clear from the description of extension field type
subgroups in~\cite{KlL90}.

\subsection{$U_4(4)$}\label{U44}

We show that $S = U_4(4) = \SU(4,4)$ satisfies the following.
\begin{enumerate}
\item[(a)]
    For $s \in S$ of the type $1 \perp 3$
    (i.~e., $s$ decomposes the natural $4$-dimensional module for
    $\SU(4,4)$ into an orthogonal sum of two irreducible modules
    of the dimensions $1$ and $3$, respectively) and of order $4^3+1 = 65$,
    $\M(S,s)$ consists of one group of the type
    $G_1 = 5 \times U_3(4) = \GU(3,4)$.
\item[(b)]
    $\total(S,s) = 209/3\,264$.
\end{enumerate}

By~\cite{MSW94}, the only maximal subgroups of $S$ that contain $s$
are one class of stabilizers $H \cong 5 \times U_3(4)$
of this decomposition,
and clearly there is only one such group containing $s$.

Note that $H$ has index $3\,264$ in $S$,
since $S$ has two orbits on the $1$-dimensional subspaces,
of lengths $1\,105$ and $3\,264$, respectively,
and elements of order $13 = 65/5$ lie in the stabilizers of points
in the latter orbit.

\begin{verbatim}
    gap> g:= SU(4,4);;
    gap> orbs:= Orbits( g, NormedRowVectors( GF(16)^4 ), OnLines );;
    gap> orblen:= List( orbs, Length );
    [ 1105, 3264 ]
    gap> List( orblen, x -> x mod 13 );
    [ 0, 1 ]
\end{verbatim}

We compute the permutation character $1_{G_1}^S$;
there is exactly one combinatorially possible permutation character
of degree $3\,264$ (cf.~\cite{BP98copy}).

\begin{verbatim}
    gap> t:= CharacterTable( "U4(4)" );;
    gap> pi:= PermChars( t, rec( torso:= [ orblen[2] ] ) );;
    gap> Length( pi );
    1
\end{verbatim}

Now we compute $\total(S,s)$.

\begin{verbatim}
    gap> spos:= Position( OrdersClassRepresentatives( t ), 65 );;
    gap> Maximum( ApproxP( pi, spos ) );
    209/3264
\end{verbatim}

\subsection{$U_6(2)$}\label{U62}

We show that $S = U_6(2) = \PSU(6,2)$ satisfies the following.
\begin{enumerate}
\item[(a)]
    For $s \in S$ of order $11$,
    $\M(S,s)$ consists of one group of the type $U_5(2) = \SU(5,2)$
    and three groups of the type $M_{22}$.
\item[(b)]
    $\total(S,s) = 5/21$.
\item[(c)]
    The preimage of $s$ in the matrix group $\SU(6,2) = 3.U_6(2)$
    can be chosen of order $33$,
    and the preimages of the groups in $\M(S,s)$ have the structures
    $3 \times U_5(2) \cong \GU(5,2)$ and $3.M_{22}$, respectively.
\item[(d)]
    With $s$ as in~(a), the automorphic extensions $S.2$, $S.3$ of $S$
    satisfy $\total^{\prime}(S.2,s) = 5/96$ and
    $\total^{\prime}(S.3,s) = 59/224$.
\end{enumerate}

According to the list of maximal subgroups of $S$ in~\cite[p.~115]{CCN85},
$s$ is contained exactly in maximal subgroups of the types $U_5(2)$
(one class) and $M_{22}$ (three classes).

The permutation character of the action on the cosets of $U_5(2)$ type
subgroups is uniquely determined by the character tables.
We get three possibilities for the permutation character on the cosets of
$M_{22}$ type subgroups; they correspond to the three classes of such
subgroups, because each of these classes contains elements in exactly one
of the conjugacy classes {\tt 4C}, {\tt 4D}, and {\tt 4E} of elements in
$S$, and these classes are fused under the outer automorphism of $S$
of order three.

\begin{verbatim}
    gap> t:= CharacterTable( "U6(2)" );;
    gap> s1:= CharacterTable( "U5(2)" );;
    gap> pi1:= PossiblePermutationCharacters( s1, t );;
    gap> Length( pi1 );
    1
    gap> s2:= CharacterTable( "M22" );;
    gap> pi2:= PossiblePermutationCharacters( s2, t );
    [ Character( CharacterTable( "U6(2)" ), [ 20736, 0, 384, 0, 0, 0, 54, 0, 0, 
          0, 0, 48, 0, 16, 6, 0, 0, 0, 0, 0, 0, 6, 0, 2, 0, 0, 0, 4, 0, 0, 0, 0, 
          1, 1, 0, 0, 0, 0, 0, 0, 0, 0, 0, 0, 0, 0 ] ), 
      Character( CharacterTable( "U6(2)" ), [ 20736, 0, 384, 0, 0, 0, 54, 0, 0, 
          0, 48, 0, 0, 16, 6, 0, 0, 0, 0, 0, 0, 6, 0, 2, 0, 0, 4, 0, 0, 0, 0, 0, 
          1, 1, 0, 0, 0, 0, 0, 0, 0, 0, 0, 0, 0, 0 ] ), 
      Character( CharacterTable( "U6(2)" ), [ 20736, 0, 384, 0, 0, 0, 54, 0, 0, 
          48, 0, 0, 0, 16, 6, 0, 0, 0, 0, 0, 0, 6, 0, 2, 0, 4, 0, 0, 0, 0, 0, 0, 
          1, 1, 0, 0, 0, 0, 0, 0, 0, 0, 0, 0, 0, 0 ] ) ]
    gap> imgs:= Set( List( pi2, x -> Position( x, 48 ) ) );
    [ 10, 11, 12 ]
    gap> AtlasClassNames( t ){ imgs };
    [ "4C", "4D", "4E" ]
    gap> GetFusionMap( t, CharacterTable( "U6(2).3" ) ){ imgs };
    [ 10, 10, 10 ]
    gap> prim:= Concatenation( pi1, pi2 );;
\end{verbatim}

Now statement~(a) follows from the fact that the permutation characters
have the value $1$ on $s$.

\begin{verbatim}
    gap> spos:= Position( OrdersClassRepresentatives( t ), 11 );;
    gap> List( prim, x -> x[ spos ] );
    [ 1, 1, 1, 1 ]
\end{verbatim}

For statement~(b), we compute $\total(S,s)$.

\begin{verbatim}
    gap> Maximum( ApproxP( prim, spos ) );
    5/21
\end{verbatim}

Statement~(c) follows from~\cite{CCN85},
plus the information that $3.U_6(2)$ does not contain groups of the structure
$3 \times M_{22}$.

\begin{verbatim}
    gap> PossibleClassFusions(
    >        CharacterTable( "Cyclic", 3 ) * CharacterTable( "M22" ),
    >        CharacterTable( "3.U6(2)" ) );
    [  ]
\end{verbatim}

For statement~(d), we need that the relevant maximal subgroups of
$S.2$ are $U_5(2).2$ and one subgroup $M_{22}.2$,
and that the relevant maximal subgroup of $S.3$ is $U_5(2) \times 3$,
see~\cite[p.~115]{CCN85}.

\begin{verbatim}
    gap> SigmaFromMaxes( CharacterTable( "U6(2).2" ), "11AB",
    >        [ CharacterTable( "U5(2).2" ), CharacterTable( "M22.2" ) ],
    >        [ 1, 1 ], "outer" );
    5/96
    gap> SigmaFromMaxes( CharacterTable( "U6(2).3" ), "11A",
    >        [ CharacterTable( "U5(2)" ) * CharacterTable( "Cyclic", 3 ) ],
    >        [ 1 ], "outer" );
    59/224
\end{verbatim}

\section{Computations using Groups}\label{hard}

Before we start the computations using groups, we clean the workspace.

\begin{verbatim}
    gap> CleanWorkspace();
\end{verbatim}

\subsection{$A_{2m+1}$, $2 \leq m \leq 11$}\label{Aodd}

For alternating groups of odd degree $n = 2m+1$,
we choose $s$ to be an $n$-cycle.
The interesting cases in~\cite[Proposition~6.7]{BGK} are $5 \leq n \leq 23$.

In each case, we compute representatives of the maximal subgroups of $A_n$,
consider only those that contain an $n$-cycle, and then compute the
permutation characters.
Additionally, we show also the names that are used for the subgroups in
the {\GAP} Library of Transitive Groups,
see~\cite{HulpkeTG}
and the documentation of this library in the {\GAP} Reference Manual.

\begin{verbatim}
    gap> PrimitivesInfoForOddDegreeAlternatingGroup:= function( n )
    >     local G, max, cycle, spos, prim, nonz;
    > 
    >     G:= AlternatingGroup( n );
    > 
    >     # Compute representatives of the classes of maximal subgroups.
    >     max:= MaximalSubgroupClassReps( G );
    > 
    >     # Omit subgroups that cannot contain an `n'-cycle.
    >     max:= Filtered( max, m -> IsTransitive( m, [ 1 .. n ] ) );
    > 
    >     # Compute the permutation characters.
    >     cycle:= [];
    >     cycle[ n-1 ]:= 1;
    >     spos:= PositionProperty( ConjugacyClasses( CharacterTable( G ) ),
    >                c -> CycleStructurePerm( Representative( c ) ) = cycle );
    >     prim:= List( max, m -> TrivialCharacter( m )^G );
    >     nonz:= PositionsProperty( prim, x -> x[ spos ] <> 0 );
    > 
    >     # Compute the subgroup names and the multiplicities.
    >     return rec( spos := spos,
    >                 prim := prim{ nonz },
    >                 grps := List( max{ nonz },
    >                               m -> TransitiveGroup( n,
    >                                        TransitiveIdentification( m ) ) ),
    >                 mult := List( prim{ nonz }, x -> x[ spos ] ) );
    > end;;
\end{verbatim}

The sets $\MM(s)$ and the values $\total(A_n,s)$ are as follows.
For each degree in question, the first list shows names for representatives
of the conjugacy classes of maximal subgroups containing a fixed $n$-cycle,
and the second list shows the number of conjugates in each class.

\begin{verbatim}
    gap> for n in [ 5, 7 .. 23 ] do
    >      prim:= PrimitivesInfoForOddDegreeAlternatingGroup( n );
    >      bound:= Maximum( ApproxP( prim.prim, prim.spos ) );
    >      Print( n, ": ", prim.grps, ", ", prim.mult, ", ", bound, "\n" );
    > od;
    5: [ D(5) = 5:2 ], [ 1 ], 1/3
    7: [ L(7) = L(3,2), L(7) = L(3,2) ], [ 1, 1 ], 2/5
    9: [ 1/2[S(3)^3]S(3), L(9):3=P|L(2,8) ], [ 1, 3 ], 9/35
    11: [ M(11), M(11) ], [ 1, 1 ], 2/105
    13: [ F_78(13)=13:6, L(13)=PSL(3,3), L(13)=PSL(3,3) ], [ 1, 2, 2 ], 4/1155
    15: [ 1/2[S(3)^5]S(5), 1/2[S(5)^3]S(3), L(15)=A_8(15)=PSL(4,2), 
      L(15)=A_8(15)=PSL(4,2) ], [ 1, 1, 1, 1 ], 29/273
    17: [ L(17):4=PYL(2,16), L(17):4=PYL(2,16) ], [ 1, 1 ], 2/135135
    19: [ F_171(19)=19:9 ], [ 1 ], 1/6098892800
    21: [ t21n150, t21n161, t21n91 ], [ 1, 1, 2 ], 29/285
    23: [ M(23), M(23) ], [ 1, 1 ], 2/130945815
\end{verbatim}

In the above output,
a subgroup printed as \verb|1/2[S(|$n_1$\verb|)^|$n_2$\verb|]S(|$n_2$\verb|)|,
where $n = n_1 n_2$ holds, denotes the intersection of $A_n$ with the
wreath product
\tthdump{$S_{n_1} \wr S_{n_2} \leq S_n$.}
(Note that the {\ATLAS} denotes the subgroup \verb|1/2[S(3)^3]S(3)| of $A_9$
as $3^3:S_4$.)
The groups printed as \verb!P|L(2,8)! and \verb!PYL(2,16)! denote
$\PGammaL(2,8)$ and $\PGammaL(2,16)$, respectively.
And the three subgroups of $A_{21}$ have the structures
\tthdump{$(S_3 \wr S_7) \cap A_{21}$, $(S_7 \wr S_3) \cap A_{21}$,}
and $\PGL(3,4)$, respectively.

Note that $A_9$ contains two conjugacy classes of maximal subgroups of
the type $\PGammaL(2,8) \cong L_2(8):3$, and that each $9$-cycle in $A_9$
is contained in exactly three \emph{conjugate} subgroups of this type.
For $n \in \{ 13, 15, 17 \}$, $A_n$ contains two conjugacy classes of
isomorphic maximal subgroups of linear type, and each $n$-cycle is contained
in subgroups from each class.
Finally, $A_{21}$ contains only one class of maximal subgroups of linear type.

For the two groups $A_5$ and $A_7$,
the values computed above are not sufficient.
See Section~\ref{A5} and~\ref{A7} for a further treatment.

The above computations look like a brute-force approach,
but note that the computation of the maximal subgroups of alternating
and symmetric groups in {\GAP} uses the classification of these
subgroups, and also the conjugacy classes of elements in alternating and
symmetric groups can be computed cheaply.

Alternative (character-theoretic) computations
for $n \in \{ 5, 7, 9, 11, 13 \}$ were shown in Section~\ref{easyloop}.
(A hand calculation for the case $n = 19$ can be found in~\cite{BW1}.)

\subsection{$A_5$}\label{A5}

We show that $S = A_5$ satisfies the following.
\begin{enumerate}
\item[(a)]
    $\total(S) = 1/3$,
    and this value is attained exactly for $\total(S,s)$
    with $s$ of order $5$.
\item[(b)]
    For $s \in S$ of order $5$,
    $\M(S,s)$ consists of one group of the type $D_{10}$.
\item[(c)]
    $\prop(S) = 1/3$,
    and this value is attained exactly for $\prop(S,s)$
    with $s$ of order $5$.
\item[(d)]
    Each element in $S$ together with one of
    $(1,2)(3,4)$, $(1,3)(2,4)$, $(1,4)(2,3)$
    generates a proper subgroup of $S$.
\item[(e)]
    Both the spread and the uniform spread of $S$ is exactly two
    (see~\cite{BW1}),
    with $s$ of order $5$.
\end{enumerate}

Statement~(a) follows from inspection of the primitive permutation
characters, cf.~Section~\ref{easyloop}.

\begin{verbatim}
    gap> t:= CharacterTable( "A5" );;
    gap> ProbGenInfoSimple( t );
    [ "A5", 1/3, 2, [ "5A" ], [ 1 ] ]
\end{verbatim}

Statement~(b) can be read off from the primitive permutation characters,
and the fact that the unique class of maximal subgroups that contain
elements of order $5$ consists of groups of the structure $D_{10}$,
see~\cite[p.~2]{CCN85}.

\begin{verbatim}
    gap> OrdersClassRepresentatives( t );
    [ 1, 2, 3, 5, 5 ]
    gap> PrimitivePermutationCharacters( t );
    [ Character( CharacterTable( "A5" ), [ 5, 1, 2, 0, 0 ] ), 
      Character( CharacterTable( "A5" ), [ 6, 2, 0, 1, 1 ] ), 
      Character( CharacterTable( "A5" ), [ 10, 2, 1, 0, 0 ] ) ]
\end{verbatim}

For statement~(c), we compute that for all nonidentity elements $s \in S$
and involutions $g \in S$,
$\prop(g,s) \geq 1/3$ holds,
with equality if and only if $s$ has order $5$.
We actually compute, for class representatives $s$,
the proportion of involutions $g$ such that
$\langle g, s \rangle \not= S$ holds.

\begin{verbatim}
    gap> g:= AlternatingGroup( 5 );;
    gap> inv:= g.1^2 * g.2;
    (1,4)(2,5)
    gap> cclreps:= List( ConjugacyClasses( g ), Representative );;
    gap> SortParallel( List( cclreps, Order ), cclreps );
    gap> List( cclreps, Order );
    [ 1, 2, 3, 5, 5 ]
    gap> Size( ConjugacyClass( g, inv ) );
    15
    gap> prop:= List( cclreps,
    >                 r -> RatioOfNongenerationTransPermGroup( g, inv, r ) );
    [ 1, 1, 3/5, 1/3, 1/3 ]
    gap> Minimum( prop );
    1/3
\end{verbatim}

Statement~(d) follows by explicit computations.

\begin{verbatim}
    gap> triple:= [ (1,2)(3,4), (1,3)(2,4), (1,4)(2,3) ];;
    gap> CommonGeneratorWithGivenElements( g, cclreps, triple );
    fail
\end{verbatim}

As for statement~(e), we know from~(a) that the uniform spread of $S$
is at least two, and from~(d) that the spread is less than three.

\subsection{$A_6$}\label{A6}

We show that $S = A_6$ satisfies the following.
\begin{enumerate}
\item[(a)]
    $\total(S) = 2/3$,
    and this value is attained exactly for $\total(S,s)$
    with $s$ of order $5$.
\item[(b)]
    For $s$ of order $5$,
    $\M(S,s)$ consists of two nonconjugate groups of the type $A_5$.
\item[(c)]
    $\prop(S) = 5/9$,
    and this value is attained exactly for $\prop(S,s)$
    with $s$ of order $5$.
\item[(d)]
    Each element in $S$ together with one of
    $(1,2)(3,4)$, $(1,3)(2,4)$, $(1,4)(2,3)$
    generates a proper subgroup of $S$.
\item[(e)]
    Both the spread and the uniform spread of $S$ is exactly two
    (see~\cite{BW1}),
    with $s$ of order $4$.
\item[(f)]
    For $x$, $y \in S_6^{\times}$, there is $s \in S_6$
    such that $S \subseteq \langle x, s \rangle \cap \langle y, s \rangle$.
    It is \emph{not} possible to find $s \in S$ with this property,
    or $s$ in a prescribed conjugacy class of $S_6$.
\item[(g)]
    $\total( \PGL(2,9) ) = 1/6$ and $\total( M_{10} ) = 1/9$,
    with $s$ of order $10$ and $8$, respectively.
\end{enumerate}

(Note that in this example, the optimal choice of $s$ for $\prop(S)$ cannot be
used to obtain the result on the exact spread.)

Statement~(a) follows from inspection of the primitive permutation
characters, cf.~Section~\ref{easyloop}.

\begin{verbatim}
    gap> t:= CharacterTable( "A6" );;
    gap> ProbGenInfoSimple( t );
    [ "A6", 2/3, 1, [ "5A" ], [ 2 ] ]
\end{verbatim}

Statement~(b) can be read off from the permutation characters,
and the fact that the two classes of maximal subgroups that contain
elements of order $5$ consist of groups of the structure $A_5$,
see~\cite[p.~4]{CCN85}.

\begin{verbatim}
    gap> OrdersClassRepresentatives( t );
    [ 1, 2, 3, 3, 4, 5, 5 ]
    gap> prim:= PrimitivePermutationCharacters( t );
    [ Character( CharacterTable( "A6" ), [ 6, 2, 3, 0, 0, 1, 1 ] ), 
      Character( CharacterTable( "A6" ), [ 6, 2, 0, 3, 0, 1, 1 ] ), 
      Character( CharacterTable( "A6" ), [ 10, 2, 1, 1, 2, 0, 0 ] ), 
      Character( CharacterTable( "A6" ), [ 15, 3, 3, 0, 1, 0, 0 ] ), 
      Character( CharacterTable( "A6" ), [ 15, 3, 0, 3, 1, 0, 0 ] ) ]
\end{verbatim}

For statement~(c),
we first compute that for all nonidentity elements $s \in S$
and involutions $g \in S$,
$\prop(g,s) \geq 5/9$ holds,
with equality if and only if $s$ has order $5$.
We actually compute, for class representatives $s$,
the proportion of involutions $g$ such that
$\langle g, s \rangle \not= S$ holds.

\begin{verbatim}
    gap> S:= AlternatingGroup( 6 );;
    gap> inv:= (S.1*S.2)^2;
    (1,3)(2,5)
    gap> cclreps:= List( ConjugacyClasses( S ), Representative );;
    gap> SortParallel( List( cclreps, Order ), cclreps );
    gap> List( cclreps, Order );
    [ 1, 2, 3, 3, 4, 5, 5 ]
    gap> C:= ConjugacyClass( S, inv );;
    gap> Size( C );
    45
    gap> prop:= List( cclreps,
    >                 r -> RatioOfNongenerationTransPermGroup( S, inv, r ) );
    [ 1, 1, 1, 1, 29/45, 5/9, 5/9 ]
    gap> Minimum( prop );
    5/9
\end{verbatim}

Now statement~(c) follows from the fact that for $g \in S$ of order larger
than two, $\total(S,g) \leq 1/2 < 5/9$ holds.

\begin{verbatim}
    gap> ApproxP( prim, 6 );
    [ 0, 2/3, 1/2, 1/2, 0, 1/3, 1/3 ]
\end{verbatim}

Statement~(d) follows by explicit computations.

\begin{verbatim}
    gap> triple:= [ (1,2)(3,4), (1,3)(2,4), (1,4)(2,3) ];;
    gap> CommonGeneratorWithGivenElements( S, cclreps, triple );
    fail
\end{verbatim}

An alternative triple to that in statement~(d) is the one given
in~\cite{BW1}.

\begin{verbatim}
    gap> triple:= [ (1,3)(2,4), (1,5)(2,6), (3,6)(4,5) ];;
    gap> CommonGeneratorWithGivenElements( S, cclreps, triple );
    fail
\end{verbatim}

Of course we can also construct such a triple, as follows.

\begin{verbatim}
    gap> TripleWithProperty( [ [ inv ], C, C ],
    >        l -> ForAll( S, elm ->
    >   ForAny( l, x -> not IsGeneratorsOfTransPermGroup( S, [ elm, x ] ) ) ) );
    [ (1,3)(2,5), (1,3)(2,6), (1,3)(2,4) ]
\end{verbatim}

For statement~(e), we use the random approach described in
Section~\ref{groups}.

\begin{verbatim}
    gap> s:= (1,2,3,4)(5,6);;
    gap> reps:= Filtered( cclreps, x -> Order( x ) > 1 );;
    gap> ResetGlobalRandomNumberGenerators();
    gap> for pair in UnorderedTuples( reps, 2 ) do
    >      if RandomCheckUniformSpread( S, pair, s, 40 ) <> true then
    >        Print( "#E  nongeneration!\n" );
    >      fi;
    >    od;
\end{verbatim}

We get no output, so a suitable element of order $4$ works in all cases.
Note that we cannot use an element of order $5$,
because it fixes a point in the natural permutation representation,
and we may take $x_1 = (1,2,3)$ and $x_2 = (4,5,6)$.
With this argument, only elements of order $4$ and double $3$-cycles
are possible choices for $s$, and the latter are excluded by the fact
that an outer automorphism maps the class of double $s$-cycles in $A_6$
to the class of $3$-cycles.
So no element in $A_6$ of order different from $4$ works.

Next we show statement~(f).
Already in $A_6.2_1 = S_6$, elements $s$ of order $4$ do in general not work
because they do not generate with transpositions.

\begin{verbatim}
    gap> G:= SymmetricGroup( 6 );;
    gap> RatioOfNongenerationTransPermGroup( G, s, (1,2) );
    1
\end{verbatim}

Also, choosing $s$ from a prescribed conjugacy class of $S_6$ (that is,
also $s$ outside $A_6$ is allowed) with the property that
$A_6 \subseteq \langle x, s \rangle \cap \langle y, s \rangle$
is not possible.
Note that only $6$-cycles are possible for $s$ if $x$ and $y$ are commuting
transpositions, and --applying the outer automorphism--
no $6$-cycle works for two commuting fixed-point free involutions.
(The group is small enough for a brute force test.)

\begin{verbatim}
    gap> goods:= Filtered( Elements( G ),
    >      s -> IsGeneratorsOfTransPermGroup( G, [ s, (1,2) ] ) and
    >           IsGeneratorsOfTransPermGroup( G, [ s, (3,4) ] ) );;
    gap> Collected( List( goods, CycleStructurePerm ) );
    [ [ [ ,,,, 1 ], 24 ] ]
    gap> goods:= Filtered( Elements( G ),
    >      s -> IsGeneratorsOfTransPermGroup( G, [ s, (1,2)(3,4)(5,6) ] ) and
    >           IsGeneratorsOfTransPermGroup( G, [ s, (1,3)(2,4)(5,6) ] ) );;
    gap> Collected( List( goods, CycleStructurePerm ) );
    [ [ [ 1, 1 ], 24 ] ]
\end{verbatim}

However, for each pair of nonidentity element $x$, $y \in S_6$,
there is $s \in S_6$ such that
$\langle x, s \rangle$ and $\langle y, s \rangle$ both contain $A_6$.
(If $s$ works for the pair $(x,y)$ then $s^g$ works for $(x^g,y^g)$,
so it is sufficient to consider only orbit representatives $(x,y)$ under
the conjugation action of $G$ on pairs.
Thus we check conjugacy class representatives $x$ and, for fixed $x$,
representatives of orbits of $C_G(x)$ on the classes $y^G$,
i.~e., representatives of $C_G(y)$-$C_G(x)$-double cosets in $G$.
Moreover, clearly we can restrict the checks to elements $x$, $y$ of
prime order.)

\begin{verbatim}
    gap> Sgens:= GeneratorsOfGroup( S );;
    gap> primord:= Filtered( List( ConjugacyClasses( G ), Representative ),
    >                        x -> IsPrimeInt( Order( x ) ) );;
    gap> for x in primord do
    >      for y in primord do
    >        for pair in DoubleCosetRepsAndSizes( G, Centralizer( G, y ),
    >                        Centralizer( G, x ) ) do
    >          if not ForAny( G, s -> IsSubset( Group( x,s ), S ) and 
    >                                 IsSubset( Group( y^pair[1], s ), S ) ) then
    >            Error( [ x, y ] );
    >          fi;
    >        od;
    >      od;
    >    od;
\end{verbatim}

In other words, the spread of $S_6$ is $2$ but the uniform spread of $S_6$
is not $2$ but only $1$.

We cannot always find $s \in A_6$ with the required property:
If $x$ is a transposition then any $s$ with $S \subseteq\langle x, s \rangle$
must be a $5$-cycle.

\begin{verbatim}
    gap> filt:= Filtered( S, s -> IsSubset( Group( (1,2), s ), S ) );;
    gap> Collected( List( filt, Order ) );
    [ [ 5, 48 ] ]
\end{verbatim}

Moreover, clearly such $s$ fixes one of the moved points of $x$,
so we may prescribe a transposition $y \not= x$ that commutes with $x$,
it satisfies $S \not\subseteq\langle y, s \rangle$.

For the other two automorphic extensions $A_6.2_2 = \PGL(2,9)$
and $A_6.2_3 = M_{10}$,
we compute the character-theoretic bounds $\total(A_6.2_2) = 1/6$
and $\total(A_6.2_3) = 1/9$,
which shows statement~(g).

\begin{verbatim}
    gap> ProbGenInfoSimple( CharacterTable( "A6.2_2" ) );
    [ "A6.2_2", 1/6, 5, [ "10A" ], [ 1 ] ]
    gap> ProbGenInfoSimple( CharacterTable( "A6.2_3" ) );
    [ "A6.2_3", 1/9, 8, [ "8C" ], [ 1 ] ]
\end{verbatim}

Note that $\total^{\prime}( \PGL(2,9), s ) = 1/6$,
with $s$ of order $5$,
and $\total^{\prime}( M_{10}, s ) = 0$ for any $s \in A_6$
since $M_{10}$ is a non-split extension of $A_6$.

\begin{verbatim}
    gap> t:= CharacterTable( "A6" );;
    gap> t2:= CharacterTable( "A6.2_2" );;
    gap> spos:= PositionsProperty( OrdersClassRepresentatives( t ), x -> x = 5 );;
    gap> ProbGenInfoAlmostSimple( t, t2, spos );
    [ "A6.2_2", 1/6, [ "5A", "5B" ], [ 1, 1 ] ]
\end{verbatim}

\subsection{$A_7$}\label{A7}

We show that $S = A_7$ satisfies the following.
\begin{enumerate}
\item[(a)]
    $\total(S) = 2/5$,
    and this value is attained exactly for $\total(S,s)$
    with $s$ of order $7$.
\item[(b)]
    For $s$ of order $7$,
    $\M(S,s)$ consists of two nonconjugate subgroups of the type $L_2(7)$.
\item[(c)]
    $\prop(S) = 2/5$,
    and this value is attained exactly for $\prop(S,s)$
    with $s$ of order $7$.
\item[(d)]
    The uniform spread of $S$ is exactly three,
    with $s$ of order $7$.
\end{enumerate}

Statement~(a) follows from inspection of the primitive permutation
characters, cf.~Section~\ref{easyloop}.

\begin{verbatim}
    gap> t:= CharacterTable( "A7" );;
    gap> ProbGenInfoSimple( t );
    [ "A7", 2/5, 2, [ "7A" ], [ 2 ] ]
\end{verbatim}

Statement~(b) can be read off from the permutation characters,
and the fact that the two classes of maximal subgroups that contain
elements of order $7$ consist of groups of the structure $L_2(7)$,
see~\cite[p.~10]{CCN85}.

\begin{verbatim}
    gap> OrdersClassRepresentatives( t );
    [ 1, 2, 3, 3, 4, 5, 6, 7, 7 ]
    gap> prim:= PrimitivePermutationCharacters( t );
    [ Character( CharacterTable( "A7" ), [ 7, 3, 4, 1, 1, 2, 0, 0, 0 ] ), 
      Character( CharacterTable( "A7" ), [ 15, 3, 0, 3, 1, 0, 0, 1, 1 ] ), 
      Character( CharacterTable( "A7" ), [ 15, 3, 0, 3, 1, 0, 0, 1, 1 ] ), 
      Character( CharacterTable( "A7" ), [ 21, 5, 6, 0, 1, 1, 2, 0, 0 ] ), 
      Character( CharacterTable( "A7" ), [ 35, 7, 5, 2, 1, 0, 1, 0, 0 ] ) ]
\end{verbatim}

For statement~(c), we compute that for all nonidentity elements $s \in S$
and involutions $g \in S$,
$\prop(g,s) \geq 2/5$ holds,
with equality if and only if $s$ has order $7$.
We actually compute, for class representatives $s$,
the proportion of involutions $g$ such that
$\langle g, s \rangle \not= S$ holds.

\begin{verbatim}
    gap> g:= AlternatingGroup( 7 );;
    gap> inv:= (g.1^3*g.2)^3;
    (2,6)(3,7)
    gap> ccl:= List( ConjugacyClasses( g ), Representative );;
    gap> SortParallel( List( ccl, Order ), ccl );
    gap> List( ccl, Order );
    [ 1, 2, 3, 3, 4, 5, 6, 7, 7 ]
    gap> Size( ConjugacyClass( g, inv ) );
    105
    gap> prop:= List( ccl, r -> RatioOfNongenerationTransPermGroup( g, inv, r ) );
    [ 1, 1, 1, 1, 89/105, 17/21, 19/35, 2/5, 2/5 ]
    gap> Minimum( prop );
    2/5
\end{verbatim}

For statement~(d),
we use the random approach described in Section~\ref{groups}.
By the character-theoretic bounds, it suffices to consider triples
of elements in the classes {\tt 2A} or {\tt 3B}.

\begin{verbatim}
    gap> OrdersClassRepresentatives( t );
    [ 1, 2, 3, 3, 4, 5, 6, 7, 7 ]
    gap> spos:= Position( OrdersClassRepresentatives( t ), 7 );;
    gap> SizesCentralizers( t );
    [ 2520, 24, 36, 9, 4, 5, 12, 7, 7 ]
    gap> ApproxP( prim, spos );
    [ 0, 2/5, 0, 2/5, 2/15, 0, 0, 2/15, 2/15 ]
    gap> s:= (1,2,3,4,5,6,7);;
    gap> 3B:= (1,2,3)(4,5,6);;
    gap> C3B:= ConjugacyClass( g, 3B );;
    gap> Size( C3B );
    280
    gap> ResetGlobalRandomNumberGenerators();
    gap> for triple in UnorderedTuples( [ inv, 3B ], 3 ) do
    >      if RandomCheckUniformSpread( g, triple, s, 80 ) <> true then
    >        Print( "#E  nongeneration!\n" );
    >      fi;
    >    od;
\end{verbatim}

We get no output, so the uniform spread of $S$ is at least three.

Alternatively, we can use Lemma~\ref{existsgoodconjugate};
this approach is technically more involved but faster.
We work with the diagonal product of the two degree $15$ representations
of $S$,
which is constructed from the information stored in the {\GAP} Library
of Tables of Marks.

\begin{verbatim}
    gap> tom:= TableOfMarks( "A7" );;
    gap> a7:= UnderlyingGroup( tom );;
    gap> tommaxes:= MaximalSubgroupsTom( tom );
    [ [ 39, 38, 37, 36, 35 ], [ 7, 15, 15, 21, 35 ] ]
    gap> index15:= List( tommaxes[1]{ [ 2, 3 ] },
    >                    i -> RepresentativeTom( tom, i ) );
    [ Group([ (1,3)(2,7), (1,5,7)(3,4,6) ]), 
      Group([ (1,4)(2,3), (2,4,6)(3,5,7) ]) ]
    gap> deg15:= List( index15, s -> RightTransversal( a7, s ) );;
    gap> reps:= List( deg15, l -> Action( a7, l, OnRight ) );
    [ Group([ (1,5,7)(2,9,10)(3,11,4)(6,12,8)(13,14,15), 
          (1,8,15,5,12)(2,13,11,3,10)(4,14,9,7,6) ]), 
      Group([ (1,2,3)(4,6,5)(7,8,9)(10,12,11)(13,15,14), 
          (1,12,3,13,10)(2,9,15,4,11)(5,6,14,7,8) ]) ]
    gap> g:= DiagonalProductOfPermGroups( reps );;
    gap> ResetGlobalRandomNumberGenerators();
    gap> repeat s:= Random( g );
    >    until Order( s ) = 7;
    gap> NrMovedPoints( s );
    28
    gap> mpg:= MovedPoints( g );;
    gap> fixs:= Difference( mpg, MovedPoints( s ) );;
    gap> orb_s:= Orbit( g, fixs, OnSets );;
    gap> Length( orb_s );
    120
    gap> SizesCentralizers( t );
    [ 2520, 24, 36, 9, 4, 5, 12, 7, 7 ]
    gap> repeat 2a:= Random( g ); until Order( 2a ) = 2;
    gap> repeat 3b:= Random( g );
    >    until Order( 3b ) = 3 and Size( Centralizer( g, 3b ) ) = 9;
    gap> orb2a:= Orbit( g, Difference( mpg, MovedPoints( 2a ) ), OnSets );;
    gap> orb3b:= Orbit( g, Difference( mpg, MovedPoints( 3b ) ), OnSets );;
    gap> orb2aor3b:= Union( orb2a, orb3b );;
    gap> TripleWithProperty( [ [ orb2a[1], orb3b[1] ], orb2aor3b, orb2aor3b ],
    >        l -> ForAll( orb_s,
    >                 f -> not IsEmpty( Intersection( Union( l ), f ) ) ) );
    fail
\end{verbatim}

It remains to show that for any choice of $s \in S$,
a quadruple of elements in $S^{\times}$ exists such that $s$ generates
a proper subgroup of $S$ together with at least one of these elements.

First we observe (without using {\GAP}) that there is a pair of $3$-cycles
whose fixed points cover the seven points of the natural permutation
representation.
This implies the statement for all elements $s \in S$
that fix a point in this representation.
So it remains to consider elements $s$ of the orders six and seven.

For the order seven element, the above setup and
Lemma~\ref{existsgoodconjugate} can be used.

\begin{verbatim}
    gap> QuadrupleWithProperty( [ [ orb2a[1] ], orb2a, orb2a, orb2a ],
    >        l -> ForAll( orb_s,
    >                 f -> not IsEmpty( Intersection( Union( l ), f ) ) ) );
    [ [ 2, 11, 15, 19, 24, 29 ], [ 4, 9, 13, 21, 22, 28 ], 
      [ 3, 10, 14, 20, 23, 30 ], [ 1, 5, 7, 25, 26, 27 ] ]
\end{verbatim}

For the order six element, we use the diagonal product of the
primitive permutation representations of the degrees $21$ and $35$.

\begin{verbatim}
    gap> has6A:= List( tommaxes[1]{ [ 4, 5 ] },
    >                  i -> RepresentativeTom( tom, i ) );
    [ Group([ (1,2)(3,7), (2,6,5,4)(3,7) ]), 
      Group([ (2,3)(5,7), (1,2)(4,5,6,7), (2,3)(5,6) ]) ]
    gap> trans:= List( has6A, s -> RightTransversal( a7, s ) );;
    gap> reps:= List( trans, l -> Action( a7, l, OnRight ) );
    [ Group([ (1,16,12)(2,17,13)(3,18,11)(4,19,14)(15,20,21), 
          (1,4,7,9,10)(2,5,8,3,6)(11,12,15,14,13)(16,20,19,17,18) ]), 
      Group([ (2,16,6)(3,17,7)(4,18,8)(5,19,9)(10,20,26)(11,21,27)(12,22,28)(13,
            23,29)(14,24,30)(15,25,31), (1,2,3,4,5)(6,10,13,15,9)(7,11,14,8,
            12)(16,20,23,25,19)(17,21,24,18,22)(26,32,35,31,28)(27,33,29,34,30) 
         ]) ]
    gap> g:= DiagonalProductOfPermGroups( reps );;
    gap> repeat s:= Random( g );
    >    until Order( s ) = 6;
    gap> NrMovedPoints( s );
    53
    gap> mpg:= MovedPoints( g );;
    gap> fixs:= Difference( mpg, MovedPoints( s ) );;
    gap> orb_s:= Orbit( g, fixs, OnSets );;
    gap> Length( orb_s );
    105
    gap> repeat 3a:= Random( g );
    >    until Order( 3a ) = 3 and Size( Centralizer( g, 3a ) ) = 36;
    gap> orb3a:= Orbit( g, Difference( mpg, MovedPoints( 3a ) ), OnSets );;
    gap> Length( orb3a );
    35
    gap> TripleWithProperty( [ [ orb3a[1] ], orb3a, orb3a ],
    >        l -> ForAll( orb_s,
    >                 f -> not IsEmpty( Intersection( Union( l ), f ) ) ) );
    [ [ 5, 7, 9, 17, 18, 19, 31, 32, 34, 40, 53 ], 
      [ 5, 7, 9, 11, 13, 14, 30, 41, 42, 44, 53 ], 
      [ 2, 3, 4, 5, 7, 9, 26, 47, 48, 50, 53 ] ]
\end{verbatim}

So we have found not only a quadruple but even a triple of $3$-cycles
that excludes candidates $s$ of order six.

\subsection{$L_d(q)$}\label{SL}

In the treatment of small dimensional linear groups $S = \SL(d,q)$,
\cite{BGK} uses a Singer element $s$ of order $(q^d-1)/(q-1)$.
(So the order of the corresponding element in
$\PSL(d,q) = (q^d-1)/[(q-1) \gcd(d,q-1)]$.)
By~\cite{Be00}, $\M(S,s)$ consists of extension field type subgroups,
except in the cases $d = 2$, $q \in \{ 2, 5, 7, 9 \}$, and $(d,q) = (3,4)$.
These subgroups have the structure $\GL(d/p,q^p):\alpha_q \cap S$,
for prime divisors $p$ of $d$,
where $\alpha_q$ denotes the Frobenius automorphism that acts on
matrices by raising each entry to the $q$-th power.
(If $q$ is a prime then we have $\GL(d/p,q^p):\alpha_q = \GammaL(d/p,q^p)$.)
Since $s$ acts irreducibly,
it is contained in at most one conjugate of each class of
extension field type subgroups (cf.~\cite[Lemma~2.12]{BGK}).

First we write a {\GAP} function \verb|RelativeSigmaL| that takes
a positive integer $d$ and a basis $B$ of the field
extension of degree $n$ over the field with $q$ elements,
and returns the group $\GL(d,q^n):\alpha_q$, as a subgroup of $\GL(dn,q)$.

\begin{verbatim}
    gap> RelativeSigmaL:= function( d, B )
    >     local n, F, q, glgens, diag, pi, frob, i;
    > 
    >     n:= Length( B );
    >     F:= LeftActingDomain( UnderlyingLeftModule( B ) );
    >     q:= Size( F );
    > 
    >     # Create the generating matrices inside the linear subgroup.
    >     glgens:= List( GeneratorsOfGroup( SL( d, q^n ) ),
    >                    m -> BlownUpMat( B, m ) );
    > 
    >     # Create the matrix of a diagonal part that maps to determinant 1.
    >     diag:= IdentityMat( d*n, F );
    >     diag{ [ 1 .. n ] }{ [ 1 .. n ] }:= BlownUpMat( B, [ [ Z(q^n)^(q-1) ] ] );
    >     Add( glgens, diag );
    > 
    >     # Create the matrix that realizes the Frobenius action,
    >     # and adjust the determinant.
    >     pi:= List( B, b -> Coefficients( B, b^q ) );
    >     frob:= NullMat( d*n, d*n, F );
    >     for i in [ 0 .. d-1 ] do
    >       frob{ [ 1 .. n ] + i*n }{ [ 1 .. n ] + i*n }:= pi;
    >     od;
    >     diag:= IdentityMat( d*n, F );
    >     diag{ [ 1 .. n ] }{ [ 1 .. n ] }:= BlownUpMat( B, [ [ Z(q^n) ] ] );
    >     diag:= diag^LogFFE( Inverse( Determinant( frob ) ), Determinant( diag ) );
    > 
    >     # Return the result.
    >     return Group( Concatenation( glgens, [ diag * frob ] ) );
    > end;;
\end{verbatim}

The next function computes $\total(\SL(d,q),s)$,
by computing the sum of $\fpr(g,S/(\GL(d/p,q^p):\alpha_q \cap S))$,
for prime divisors $p$ of $d$, and taking the maximum over
$g \in S^{\times}$.
The computations take place in a permutation representation of $\PSL(d,q)$.


\begin{verbatim}
    gap> ApproxPForSL:= function( d, q )
    >     local G, epi, PG, primes, maxes, names, ccl;
    > 
    >     # Check whether this is an admissible case (see [Be00]).
    >     if ( d = 2 and q in [ 2, 5, 7, 9 ] ) or ( d = 3 and q = 4 ) then
    >       return fail;
    >     fi;
    > 
    >     # Create the group SL(d,q), and the map to PSL(d,q).
    >     G:= SL( d, q );
    >     epi:= ActionHomomorphism( G, NormedRowVectors( GF(q)^d ), OnLines );
    >     PG:= ImagesSource( epi );
    > 
    >     # Create the subgroups corresponding to the prime divisors of `d'.
    >     primes:= Set( Factors( d ) );
    >     maxes:= List( primes, p -> RelativeSigmaL( d/p,
    >                                  Basis( AsField( GF(q), GF(q^p) ) ) ) );
    >     names:= List( primes, p -> Concatenation( "GL(", String( d/p ), ",",
    >                                  String( q^p ), ").", String( p ) ) );
    >     if 2 < q then
    >       names:= List( names, name -> Concatenation( name, " cap G" ) );
    >     fi;
    > 
    >     # Compute the conjugacy classes of prime order elements in the maxes.
    >     # (In order to avoid computing all conjugacy classes of these subgroups,
    >     # we work in Sylow subgroups.)
    >     ccl:= List( List( maxes, x -> ImagesSet( epi, x ) ),
    >             M -> ClassesOfPrimeOrder( M, Set( Factors( Size( M ) ) ),
    >                                       TrivialSubgroup( M ) ) );
    > 
    >     return [ names, UpperBoundFixedPointRatios( PG, ccl, true )[1] ];
    > end;;
\end{verbatim}


We apply this function to the cases that are interesting
in~\cite[Section~5.12]{BGK}.

\begin{verbatim}
    gap> pairs:= [ [ 3, 2 ], [ 3, 3 ], [ 4, 2 ], [ 4, 3 ], [ 4, 4 ],
    >            [ 6, 2 ], [ 6, 3 ], [ 6, 4 ], [ 6, 5 ], [ 8, 2 ], [ 10, 2 ] ];;
    gap> array:= [];;
    gap> for pair in pairs do
    >      d:= pair[1];  q:= pair[2];
    >      approx:= ApproxPForSL( d, q );
    >      Add( array, [ Concatenation( "SL(", String(d), ",", String(q), ")" ),
    >                    (q^d-1)/(q-1),
    >                    approx[1], approx[2] ] );
    >    od;
    gap> PrintFormattedArray( array );
       SL(3,2)    7                             [ "GL(1,8).3" ]             1/4
       SL(3,3)   13                      [ "GL(1,27).3 cap G" ]            1/24
       SL(4,2)   15                             [ "GL(2,4).2" ]            3/14
       SL(4,3)   40                       [ "GL(2,9).2 cap G" ]         53/1053
       SL(4,4)   85                      [ "GL(2,16).2 cap G" ]           1/108
       SL(6,2)   63                [ "GL(3,4).2", "GL(2,8).3" ]       365/55552
       SL(6,3)  364   [ "GL(3,9).2 cap G", "GL(2,27).3 cap G" ] 22843/123845436
       SL(6,4) 1365  [ "GL(3,16).2 cap G", "GL(2,64).3 cap G" ]         1/85932
       SL(6,5) 3906 [ "GL(3,25).2 cap G", "GL(2,125).3 cap G" ]        1/484220
       SL(8,2)  255                             [ "GL(4,4).2" ]          1/7874
      SL(10,2) 1023               [ "GL(5,4).2", "GL(2,32).5" ]        1/129794
\end{verbatim}


The only missing case for~\cite{BGK} is $S = L_3(4)$,
for which $\M(S,s)$ consists of three groups of the type $L_3(2)$
(see~\cite[p.~23]{CCN85}).
The group $L_3(4)$ has been considered already in Section~\ref{easyloop},
where $\total(S,s) = 1/5$ has been proved.
Also the cases $\SL(3,3)$, $\SL(4,2) \cong A_8$, and $\SL(4,3)$ have been
handled there.

An alternative character-theoretic proof for $S = L_6(2)$ looks as follows.
In this case, the subgroups in $\M(S,s)$ have the types
$\GammaL(3,4) \cong \GL(3,4).2 \cong 3.L_3(4).3.2_2$
and $\GammaL(2,8) \cong \GL(2,8).3 \cong (7 \times L_2(8)).3$.

\begin{verbatim}
    gap> t:= CharacterTable( "L6(2)" );;
    gap> s1:= CharacterTable( "3.L3(4).3.2_2" );;
    gap> s2:= CharacterTable( "(7xL2(8)).3" );;
    gap> SigmaFromMaxes( t, "63A", [ s1, s2 ], [ 1, 1 ] );
    365/55552
\end{verbatim}


\subsection{$\ast$~$L_d(q)$ with prime $d$}

For $S = \SL(d,q)$ with \emph{prime} dimension $d$,
and $s \in S$ a Singer cycle,
we have $\M(S,s) = \{ M \}$,
where $M = N_S(\langle s \rangle) \cong \GammaL(1,q^d) \cap S$.
So
\[
   \total(g,s) = \fpr(g,S/M) = |g^S \cap M|/|g^S|
               < |M|/|g^S| \leq (q^d-1) \cdot d/|g^S|
\]
holds for any $g \in S \setminus Z(S)$,
which implies
$\total( S, s ) < \max\{ (q^d-1) \cdot d/|g^S|; g \in S \setminus Z(S) \}$.
The right hand side of this inequality is returned by the following function.
In~\cite[Lemma~3.8]{BGK},
the global upper bound $1/q^d$ is derived for primes $d \geq 5$.

\begin{verbatim}
    gap> UpperBoundForSL:= function( d, q )
    >     local G, Msize, ccl;
    > 
    >     if not IsPrimeInt( d ) then
    >       Error( "<d> must be a prime" );
    >     fi;
    > 
    >     G:= SL( d, q );
    >     Msize:= (q^d-1) * d;
    >     ccl:= Filtered( ConjugacyClasses( G ),
    >                     c ->     Msize mod Order( Representative( c ) ) = 0
    >                          and Size( c ) <> 1 );
    > 
    >     return Msize / Minimum( List( ccl, Size ) );
    > end;;
\end{verbatim}

The interesting values are $(d,q)$ with $d \in \{ 5, 7, 11 \}$ and
$q \in \{ 2, 3, 4 \}$, and perhaps also $(d,q) \in \{ (3,2), (3,3) \}$.
(Here we exclude $\SL(11,4)$ because writing down the conjugacy classes
of this group would exceed the permitted memory.)

\begin{verbatim}
    gap> NrConjugacyClasses( SL(11,4) );
    1397660
    gap> pairs:= [ [ 3, 2 ], [ 3, 3 ], [ 5, 2 ], [ 5, 3 ], [ 5, 4 ],
    >              [ 7, 2 ], [ 7, 3 ], [ 7, 4 ],
    >              [ 11, 2 ], [ 11, 3 ] ];;
    gap> array:= [];;
    gap> for pair in pairs do
    >      d:= pair[1];  q:= pair[2];
    >      approx:= UpperBoundForSL( d, q );
    >      Add( array, [ Concatenation( "SL(", String(d), ",", String(q), ")" ),
    >                    (q^d-1)/(q-1),
    >                    approx ] );
    >    od;
    gap> PrintFormattedArray( array );
       SL(3,2)     7                                   7/8
       SL(3,3)    13                                   3/4
       SL(5,2)    31                              31/64512
       SL(5,3)   121                                 10/81
       SL(5,4)   341                                15/256
       SL(7,2)   127                             7/9142272
       SL(7,3)  1093                                14/729
       SL(7,4)  5461                               21/4096
      SL(11,2)  2047 2047/34112245508649716682268134604800
      SL(11,3) 88573                              22/59049
\end{verbatim}

The exact values are clearly better than the above bounds.
We compute them for $L_5(2)$ and $L_7(2)$.
In the latter case, the class fusion of the $127:7$ type subgroup $M$
is not uniquely determined by the character tables;
here we use the additional information that the elements of order $7$ in $M$
have centralizer order $49$ in $L_7(2)$.
(See Section~\ref{easyloop} for the examples with $d = 3$.)

\begin{verbatim}
    gap> SigmaFromMaxes( CharacterTable( "L5(2)" ), "31A",
    >        [ CharacterTable( "31:5" ) ], [ 1 ] );
    1/5376
    gap> t:= CharacterTable( "L7(2)" );;
    gap> s:= CharacterTable( "P:Q", [ 127, 7 ] );;
    gap> pi:= PossiblePermutationCharacters( s, t );;
    gap> Length( pi );
    2
    gap> ord7:= PositionsProperty( OrdersClassRepresentatives( t ), x -> x = 7 );
    [ 38, 45, 76, 77, 83 ]
    gap> sizes:= SizesCentralizers( t ){ ord7 };
    [ 141120, 141120, 3528, 3528, 49 ]
    gap> List( pi, x -> x[83] );
    [ 42, 0 ]
    gap> spos:= Position( OrdersClassRepresentatives( t ), 127 );;
    gap> Maximum( ApproxP( pi{ [ 1 ] }, spos ) );
    1/4388290560
\end{verbatim}


\subsection{Automorphic Extensions of $L_d(q)$}\label{SLaut}

For the following values of $d$ and $q$,
automorphic extensions $G$ of $L_d(q)$ had to be checked
for~\cite[Section~5.12]{BGK}.
\[
   (d,q) \in \{ (3,4), (6,2), (6,3), (6,4), (6,5), (10,2) \}
\]
The first case has been treated in Section~\ref{easyloopaut}.
For the other cases, we compute $\total^{\prime}(G,s)$ below.

In any case, the extension by a \emph{graph} automorphism occurs,
which can be described by mapping each matrix in $\SL(d,q)$ to its inverse
transpose.
If $q > 2$, also extensions by \emph{diagonal} automorphisms occur,
which are induced by conjugation with elements in $\GL(d,q)$.
If $q$ is nonprime then also extensions by \emph{field} automorphisms occur,
which can be described by powering the matrix entries by roots of $q$.
Finally, products (of prime order) of these three kinds of automorphisms
have to be considered.

We start with the extension $G$ of $S = \SL(d,q)$ by a graph automorphism.
$G$ can be embedded into $\GL(2d,q)$ by representing the matrix $A \in S$
as a block diagonal matrix with diagonal blocks equal to $A$ and $A^{-tr}$,
and representing the graph automorphism by a permutation matrix that
interchanges the two blocks.
In order to construct the field extension type subgroups of $G$,
we have to choose the basis of the field extension in such a way that the
subgroup is normalized by the permutation matrix;
a sufficient condition is that the matrices of the $\F_q$-linear mappings
induced by the basis elements are symmetric.

(We do not give a function that computes a basis with this property from
the parameters $d$ and $q$.
Instead, we only write down the bases that we will need.)

\begin{verbatim}
    gap> SymmetricBasis:= function( q, n )
    >     local vectors, B, issymmetric;
    > 
    >     if   q = 2 and n = 2 then
    >       vectors:= [ Z(2)^0, Z(2^2) ];
    >     elif q = 2 and n = 3 then
    >       vectors:= [ Z(2)^0, Z(2^3), Z(2^3)^5 ];
    >     elif q = 2 and n = 5 then
    >       vectors:= [ Z(2)^0, Z(2^5), Z(2^5)^4, Z(2^5)^25, Z(2^5)^26 ];
    >     elif q = 3 and n = 2 then
    >       vectors:= [ Z(3)^0, Z(3^2) ];
    >     elif q = 3 and n = 3 then
    >       vectors:= [ Z(3)^0, Z(3^3)^2, Z(3^3)^7 ];
    >     elif q = 4 and n = 2 then
    >       vectors:= [ Z(2)^0, Z(2^4)^3 ];
    >     elif q = 4 and n = 3 then
    >       vectors:= [ Z(2)^0, Z(2^3), Z(2^3)^5 ];
    >     elif q = 5 and n = 2 then
    >       vectors:= [ Z(5)^0, Z(5^2)^2 ];
    >     elif q = 5 and n = 3 then
    >       vectors:= [ Z(5)^0, Z(5^3)^9, Z(5^3)^27 ];
    >     else
    >       Error( "sorry, no basis for <q> and <n> stored" );
    >     fi;
    > 
    >     B:= Basis( AsField( GF(q), GF(q^n) ), vectors );
    > 
    >     # Check that the basis really has the required property.
    >     issymmetric:= M -> M = TransposedMat( M );
    >     if not ForAll( B, b -> issymmetric( BlownUpMat( B, [ [ b ] ] ) ) ) then
    >       Error( "wrong basis!" );
    >     fi;
    > 
    >     # Return the result.
    >     return B;
    > end;;
\end{verbatim}

%
%

In later examples, we will need similar embeddings of matrices.
Therefore, we provide a more general function \verb|EmbeddedMatrix|
that takes a field \verb|F|, a matrix \verb|mat|, and a function \verb|func|,
and returns a block diagonal matrix over \verb|F| whose diagonal blocks are
\verb|mat| and \verb|func( mat )|.

\begin{verbatim}
    gap> BindGlobal( "EmbeddedMatrix", function( F, mat, func )
    >   local d, result;
    > 
    >   d:= Length( mat );
    >   result:= NullMat( 2*d, 2*d, F );
    >   result{ [ 1 .. d ] }{ [ 1 .. d ] }:= mat;
    >   result{ [ d+1 .. 2*d ] }{ [ d+1 .. 2*d ] }:= func( mat );
    > 
    >   return result;
    > end );
\end{verbatim}

The following function is similar to \verb|ApproxPForSL|,
the differences are that the group $G$ in question is not $\SL(d,q)$ but
the extension of this group by a graph automorphism,
and that $\total^{\prime}(G,s)$ is computed not $\total(G,s)$.

\begin{verbatim}
    gap> ApproxPForOuterClassesInExtensionOfSLByGraphAut:= function( d, q )
    >     local embedG, swap, G, orb, epi, PG, Gprime, primes, maxes, ccl, names;
    > 
    >     # Check whether this is an admissible case (see [Be00],
    >     # note that a graph automorphism exists only for `d > 2').
    >     if d = 2 or ( d = 3 and q = 4 ) then
    >       return fail;
    >     fi;
    > 
    >     # Provide a function that constructs a block diagonal matrix.
    >     embedG:= mat -> EmbeddedMatrix( GF( q ), mat,
    >                                     M -> TransposedMat( M^-1 ) );
    > 
    >     # Create the matrix that exchanges the two blocks.
    >     swap:= NullMat( 2*d, 2*d, GF(q) );
    >     swap{ [ 1 .. d ] }{ [ d+1 .. 2*d ] }:= IdentityMat( d, GF(q) );
    >     swap{ [ d+1 .. 2*d ] }{ [ 1 .. d ] }:= IdentityMat( d, GF(q) );
    > 
    >     # Create the group SL(d,q).2, and the map to the projective group.
    >     G:= ClosureGroupDefault( Group( List( GeneratorsOfGroup( SL( d, q ) ),
    >                                           embedG ) ),
    >                       swap );
    >     orb:= Orbit( G, One( G )[1], OnLines );
    >     epi:= ActionHomomorphism( G, orb, OnLines );
    >     PG:= ImagesSource( epi );
    >     Gprime:= DerivedSubgroup( PG );
    > 
    >     # Create the subgroups corresponding to the prime divisors of `d'.
    >     primes:= Set( Factors( d ) );
    >     maxes:= List( primes,
    >               p -> ClosureGroupDefault( Group( List( GeneratorsOfGroup(
    >                          RelativeSigmaL( d/p, SymmetricBasis( q, p ) ) ),
    >                          embedG ) ),
    >                      swap ) );
    > 
    >     # Compute conjugacy classes of outer involutions in the maxes.
    >     # (In order to avoid computing all conjugacy classes of these subgroups,
    >     # we work in the Sylow $2$ subgroups.)
    >     maxes:= List( maxes, M -> ImagesSet( epi, M ) );
    >     ccl:= List( maxes, M -> ClassesOfPrimeOrder( M, [ 2 ], Gprime ) );
    >     names:= List( primes, p -> Concatenation( "GL(", String( d/p ), ",",
    >                                    String( q^p ), ").", String( p ) ) );
    > 
    >     return [ names, UpperBoundFixedPointRatios( PG, ccl, true )[1] ];
    > end;;
\end{verbatim}

And these are the results for the groups we are interested in
(and others).

\begin{verbatim}
    gap> ApproxPForOuterClassesInExtensionOfSLByGraphAut( 4, 3 );
    [ [ "GL(2,9).2" ], 17/117 ]
    gap> ApproxPForOuterClassesInExtensionOfSLByGraphAut( 4, 4 );
    [ [ "GL(2,16).2" ], 73/1008 ]
    gap> ApproxPForOuterClassesInExtensionOfSLByGraphAut( 6, 2 );
    [ [ "GL(3,4).2", "GL(2,8).3" ], 41/1984 ]
    gap> ApproxPForOuterClassesInExtensionOfSLByGraphAut( 6, 3 );
    [ [ "GL(3,9).2", "GL(2,27).3" ], 541/352836 ]
    gap> ApproxPForOuterClassesInExtensionOfSLByGraphAut( 6, 4 );
    [ [ "GL(3,16).2", "GL(2,64).3" ], 3265/12570624 ]
    gap> ApproxPForOuterClassesInExtensionOfSLByGraphAut( 6, 5 );
    [ [ "GL(3,25).2", "GL(2,125).3" ], 13001/195250000 ]
    gap> ApproxPForOuterClassesInExtensionOfSLByGraphAut( 8, 2 );
    [ [ "GL(4,4).2" ], 367/1007872 ]
    gap> ApproxPForOuterClassesInExtensionOfSLByGraphAut( 10, 2 );
    [ [ "GL(5,4).2", "GL(2,32).5" ], 609281/476346056704 ]
\end{verbatim}

Now we consider diagonal automorphisms.
We modify the approach for $\SL(d,q)$ by constructing the field extension
type subgroups of $\GL(d,q)$ \ldots

\begin{verbatim}
    gap> RelativeGammaL:= function( d, B )
    >     local n, F, q, diag;
    > 
    >     n:= Length( B );
    >     F:= LeftActingDomain( UnderlyingLeftModule( B ) );
    >     q:= Size( F );
    >     diag:= IdentityMat( d * n, F );
    >     diag{[ 1 .. n ]}{[ 1 .. n ]}:= BlownUpMat( B, [ [ Z(q^n) ] ] );
    >     return ClosureGroup( RelativeSigmaL( d, B ),  diag );
    > end;;
\end{verbatim}

\ldots and counting the elements of prime order outside the simple group.

\begin{verbatim}
    gap> ApproxPForOuterClassesInGL:= function( d, q )
    >     local G, epi, PG, Gprime, primes, maxes, names;
    > 
    >     # Check whether this is an admissible case (see [Be00]).
    >     if ( d = 2 and q in [ 2, 5, 7, 9 ] ) or ( d = 3 and q = 4 ) then
    >       return fail;
    >     fi;
    > 
    >     # Create the group GL(d,q), and the map to PGL(d,q).
    >     G:= GL( d, q );
    >     epi:= ActionHomomorphism( G, NormedRowVectors( GF(q)^d ), OnLines );
    >     PG:= ImagesSource( epi );
    >     Gprime:= ImagesSet( epi, SL( d, q ) );
    > 
    >     # Create the subgroups corresponding to the prime divisors of `d'.
    >     primes:= Set( Factors( d ) );
    >     maxes:= List( primes, p -> RelativeGammaL( d/p,
    >                                    Basis( AsField( GF(q), GF(q^p) ) ) ) );
    >     maxes:= List( maxes, M -> ImagesSet( epi, M ) );
    >     names:= List( primes, p -> Concatenation( "M(", String( d/p ), ",",
    >                                    String( q^p ), ")" ) );
    > 
    >     return [ names,
    >              UpperBoundFixedPointRatios( PG, List( maxes,
    >                  M -> ClassesOfPrimeOrder( M,
    >                           Set( Factors( Index( PG, Gprime ) ) ), Gprime ) ),
    >                  true )[1] ];
    > end;;
\end{verbatim}

Here are the required results.

\begin{verbatim}
    gap> ApproxPForOuterClassesInGL( 6, 3 );
    [ [ "M(3,9)", "M(2,27)" ], 41/882090 ]
    gap> ApproxPForOuterClassesInGL( 4, 3 );
    [ [ "M(2,9)" ], 0 ]
    gap> ApproxPForOuterClassesInGL( 6, 4 );
    [ [ "M(3,16)", "M(2,64)" ], 1/87296 ]
    gap> ApproxPForOuterClassesInGL( 6, 5 );
    [ [ "M(3,25)", "M(2,125)" ], 821563/756593750000 ]
\end{verbatim}

(Note that the extension field type subgroup in $\PGL(4,3) = L_4(3).2_1$
is a \emph{non-split} extension of its intersection with $L_4(3)$,
hence the zero value.)

%
%

Concerning extensions by Frobenius automorphisms,
only the case $(d,q) = (6,4)$ is interesting in~\cite{BGK}.
In fact, we would not need to compute anything for the extension $G$ of
$S = \SL(6,4)$ by the Frobenius map that squares each matrix entry.
This is because $\M^{\prime}(G,s)$ consists of the normalizers
of the two subgroups of the types $\SL(3,16)$ and $\SL(2,64)$,
and the former maximal subgroup is a \emph{non-split} extension of its
intersection with $S$,
so only one maximal subgroup can contribute to $\total^{\prime}(G,s)$,
which is thus smaller than $1/2$, by~\cite[Prop.~2.6]{BGK}.

However, it is easy enough to compute the exact value of
$\total^{\prime}(G,s)$.
We work with the projective action of $S$ on its natural module,
and compute the permutation induced by the Frobenius map as the
Frobenius action on the normed row vectors.

\begin{verbatim}
    gap> matgrp:= SL(6,4);;
    gap> dom:= NormedRowVectors( GF(4)^6 );;
    gap> Gprime:= Action( matgrp, dom, OnLines );;
    gap> pi:= PermList( List( dom, v -> Position( dom, List( v, x -> x^2 ) ) ) );;
    gap> G:= ClosureGroup( Gprime, pi );;
\end{verbatim}

Then we compute the maximal subgroups, the classes of outer involutions,
and the bound, similar to the situation with graph automorphisms.


\begin{verbatim}
    gap> maxes:= List( [ 2, 3 ], p -> Normalizer( G,
    >              Action( RelativeSigmaL( 6/p,
    >                Basis( AsField( GF(4), GF(4^p) ) ) ), dom, OnLines ) ) );;
    gap> ccl:= List( maxes, M -> ClassesOfPrimeOrder( M, [ 2 ], Gprime ) );;
    gap> List( ccl, Length );
    [ 0, 1 ]
    gap> UpperBoundFixedPointRatios( G, ccl, true );
    [ 1/34467840, true ]
\end{verbatim}

For $(d,q) = (6,4)$,
we have to consider also the extension $G$ of $S = \SL(6,4)$
by the product $\alpha$ of the Frobenius map and the graph automorphism.
We use the same approach as for the graph automorphism,
i.~e., we embed $\SL(6,4)$ into a $12$-dimensional group of $6 \times 6$
block matrices,
where the second block is the image of the first block under $\alpha$,
and describe $\alpha$ by the transposition of the two blocks.

First we construct the projective actions of $S$ and $G$ on an orbit of
$1$-spaces.

\begin{verbatim}
    gap> embedFG:= function( F, mat )
    >      return EmbeddedMatrix( F, mat,
    >                 M -> List( TransposedMat( M^-1 ),
    >                            row -> List( row, x -> x^2 ) ) );
    >    end;;
    gap> d:= 6;;  q:= 4;;
    gap> alpha:= NullMat( 2*d, 2*d, GF(q) );;
    gap> alpha{ [ 1 .. d ] }{ [ d+1 .. 2*d ] }:= IdentityMat( d, GF(q) );;
    gap> alpha{ [ d+1 .. 2*d ] }{ [ 1 .. d ] }:= IdentityMat( d, GF(q) );;
    gap> Gprime:= Group( List( GeneratorsOfGroup( SL(d,q) ),
    >                          mat -> embedFG( GF(q), mat ) ) );;
    gap> G:= ClosureGroupDefault( Gprime, alpha );;
    gap> orb:= Orbit( G, One( G )[1], OnLines );;
    gap> G:= Action( G, orb, OnLines );;
    gap> Gprime:= Action( Gprime, orb, OnLines );;
\end{verbatim}

Next we construct the maximal subgroups, the classes of outer involutions,
and the bound.

\begin{verbatim}
    gap> maxes:= List( Set( Factors( d ) ), p -> Group( List( GeneratorsOfGroup(
    >              RelativeSigmaL( d/p, Basis( AsField( GF(q), GF(q^p) ) ) ) ),
    >                mat -> embedFG( GF(q), mat ) ) ) );;
    gap> maxes:= List( maxes, x -> Action( x, orb, OnLines ) );;
    gap> maxes:= List( maxes, x -> Normalizer( G, x ) );;
    gap> ccl:= List( maxes, M -> ClassesOfPrimeOrder( M, [ 2 ], Gprime ) );;
    gap> List( ccl, Length );
    [ 0, 1 ]
    gap> UpperBoundFixedPointRatios( G, ccl, true );
    [ 1/10792960, true ]
\end{verbatim}

The only missing cases are the extensions of $\SL(6,3)$ and $\SL(6,5)$
by the involutory outer automorphism that acts as the product of a diagonal
and a graph automorphism.

In the case $S = \SL(6,3)$, we can directly write down the extension $G$.

\begin{verbatim}
    gap> d:= 6;;  q:= 3;;
    gap> diag:= IdentityMat( d, GF(q) );;
    gap> diag[1][1]:= Z(q);;
    gap> embedDG:= mat -> EmbeddedMatrix( GF(q), mat,
    >                                     M -> TransposedMat( M^-1 )^diag );;
    gap> Gprime:= Group( List( GeneratorsOfGroup( SL(d,q) ), embedDG ) );;
    gap> alpha:= NullMat( 2*d, 2*d, GF(q) );;
    gap> alpha{ [ 1 .. d ] }{ [ d+1 .. 2*d ] }:= IdentityMat( d, GF(q) );;
    gap> alpha{ [ d+1 .. 2*d ] }{ [ 1 .. d ] }:= IdentityMat( d, GF(q) );;
    gap> G:= ClosureGroupDefault( Gprime, alpha );;
\end{verbatim}

The maximal subgroups are constructed as the normalizers in $G$ of the
extension field type subgroups in $S$.
We work with a permutation representation of $G$.

\begin{verbatim}
    gap> maxes:= List( Set( Factors( d ) ), p -> Group( List( GeneratorsOfGroup(
    >              RelativeSigmaL( d/p, Basis( AsField( GF(q), GF(q^p) ) ) ) ),
    >                embedDG ) ) );;
    gap> orb:= Orbit( G, One( G )[1], OnLines );;
    gap> G:= Action( G, orb, OnLines );;
    gap> Gprime:= Action( Gprime, orb, OnLines );;
    gap> maxes:= List( maxes, M -> Normalizer( G, Action( M, orb, OnLines ) ) );;
    gap> ccl:= List( maxes, M -> ClassesOfPrimeOrder( M, [ 2 ], Gprime ) );;
    gap> List( ccl, Length );
    [ 1, 1 ]
    gap> UpperBoundFixedPointRatios( G, ccl, true );
    [ 25/352836, true ]
\end{verbatim}

For $S = \SL(6,5)$, this approach does not work because we cannot
realize the diagonal involution by an involutory matrix.
Instead, we consider the extension of $\GL(6,5) \cong 2.(2 \times L_6(5)).2$
by the graph automorphism $\alpha$, which can be embedded into $\GL(12,5)$.

\begin{verbatim}
    gap> d:= 6;;  q:= 5;;
    gap> embedG:= mat -> EmbeddedMatrix( GF(q),
    >                                    mat, M -> TransposedMat( M^-1 ) );;
    gap> Gprime:= Group( List( GeneratorsOfGroup( SL(d,q) ), embedG ) );;
    gap> maxes:= List( Set( Factors( d ) ), p -> Group( List( GeneratorsOfGroup(
    >              RelativeSigmaL( d/p, Basis( AsField( GF(q), GF(q^p) ) ) ) ),
    >                embedG ) ) );;
    gap> diag:= IdentityMat( d, GF(q) );;
    gap> diag[1][1]:= Z(q);;
    gap> diag:= embedG( diag );;
    gap> alpha:= NullMat( 2*d, 2*d, GF(q) );;
    gap> alpha{ [ 1 .. d ] }{ [ d+1 .. 2*d ] }:= IdentityMat( d, GF(q) );;
    gap> alpha{ [ d+1 .. 2*d ] }{ [ 1 .. d ] }:= IdentityMat( d, GF(q) );;
    gap> G:= ClosureGroupDefault( Gprime, alpha * diag );;
\end{verbatim}

Now we switch to the permutation action of this group on the $1$-dimensional
subspaces, thus factoring out the cyclic normal subgroup of order four.
In this action, the involutory diagonal automorphism is represented by an
involution, and we can proceed as above.

\begin{verbatim}
    gap> orb:= Orbit( G, One( G )[1], OnLines );;
    gap> Gprime:= Action( Gprime, orb, OnLines );;
    gap> G:= Action( G, orb, OnLines );;
    gap> maxes:= List( maxes, M -> Action( M, orb, OnLines ) );;
    gap> extmaxes:= List( maxes, M -> Normalizer( G, M ) );;
    gap> ccl:= List( extmaxes, M -> ClassesOfPrimeOrder( M, [ 2 ], Gprime ) );;
    gap> List( ccl, Length );
    [ 2, 1 ]
    gap> UpperBoundFixedPointRatios( G, ccl, true );
    [ 3863/6052750000, true ]
\end{verbatim}

In the same way, we can recheck the values for the extensions of $\SL(6,5)$
by the diagonal or by the graph automorphism.

\begin{verbatim}
    gap> diag:= Permutation( diag, orb, OnLines );;
    gap> G:= ClosureGroupDefault( Gprime, diag );;
    gap> extmaxes:= List( maxes, M -> Normalizer( G, M ) );;
    gap> ccl:= List( extmaxes, M -> ClassesOfPrimeOrder( M, [ 2 ], Gprime ) );;
    gap> List( ccl, Length );
    [ 3, 1 ]
    gap> UpperBoundFixedPointRatios( G, ccl, true );
    [ 821563/756593750000, true ]
    gap> alpha:= Permutation( alpha, orb, OnLines );;
    gap> G:= ClosureGroupDefault( Gprime, alpha );;
    gap> extmaxes:= List( maxes, M -> Normalizer( G, M ) );;
    gap> ccl:= List( extmaxes, M -> ClassesOfPrimeOrder( M, [ 2 ], Gprime ) );;
    gap> List( ccl, Length );
    [ 2, 2 ]
    gap> UpperBoundFixedPointRatios( G, ccl, true );
    [ 13001/195250000, true ]
\end{verbatim}

\subsection{$L_3(2)$}\label{L32}

We show that $S = L_3(2) = \SL(3,2)$ satisfies the following.
\begin{enumerate}
\item[(a)]
    $\total(S) = 1/4$,
    and this value is attained exactly for $\total(S,s)$
    with $s$ of order $7$.
\item[(b)]
    For $s$ of order $7$,
    $\M(S,s)$ consists of one group of the type $7:3$.
\item[(c)]
    $\prop(S) = 1/4$,
    and this value is attained exactly for $\prop(S,s)$
    with $s$ of order $7$.
\item[(d)]
    The uniform spread of $S$ is at exactly three,
    with $s$ of order $7$,
    and the spread of $S$ is exactly four.
    (This had been left open in~\cite{BW1}.)
\end{enumerate}

(Note that in this example, the spread and the uniform spread differ.)

Statement~(a) follows from inspection of the primitive permutation
characters, cf.~Section~\ref{easyloop}.

\begin{verbatim}
    gap> t:= CharacterTable( "L3(2)" );;
    gap> ProbGenInfoSimple( t );
    [ "L3(2)", 1/4, 3, [ "7A" ], [ 1 ] ]
\end{verbatim}

Statement~(b) can be read off from the permutation characters,
and the fact that the unique class of maximal subgroups that contain
elements of order $7$ consists of groups of the structure $7:3$,
see~\cite[p.~3]{CCN85}.

\begin{verbatim}
    gap> OrdersClassRepresentatives( t );
    [ 1, 2, 3, 4, 7, 7 ]
    gap> PrimitivePermutationCharacters( t );
    [ Character( CharacterTable( "L3(2)" ), [ 7, 3, 1, 1, 0, 0 ] ), 
      Character( CharacterTable( "L3(2)" ), [ 7, 3, 1, 1, 0, 0 ] ), 
      Character( CharacterTable( "L3(2)" ), [ 8, 0, 2, 0, 1, 1 ] ) ]
\end{verbatim}

For the other statements, we will use the primitive permutation
representations on $7$ and $8$ points of $S$
(computed from the {\GAP} Library of Tables of Marks),
and their diagonal products of the degrees $14$ and $15$.

\begin{verbatim}
    gap> tom:= TableOfMarks( "L3(2)" );;
    gap> g:= UnderlyingGroup( tom );
    Group([ (2,4)(5,7), (1,2,3)(4,5,6) ])
    gap> mx:= MaximalSubgroupsTom( tom );
    [ [ 14, 13, 12 ], [ 7, 7, 8 ] ]
    gap> maxes:= List( mx[1], i -> RepresentativeTom( tom, i ) );;
    gap> tr:= List( maxes, s -> RightTransversal( g, s ) );;
    gap> acts:= List( tr, x -> Action( g, x, OnRight ) );;
    gap> g7:= acts[1];
    Group([ (3,4)(6,7), (1,3,2)(4,6,5) ])
    gap> g8:= acts[3];
    Group([ (1,6)(2,5)(3,8)(4,7), (1,7,3)(2,5,8) ])
    gap> g14:= DiagonalProductOfPermGroups( acts{ [ 1, 2 ] } );
    Group([ (3,4)(6,7)(11,13)(12,14), (1,3,2)(4,6,5)(8,11,9)(10,12,13) ])
    gap> g15:= DiagonalProductOfPermGroups( acts{ [ 2, 3 ] } );
    Group([ (4,6)(5,7)(8,13)(9,12)(10,15)(11,14), 
      (1,4,2)(3,5,6)(8,14,10)(9,12,15) ])
\end{verbatim}

First we compute that for all nonidentity elements $s \in S$
and order three elements $g \in S$,
$\prop(g,s) \geq 1/4$ holds,
with equality if and only if $s$ has order $7$;
this implies statement~(c).
We actually compute, for class representatives $s$,
the proportion of order three elements $g$ such that
$\langle g, s \rangle \not= S$ holds.

\begin{verbatim}
    gap> ccl:= List( ConjugacyClasses( g7 ), Representative );;
    gap> SortParallel( List( ccl, Order ), ccl );
    gap> List( ccl, Order );
    [ 1, 2, 3, 4, 7, 7 ]
    gap> Size( ConjugacyClass( g7, ccl[3] ) );
    56
    gap> prop:= List( ccl,
    >                 r -> RatioOfNongenerationTransPermGroup( g7, ccl[3], r ) );
    [ 1, 5/7, 19/28, 2/7, 1/4, 1/4 ]
    gap> Minimum( prop );
    1/4
\end{verbatim}

Now we show that the uniform spread of $S$ is less than four.
In any of the primitive permutation representations of degree seven,
we find three involutions whose sets of fixed points cover the
seven points.
The elements $s$ of order different from $7$ in $S$ fix a point in this
representation, so each such $s$ generates a proper subgroup of $S$
together with one of the three involutions.

\begin{verbatim}
    gap> x:= g7.1;
    (3,4)(6,7)
    gap> fix:= Difference( MovedPoints( g7 ), MovedPoints( x ) );
    [ 1, 2, 5 ]
    gap> orb:= Orbit( g7, fix, OnSets );
    [ [ 1, 2, 5 ], [ 1, 3, 4 ], [ 2, 3, 6 ], [ 2, 4, 7 ], [ 1, 6, 7 ], 
      [ 3, 5, 7 ], [ 4, 5, 6 ] ]
    gap> Union( orb{ [ 1, 2, 5 ] } ) = [ 1 .. 7 ];
    true
\end{verbatim}

So we still have to exclude elements $s$ of order $7$.
In the primitive permutation representation of $S$ on eight points,
we find four elements of order three whose sets of fixed points
cover the set of all points that are moved by $S$,
so with each element of order seven in $S$,
one of them generates an intransitive group.

\begin{verbatim}
    gap> three:= g8.2;
    (1,7,3)(2,5,8)
    gap> fix:= Difference( MovedPoints( g8 ), MovedPoints( three ) );
    [ 4, 6 ]
    gap> orb:= Orbit( g8, fix, OnSets );;
    gap> QuadrupleWithProperty( [ [ fix ], orb, orb, orb ],
    >        list -> Union( list ) = [ 1 .. 8 ] );
    [ [ 4, 6 ], [ 1, 7 ], [ 3, 8 ], [ 2, 5 ] ]
\end{verbatim}

Together with statement~(a), this proves that the uniform spread of $S$
is exactly three, with $s$ of order seven.

Each element of $S$ fixes a point
in the permutation representation on $15$ points.
So for proving that the spread of $S$ is less than five,
it is sufficient to find a quintuple of elements whose sets of fixed points
cover all $15$ points.
(From the permutation characters it is clear that four of these elements
must have order three, and the fifth must be an involution.)

\begin{verbatim}
    gap> x:= g15.1;
    (4,6)(5,7)(8,13)(9,12)(10,15)(11,14)
    gap> fixx:= Difference( MovedPoints( g15 ), MovedPoints( x ) );
    [ 1, 2, 3 ]
    gap> orbx:= Orbit( g15, fixx, OnSets );
    [ [ 1, 2, 3 ], [ 1, 4, 5 ], [ 1, 6, 7 ], [ 2, 4, 6 ], [ 3, 4, 7 ], 
      [ 3, 5, 6 ], [ 2, 5, 7 ] ]
    gap> y:= g15.2;
    (1,4,2)(3,5,6)(8,14,10)(9,12,15)
    gap> fixy:= Difference( MovedPoints( g15 ), MovedPoints( y ) );
    [ 7, 11, 13 ]
    gap> orby:= Orbit( g15, fixy, OnSets );;
    gap> QuadrupleWithProperty( [ [ fixy ], orby, orby, orby ],
    >        l -> Difference( [ 1 .. 15 ], Union( l ) ) in orbx );
    [ [ 7, 11, 13 ], [ 5, 8, 14 ], [ 1, 10, 15 ], [ 3, 9, 12 ] ]
\end{verbatim}

It remains to show that the spread of $S$ is (at least) four.
By the consideration of permutation characters,
we know that we can find a suitable order seven element
for all quadruples in question
except perhaps quadruples of order three elements.
We show that for each such case, we can choose $s$ of order four.
Since $\M(S,s)$ consists of two subgroups of the type $S_4$,
we work with the representation on $14$ points.)

First we compute $s$ and the $S$-orbit of its fixed points,
and the $S$-orbit of the fixed points of an element $x$ of order three.
Then we prove that for each quadruple of conjugates of $x$,
the union of their fixed points intersects the fixed points of
at least one conjugate of $s$ trivially.

\begin{verbatim}
    gap> ResetGlobalRandomNumberGenerators();
    gap> repeat s:= Random( g14 );
    >    until Order( s ) = 4;
    gap> s;
    (2,3,5,4)(6,7)(8,9)(11,12,14,13)
    gap> fixs:= Difference( MovedPoints( g14 ), MovedPoints( s ) );
    [ 1, 10 ]
    gap> orbs:= Orbit( g14, fixs, OnSets );;
    gap> Length( orbs );
    21
    gap> three:= g14.2;
    (1,3,2)(4,6,5)(8,11,9)(10,12,13)
    gap> fix:= Difference( MovedPoints( g14 ), MovedPoints( three ) );
    [ 7, 14 ]
    gap> orb:= Orbit( g14, fix, OnSets );;
    gap> Length( orb );
    28
    gap> QuadrupleWithProperty( [ [ fix ], orb, orb, orb ],
    >        l -> ForAll( orbs, o -> not IsEmpty( Intersection( o,
    >                        Union( l ) ) ) ) );
    fail
\end{verbatim}

By Lemma~\ref{existsgoodconjugate}, we are done.

\subsection{$M_{11}$}\label{spreadM11}

We show that $S = M_{11}$ satisfies the following.
\begin{enumerate}
\item[(a)]
    $\total(S) = 1/3$,
    and this value is attained exactly for $\total(S,s)$
    with $s$ of order $11$.
\item[(b)]
    For $s$ of order $11$,
    $\M(S,s)$ consists of one group of the type $L_2(11)$.
\item[(c)]
    $\prop(S) = 1/3$,
    and this value is attained exactly for $\prop(S,s)$
    with $s$ of order $11$.
\item[(d)]
    Both the uniform spread and the spread of $S$ is exactly three,
    with $s$ of order $11$.
\end{enumerate}

Statement~(a) follows from inspection of the primitive permutation
characters, cf.~Section~\ref{spor}.

\begin{verbatim}
    gap> t:= CharacterTable( "M11" );;
    gap> ProbGenInfoSimple( t );
    [ "M11", 1/3, 2, [ "11A" ], [ 1 ] ]
\end{verbatim}

Statement~(b) can be read off from the permutation characters,
and the fact that the unique class of maximal subgroups that contain
elements of order $11$ consists of groups of the structure $L_2(11)$,
see~\cite[p.~18]{CCN85}.

\begin{verbatim}
    gap> OrdersClassRepresentatives( t );
    [ 1, 2, 3, 4, 5, 6, 8, 8, 11, 11 ]
    gap> PrimitivePermutationCharacters( t );
    [ Character( CharacterTable( "M11" ), [ 11, 3, 2, 3, 1, 0, 1, 1, 0, 0 ] ),
      Character( CharacterTable( "M11" ), [ 12, 4, 3, 0, 2, 1, 0, 0, 1, 1 ] ),
      Character( CharacterTable( "M11" ), [ 55, 7, 1, 3, 0, 1, 1, 1, 0, 0 ] ),
      Character( CharacterTable( "M11" ), [ 66, 10, 3, 2, 1, 1, 0, 0, 0, 0 ] ),
      Character( CharacterTable( "M11" ), [ 165, 13, 3, 1, 0, 1, 1, 1, 0, 0 ] ) ]
    gap> Maxes( t );
    [ "A6.2_3", "L2(11)", "3^2:Q8.2", "A5.2", "2.S4" ]
\end{verbatim}

For the other statements, we will use the primitive permutation
representations of $S$ on $11$ and $12$ points
(which are fetched from the {\ATLAS} of Group Representations~\cite{AGR}),
and their diagonal product.

\begin{verbatim}
    gap> gens11:= OneAtlasGeneratingSet( "M11", NrMovedPoints, 11 );
    rec( charactername := "1a+10a",
      generators := [ (2,10)(4,11)(5,7)(8,9), (1,4,3,8)(2,5,6,9) ],
      groupname := "M11", id := "",
      identifier := [ "M11", [ "M11G1-p11B0.m1", "M11G1-p11B0.m2" ], 1, 11 ],
      isPrimitive := true, maxnr := 1, p := 11, rankAction := 2,
      repname := "M11G1-p11B0", repnr := 1, size := 7920, stabilizer := "A6.2_3",
      standardization := 1, transitivity := 4, type := "perm" )
    gap> g11:= GroupWithGenerators( gens11.generators );;
    gap> gens12:= OneAtlasGeneratingSet( "M11", NrMovedPoints, 12 );;
    gap> g12:= GroupWithGenerators( gens12.generators );;
    gap> g23:= DiagonalProductOfPermGroups( [ g11, g12 ] );
    Group([ (2,10)(4,11)(5,7)(8,9)(12,17)(13,20)(16,18)(19,21), 
      (1,4,3,8)(2,5,6,9)(12,17,18,15)(13,19)(14,20)(16,22,23,21) ])
\end{verbatim}

First we compute that for all nonidentity elements $s \in S$
and involutions $g \in S$,
$\prop(g,s) \geq 1/3$ holds,
with equality if and only if $s$ has order $11$;
this implies statement~(c).
We actually compute, for class representatives $s$,
the proportion of involutions $g$ such that
$\langle g, s \rangle \not= S$ holds.

\begin{verbatim}
    gap> inv:= g11.1;
    (2,10)(4,11)(5,7)(8,9)
    gap> ccl:= List( ConjugacyClasses( g11 ), Representative );;
    gap> SortParallel( List( ccl, Order ), ccl );
    gap> List( ccl, Order );
    [ 1, 2, 3, 4, 5, 6, 8, 8, 11, 11 ]
    gap> Size( ConjugacyClass( g11, inv ) );
    165
    gap> prop:= List( ccl,
    >                 r -> RatioOfNongenerationTransPermGroup( g11, inv, r ) );
    [ 1, 1, 1, 149/165, 25/33, 31/55, 23/55, 23/55, 1/3, 1/3 ]
    gap> Minimum( prop );
    1/3
\end{verbatim}

For the first part of statement~(d),
we have to deal only with the case of triples of involutions.

The $11$-cycle $s$ is contained in exactly one maximal subgroup
of $S$, of index $12$.
By Corollary~\ref{existsgoodconjugate1},
it is enough to show that in the primitive degree $12$ representation of $S$,
the fixed points of no triple $(x_1, x_2, x_3)$ of involutions in $S$
can cover all twelve points;
equivalenly (considering complements),
we show that there is no triple such that
the intersection of the sets of \emph{moved} points is empty.

\begin{verbatim}
    gap> inv:= g12.1;
    (1,6)(2,9)(5,7)(8,10)
    gap> moved:= MovedPoints( inv );
    [ 1, 2, 5, 6, 7, 8, 9, 10 ]
    gap> orb12:= Orbit( g12, moved, OnSets );;
    gap> Length( orb12 );
    165
    gap> TripleWithProperty( [ orb12{[1]}, orb12, orb12 ],
    >        list -> IsEmpty( Intersection( list ) ) );
    fail
\end{verbatim}

This implies that the uniform spread of $S$ is at least three.

Now we show that there is a quadruple consisting of one
element of order three and three involutions whose fixed points cover
all points in the degree $23$ representation constructed above;
since the permutation character of this representation is strictly positive,
this implies that $S$ does not have spread four,
by Corollary~\ref{disprovespread},
and we have proved statement~(d).

\begin{verbatim}
    gap> inv:= g23.1;
    (2,10)(4,11)(5,7)(8,9)(12,17)(13,20)(16,18)(19,21)
    gap> moved:= MovedPoints( inv );
    [ 2, 4, 5, 7, 8, 9, 10, 11, 12, 13, 16, 17, 18, 19, 20, 21 ]
    gap> orb23:= Orbit( g23, moved, OnSets );;
    gap> three:= ( g23.1*g23.2^2 )^2;
    (2,6,10)(4,8,7)(5,9,11)(12,17,23)(15,18,16)(19,21,22)
    gap> movedthree:= MovedPoints( three );;
    gap> QuadrupleWithProperty( [ [ movedthree ], orb23, orb23, orb23 ],
    >        list -> IsEmpty( Intersection( list ) ) );
    [ [ 2, 4, 5, 6, 7, 8, 9, 10, 11, 12, 15, 16, 17, 18, 19, 21, 22, 23 ], 
      [ 1, 3, 4, 5, 6, 8, 9, 10, 12, 13, 14, 16, 17, 18, 20, 21 ], 
      [ 1, 2, 3, 4, 5, 6, 7, 11, 12, 13, 14, 15, 18, 19, 20, 23 ], 
      [ 1, 2, 3, 7, 8, 9, 10, 11, 13, 14, 15, 16, 17, 20, 22, 23 ] ]
\end{verbatim}


\subsection{$M_{12}$}\label{spreadM12}

We show that $S = M_{12}$ satisfies the following.
\begin{enumerate}
\item[(a)]
    $\total(S) = 1/3$,
    and this value is attained exactly for $\total(S,s)$
    with $s$ of order $10$.
\item[(b)]
    For $s \in S$ of order $10$,
    $\M(S,s)$ consists of two nonconjugate subgroups of the type $A_6.2^2$,
    and one group of the type $2 \times S_5$.
\item[(c)]
    $\prop(S) = 31/99$,
    and this value is attained exactly for $\prop(S,s)$
    with $s$ of order $10$.
\item[(d)]
    The uniform spread of $S$ is at least three,
    with $s$ of order $10$.
\item[(e)]
    $\total^{\prime}(\Aut(S), s) = 4/99$.
\end{enumerate}

Statement~(a) follows from inspection of the primitive permutation
characters, cf.~Section~\ref{spor}.

\begin{verbatim}
    gap> t:= CharacterTable( "M12" );;
    gap> ProbGenInfoSimple( t );
    [ "M12", 1/3, 2, [ "10A" ], [ 3 ] ]
\end{verbatim}

Statement~(b) can be read off from the permutation characters,
and the fact that the only classes of maximal subgroups that contain
elements of order $10$ consist of groups of the structures
$A_6.2^2$ (two classes) and $2 \times S_5$ (one class),
see~\cite[p.~33]{CCN85}.

\begin{verbatim}
    gap> spos:= Position( OrdersClassRepresentatives( t ), 10 );
    13
    gap> prim:= PrimitivePermutationCharacters( t );;
    gap> List( prim, x -> x{ [ 1, spos ] } );
    [ [ 12, 0 ], [ 12, 0 ], [ 66, 1 ], [ 66, 1 ], [ 144, 0 ], [ 220, 0 ], 
      [ 220, 0 ], [ 396, 1 ], [ 495, 0 ], [ 495, 0 ], [ 1320, 0 ] ]
    gap> Maxes( t );
    [ "M11", "M12M2", "A6.2^2", "M12M4", "L2(11)", "3^2.2.S4", "M12M7", "2xS5", 
      "M8.S4", "4^2:D12", "A4xS3" ]
\end{verbatim}

For statement~(c) (which implies statement~(d)),
we use the primitive permutation representation on $12$ points.

\begin{verbatim}
    gap> g:= MathieuGroup( 12 );
    Group([ (1,2,3,4,5,6,7,8,9,10,11), (3,7,11,8)(4,10,5,6), 
      (1,12)(2,11)(3,6)(4,8)(5,9)(7,10) ])
\end{verbatim}

First we show that for $s$ of order $10$, $\prop(S,s) = 31/99$ holds.

\begin{verbatim}
    gap> approx:= ApproxP( prim, spos );
    [ 0, 3/11, 1/3, 1/11, 1/132, 13/99, 13/99, 13/396, 1/132, 1/33, 1/33, 1/33, 
      13/396, 0, 0 ]
    gap> 2B:= g.2^2;
    (3,11)(4,5)(6,10)(7,8)
    gap> Size( ConjugacyClass( g, 2B ) );
    495
    gap> ResetGlobalRandomNumberGenerators();
    gap> repeat s:= Random( g );
    >    until Order( s ) = 10;
    gap> prop:= RatioOfNongenerationTransPermGroup( g, 2B, s );
    31/99
    gap> Filtered( approx, x -> x >= prop );
    [ 1/3 ]
\end{verbatim}

Next we show that for $s$ of order different from $10$,
$\prop(g,s)$ is larger than $31/99$ for suitable $g \in S^{\times}$.
Except for $s$ in the class {\tt 6A} (which fixes no point in the degree $12$
representation), it suffices to consider $g$ in the class {\tt 2B}
(with four fixed points).

\begin{verbatim}
    gap> x:= g.2^2;
    (3,11)(4,5)(6,10)(7,8)
    gap> ccl:= List( ConjugacyClasses( g ), Representative );;
    gap> SortParallel( List( ccl, Order ), ccl );
    gap> prop:= List( ccl, r -> RatioOfNongenerationTransPermGroup( g, x, r ) );
    [ 1, 1, 1, 1, 39/55, 1, 1, 29/33, 7/55, 43/55, 383/495, 383/495, 31/99, 5/9,
      5/9 ]
    gap> bad:= Filtered( prop, x -> x < 31/99 );
    [ 7/55 ]
    gap> pos:= Position( prop, bad[1] );;
    gap> [ Order( ccl[ pos ] ), NrMovedPoints( ccl[ pos ] ) ];
    [ 6, 12 ]
\end{verbatim}

In the remaining case, we choose $g$ in the class {\tt 2A}
(which is fixed point free).

\begin{verbatim}
    gap> x:= g.3;
    (1,12)(2,11)(3,6)(4,8)(5,9)(7,10)
    gap> s:= ccl[ pos ];;
    gap> prop:= RatioOfNongenerationTransPermGroup( g, x, s );
    17/33
    gap> prop > 31/99;
    true
\end{verbatim}

Statement~(e) has been shown already in Section~\ref{sporaut}.


\subsection{$O_7(3)$}\label{O73}

We show that $S = O_7(3)$ satisfies the following.
\begin{enumerate}
\item[(a)]
    $\total(S) = 199/351$,
    and this value is attained exactly for $\total(S,s)$
    with $s$ of order $14$.
\item[(b)]
    For $s \in S$ of order $14$,
    $\M(S,s)$ consists of one group of the type
    $2.U_4(3).2_2 = \Omega^-(6,3).2$
    and two nonconjugate groups of the type $S_9$.
\item[(c)]
    $\prop(S) = 155/351$,
    and this value is attained exactly for $\prop(S,s)$
    with $s$ of order $14$.
\item[(d)]
    The uniform spread of $S$ is at least three,
    with $s$ of order $14$.
\item[(e)]
    $\total^{\prime}(\Aut(S), s) = 1/3$.
\end{enumerate}

Currently {\GAP} provides neither the table of marks of $S$ nor all
character tables of its maximal subgroups.
First we compute those primitive permutation characters of $S$ that have
the degrees
$351$ (point stabilizer $2.U_4(3).2_2$),
$364$ (point stabilizer $3^5:U_4(2).2$),
$378$ (point stabilizer $L_4(3).2_2$),
$1\,080$ (point stabilizer $G_2(3)$, two classes),
$1\,120$ (point stabilizer $3^{3+3}:L_3(3)$),
$3\,159$ (point stabilizer $S_6(2)$, two classes),
$12\,636$ (point stabilizer $S_9$, two classes),
$22\,113$ (point stabilizer $(2^2 \times U_4(2)).2$,
which extends to $D_8 \times U_4(2).2$ in $O_7(3).2$), and
$28\,431$ (point stabilizer $2^6:A_7$).

(So we ignore the primitive permutation characters of the degrees $3\,640$,
$265\,356$, and $331\,695$.
Note that the orders of the corresponding subgroups are not divisible by $7$.)

\begin{verbatim}
    gap> t:= CharacterTable( "O7(3)" );;
    gap> someprim:= [];;
    gap> pi:= PossiblePermutationCharacters(
    >             CharacterTable( "2.U4(3).2_2" ), t );;  Length( pi );
    1
    gap> Append( someprim, pi );
    gap> pi:= PermChars( t, rec( torso:= [ 364 ] ) );;  Length( pi );
    1
    gap> Append( someprim, pi );
    gap> pi:= PossiblePermutationCharacters(
    >             CharacterTable( "L4(3).2_2" ), t );;  Length( pi );
    1
    gap> Append( someprim, pi );
    gap> pi:= PossiblePermutationCharacters( CharacterTable( "G2(3)" ), t );
    [ Character( CharacterTable( "O7(3)" ), [ 1080, 0, 0, 24, 108, 0, 0, 0, 27, 
          18, 9, 0, 12, 4, 0, 0, 0, 0, 0, 0, 0, 0, 12, 0, 0, 0, 0, 0, 3, 6, 0, 3, 
          2, 2, 2, 0, 0, 0, 3, 0, 0, 0, 0, 0, 0, 4, 0, 3, 0, 1, 1, 0, 0, 0, 0, 0, 
          0, 0 ] ), Character( CharacterTable( "O7(3)" ), 
        [ 1080, 0, 0, 24, 108, 0, 0, 27, 0, 18, 9, 0, 12, 4, 0, 0, 0, 0, 0, 0, 0, 
          0, 12, 0, 0, 0, 0, 3, 0, 0, 6, 3, 2, 2, 2, 0, 0, 3, 0, 0, 0, 0, 0, 0, 
          0, 4, 3, 0, 0, 1, 1, 0, 0, 0, 0, 0, 0, 0 ] ) ]
    gap> Append( someprim, pi );
    gap> pi:= PermChars( t, rec( torso:= [ 1120 ] ) );;  Length( pi );
    1
    gap> Append( someprim, pi );
    gap> pi:= PossiblePermutationCharacters( CharacterTable( "S6(2)" ), t );
    [ Character( CharacterTable( "O7(3)" ), [ 3159, 567, 135, 39, 0, 81, 0, 0, 
          27, 27, 0, 15, 3, 3, 7, 4, 0, 27, 0, 0, 0, 0, 0, 9, 3, 0, 9, 0, 3, 9, 
          3, 0, 2, 1, 1, 0, 0, 0, 3, 0, 2, 0, 0, 0, 3, 0, 0, 3, 1, 0, 0, 0, 1, 0, 
          0, 0, 0, 0 ] ), Character( CharacterTable( "O7(3)" ), 
        [ 3159, 567, 135, 39, 0, 81, 0, 27, 0, 27, 0, 15, 3, 3, 7, 4, 0, 27, 0, 
          0, 0, 0, 0, 9, 3, 0, 9, 3, 0, 3, 9, 0, 2, 1, 1, 0, 0, 3, 0, 0, 2, 0, 0, 
          0, 3, 0, 3, 0, 1, 0, 0, 0, 1, 0, 0, 0, 0, 0 ] ) ]
    gap> Append( someprim, pi );
    gap> pi:= PossiblePermutationCharacters( CharacterTable( "S9" ), t );
    [ Character( CharacterTable( "O7(3)" ), [ 12636, 1296, 216, 84, 0, 81, 0, 0, 
          108, 27, 0, 6, 0, 12, 10, 1, 0, 27, 0, 0, 0, 0, 0, 9, 3, 0, 9, 0, 12, 
          9, 3, 0, 1, 0, 2, 0, 0, 0, 3, 1, 1, 0, 0, 0, 3, 0, 0, 0, 1, 0, 0, 1, 1, 
          0, 0, 0, 0, 1 ] ), Character( CharacterTable( "O7(3)" ), 
        [ 12636, 1296, 216, 84, 0, 81, 0, 108, 0, 27, 0, 6, 0, 12, 10, 1, 0, 27, 
          0, 0, 0, 0, 0, 9, 3, 0, 9, 12, 0, 3, 9, 0, 1, 0, 2, 0, 0, 3, 0, 1, 1, 
          0, 0, 0, 3, 0, 0, 0, 1, 0, 0, 1, 1, 0, 0, 0, 0, 1 ] ) ]
    gap> Append( someprim, pi );
    gap> t2:= CharacterTable( "O7(3).2" );;
    gap> s2:= CharacterTable( "Dihedral", 8 ) * CharacterTable( "U4(2).2" );
    CharacterTable( "Dihedral(8)xU4(2).2" )
    gap> pi:= PossiblePermutationCharacters( s2, t2 );;  Length( pi );
    1
    gap> pi:= RestrictedClassFunctions( pi, t );;
    gap> Append( someprim, pi );
    gap> pi:= PossiblePermutationCharacters(
    >             CharacterTable( "2^6:A7" ), t );;  Length( pi );
    1
    gap> Append( someprim, pi );
    gap> List( someprim, x -> x[1] );
    [ 351, 364, 378, 1080, 1080, 1120, 3159, 3159, 12636, 12636, 22113, 28431 ]
\end{verbatim}

Note that in the three cases where two possible permutation characters
were found, there are in fact two classes of subgroups that induce
different permutation characters.
For the subgroups of the types $G_2(3)$ and $S_6(2)$,
this is stated in~\cite[p.~109]{CCN85},
and for the subgroups of the type $S_9$,
this follows from the fact that each $S_9$ type subgroup in $S$ contains
elements in exactly one of the classes {\tt 3D} or {\tt 3E},
and these two classes are fused by the outer automorphism of $S$.

\begin{verbatim}
    gap> cl:= PositionsProperty( AtlasClassNames( t ),
    >                            x -> x in [ "3D", "3E" ] );
    [ 8, 9 ]
    gap> List( Filtered( someprim, x -> x[1] = 12636 ), pi -> pi{ cl } );
    [ [ 0, 108 ], [ 108, 0 ] ]
    gap> GetFusionMap( t, t2 ){ cl };
    [ 8, 8 ]
\end{verbatim}

Now we compute the lower bounds for $\total( S, s^{\prime} )$ that are given
by the sublist \verb|someprim| of the primitive permutation characters.

\begin{verbatim}
    gap> spos:= Position( OrdersClassRepresentatives( t ), 14 );
    52
    gap> Maximum( ApproxP( someprim, spos ) );
    199/351
\end{verbatim}

This shows that $\total( S, s ) = 199/351$ holds.
For statement~(a),
we have to show that choosing $s^{\prime}$ from another class than {\tt 14A}
yields a larger value for $\total( S, s^{\prime} )$.

\begin{verbatim}
    gap> approx:= List( [ 1 .. NrConjugacyClasses( t ) ],
    >       i -> Maximum( ApproxP( someprim, i ) ) );;
    gap> PositionsProperty( approx, x -> x <= 199/351 );
    [ 52 ]
\end{verbatim}

Statement~(b) can be read off from the permutation characters.

\begin{verbatim}
    gap> pos:= PositionsProperty( someprim, x -> x[ spos ] <> 0 );
    [ 1, 9, 10 ]
    gap> List( someprim{ pos }, x -> x{ [ 1, spos ] } );
    [ [ 351, 1 ], [ 12636, 1 ], [ 12636, 1 ] ]
\end{verbatim}

For statement~(c), we first compute $\prop(g, s)$ for $g$ in the class
{\tt 2A}, via explicit computations with the group.
For dealing with this case, we first construct a faithful permutation
representation of $O_7(3)$ from the natural matrix representation of
$\SO(7,3)$.


\begin{verbatim}
    gap> so73:= SpecialOrthogonalGroup( 7, 3 );;
    gap> o73:= DerivedSubgroup( so73 );;
    gap> orbs:= Orbits( o73, Elements( GF(3)^7 ) );;
    gap> Set( List( orbs, Length ) );
    [ 1, 702, 728, 756 ]
    gap> g:= Action( o73, First( orbs, x -> Length( x ) = 702 ) );;
    gap> Size( g ) = Size( t );
    true
\end{verbatim}

A {\tt 2A} element $g$ can be found as the $7$-th power of any element
of order $14$ in $S$.

\begin{verbatim}
    gap> ResetGlobalRandomNumberGenerators();
    gap> repeat s:= Random( g );
    >    until Order( s ) = 14;
    gap> 2A:= s^7;;
    gap> bad:= RatioOfNongenerationTransPermGroup( g, 2A, s );
    155/351
    gap> bad > 1/3;
    true
    gap> approx:= ApproxP( someprim, spos );;
    gap> PositionsProperty( approx, x -> x >= 1/3 );
    [ 2 ]
\end{verbatim}

This shows that $\prop(g,s) = 155/351 > 1/3$.
Since $\total( g, s ) < 1/3$ for all nonidentity $g$ not in the class
{\tt 2A}, we have $\prop( S, s ) = 155/351$.
For statement~(c), it remains to show that $\prop( S, s^{\prime} )$ is larger
than $155/351$ whenever $s^{\prime}$ is not of order $14$.
First we compute $\prop( g, s^{\prime} )$, for $g$ in the class {\tt 2A}.

\begin{verbatim}
    gap> consider:= RepresentativesMaximallyCyclicSubgroups( t );
    [ 18, 19, 25, 26, 27, 30, 31, 32, 34, 35, 38, 39, 41, 42, 43, 44, 45, 46, 47, 
      48, 49, 50, 52, 53, 54, 56, 57, 58 ]
    gap> Length( consider );
    28
    gap> consider:= ClassesPerhapsCorrespondingToTableColumns( g, t, consider );;
    gap> Length( consider );
    31
    gap> consider:= List( consider, Representative );;
    gap> SortParallel( List( consider, Order ), consider );
    gap> app2A:= List( consider, c ->
    >       RatioOfNongenerationTransPermGroup( g, 2A, c ) );
    [ 1, 1, 1, 1, 1, 1, 1, 1, 1, 85/117, 1, 10/13, 10/13, 1, 1, 23/39, 23/39, 1,
      1, 1, 1, 1, 1/3, 1/3, 155/351, 67/117, 1, 1, 1, 1, 191/351 ]
    gap> test:= PositionsProperty( app2A, x -> x <= 155/351 );
    [ 23, 24, 25 ]
    gap> List( test, i -> Order( consider[i] ) );
    [ 13, 13, 14 ]
\end{verbatim}

We see that only for $s^{\prime}$ in one of the two (algebraically conjugate)
classes of element order $13$,
$\prop( S, s^{\prime} )$ has a chance to be smaller than $155/351$.
This possibility is now excluded by counting elements in the class {\tt 3A}
that do not generate $S$ together with $s^{\prime}$ of order $13$.

\begin{verbatim}
    gap> C3A:= First( ConjugacyClasses( g ),
    >               c -> Order( Representative( c ) ) = 3 and Size( c ) = 7280 );;
    gap> repeat ss:= Random( g );
    >    until Order( ss ) = 13;
    gap> bad:= RatioOfNongenerationTransPermGroup( g, Representative( C3A ), ss );
    17/35
    gap> bad > 155/351;
    true
\end{verbatim}

Now we show statement~(d):
For each triple $(x_1, x_2, x_3)$ of nonidentity elements in $S$,
there is an element $s$ in the class {\tt 14A} such that
$\langle x_i, s \rangle = S$ holds for $1 \leq i \leq 3$.
We can read off from the character-theoretic data that
only those triples have to be checked for which at least two
elements are contained in the class {\tt 2A},
and the third element lies in one of the classes
{\tt 2A}, {\tt 2B}, {\tt 3B}.

\begin{verbatim}
    gap> approx:= ApproxP( someprim, spos );;
    gap> max:= Maximum( approx{ [ 3 .. Length( approx ) ] } );
    59/351
    gap> 155 + 2*59 < 351;
    true
    gap> third:= PositionsProperty( approx, x -> 2 * 155/351 + x >= 1 );
    [ 2, 3, 6 ]
    gap> ClassNames( t ){ third };
    [ "2a", "2b", "3b" ]
\end{verbatim}

We can find elements in the classes {\tt 2B} and {\tt 3B} as powers
of arbitrary elements of the orders $20$ and $15$, respectively.

\begin{verbatim}
    gap> ord20:= PositionsProperty( OrdersClassRepresentatives( t ),
    >                               x -> x = 20 );
    [ 58 ]
    gap> PowerMap( t, 10 ){ ord20 };
    [ 3 ]
    gap> repeat x:= Random( g );
    >    until Order( x ) = 20;
    gap> 2B:= x^10;;
    gap> C2B:= ConjugacyClass( g, 2B );;
    gap> ord15:= PositionsProperty( OrdersClassRepresentatives( t ),
    >                               x -> x = 15 );
    [ 53 ]
    gap> PowerMap( t, 10 ){ ord15 };
    [ 6 ]
    gap> repeat x:= Random( g );
    >    until Order( x ) = 15;
    gap> 3B:= x^5;;
    gap> C3B:= ConjugacyClass( g, 3B );;
\end{verbatim}

The existence of $s$ can be shown with the random approach
described in Section~\ref{groups}.

\begin{verbatim}
    gap> repeat s:= Random( g );
    >    until Order( s ) = 14;
    gap> RandomCheckUniformSpread( g, [ 2A, 2A, 2A ], s, 50 );
    true
    gap> RandomCheckUniformSpread( g, [ 2B, 2A, 2A ], s, 50 );
    true
    gap> RandomCheckUniformSpread( g, [ 3B, 2A, 2A ], s, 50 );
    true
\end{verbatim}

Finally, we show statement~(e).
Let $G = \Aut(S) = S.2$.
By~\cite[p.~109]{CCN85}, $\M^{\prime}(G,s)$ consists of the extension
of the $2.U_4(3).2_1$ type subgroup.
We compute the extension of the permutation character.

\begin{verbatim}
    gap> prim:= someprim{ [ 1 ] };
    [ Character( CharacterTable( "O7(3)" ), [ 351, 127, 47, 15, 27, 45, 36, 0, 0, 
          9, 0, 15, 3, 3, 7, 6, 19, 19, 10, 11, 12, 8, 3, 5, 3, 6, 1, 0, 0, 3, 3, 
          0, 1, 1, 1, 6, 3, 0, 0, 2, 2, 0, 3, 0, 3, 3, 0, 0, 1, 0, 0, 1, 0, 4, 4, 
          1, 2, 0 ] ) ]
    gap> spos:= Position( AtlasClassNames( t ), "14A" );;
    gap> t2:= CharacterTable( "O7(3).2" );;
    gap> map:= InverseMap( GetFusionMap( t, t2 ) );;
    gap> torso:= List( prim, pi -> CompositionMaps( pi, map ) );;
    gap> ext:= List( torso, x -> PermChars( t2, rec( torso:= x ) ) );
    [ [ Character( CharacterTable( "O7(3).2" ), [ 351, 127, 47, 15, 27, 45, 36, 
              0, 9, 0, 15, 3, 3, 7, 6, 19, 19, 10, 11, 12, 8, 3, 5, 3, 6, 1, 0, 
              3, 0, 1, 1, 1, 6, 3, 0, 2, 2, 0, 3, 0, 3, 3, 0, 1, 0, 0, 1, 0, 4, 
              1, 2, 0, 117, 37, 21, 45, 1, 13, 5, 1, 9, 9, 18, 15, 1, 7, 9, 6, 4, 
              0, 3, 0, 3, 3, 6, 2, 2, 9, 6, 1, 3, 1, 4, 1, 2, 1, 1, 0, 3, 1, 0, 
              0, 0, 0, 1, 1, 0, 0 ] ) ] ]
    gap> approx:= ApproxP( Concatenation( ext ),
    >        Position( AtlasClassNames( t2 ), "14A" ) );;
    gap> Maximum( approx{ Difference(
    >      PositionsProperty( OrdersClassRepresentatives( t2 ), IsPrimeInt ),
    >      ClassPositionsOfDerivedSubgroup( t2 ) ) } );
    1/3
\end{verbatim}

\subsection{$O_8^+(2)$}\label{O8p2}

We show that $S = O_8^+(2) = \Omega^+(8,2)$ satisfies the following.
\begin{enumerate}
\item[(a)]
    $\total(S) = 334/315$,
    and this value is attained exactly for $\total(S,s)$
    with $s$ of order $15$.
\item[(b)]
    For $s \in S$ of order $15$,
    $\M(S,s)$ consists of one group of the type $S_6(2)$,
    two conjugate groups of the type $2^6:A_8$,
    two conjugate groups of the type $A_9$,
    and one group of each of the types
    $(3 \times U_4(2)):2 = (3 \times \Omega^-(6,2)):2$ and
    $(A_5 \times A_5):2^2 = (\Omega^-(4,2) \times \Omega^-(4,2)):2^2$.
\item[(c)]
    $\prop(S) = 29/42$,
    and this value is attained exactly for $\prop(S,s)$
    with $s$ of order $15$.
\item[(d)]
    Let $x, y \in S$ such that $x, y, x y$ lie in the unique involution class
    of length $1\,575$ of $S$.
    (This is the class \verb|2A|.)
    Then each element in $S$ together with one of $x$, $y$, $x y$
    generates a proper subgroup of $S$.
\item[(e)]
    Both the spread and the uniform spread of $S$ is exactly two,
    with $s$ of order $15$.
\item[(f)]
    For each choice of $s \in S$, there is an extension $S.2$ such that
    for any element $g$ in the (outer) class {\tt 2F},
    $\langle s, g \rangle$ does not contain $S$.
\item[(g)]
    For an element $s$ of order $15$ in $S$,
    either $S$ is the only maximal subgroup of $S.2$ that contains $s$,
    or the maximal subgroups of $S.2$ that contain $s$
    are $S$ and the extensions of the subgroups listed in statement~(b);
    these groups have the structures $S_6(2) \times 2$, $2^6:S_8$ (twice),
    $S_9$ (twice), $S_3 \times U_4(2).2$, and
\tthdump{$S_5 \wr 2$.}
\item[(h)]
    For $s \in S$ of order $15$ and arbitrary $g \in S.3 \setminus S$,
    we have $\langle s, g \rangle = S.3$.
\item[(i)]
    If $x$, $y$ are nonidentity elements in $\Aut(S)$ then there is an
    element $s$ of order $15$ in $S$ such that
    $S \subseteq \langle x, s \rangle \cap \langle y, s \rangle$.
\end{enumerate}

Statement~(a) follows from inspection of the primitive permutation
characters, cf.~Section~\ref{easyloop}.

\begin{verbatim}
    gap> t:= CharacterTable( "O8+(2)" );;
    gap> ProbGenInfoSimple( t );
    [ "O8+(2)", 334/315, 0, [ "15A", "15B", "15C" ], [ 7, 7, 7 ] ]
\end{verbatim}

Statement~(b) can be read off from the permutation characters,
and the fact that the only classes of maximal subgroups that contain
elements of order $15$ consist of groups of the structures
as claimed, see~\cite[p.~85]{CCN85}.

\begin{verbatim}
    gap> prim:= PrimitivePermutationCharacters( t );;
    gap> spos:= Position( OrdersClassRepresentatives( t ), 15 );;
    gap> List( Filtered( prim, x -> x[ spos ] <> 0 ), l -> l{ [ 1, spos ] } );
    [ [ 120, 1 ], [ 135, 2 ], [ 960, 2 ], [ 1120, 1 ], [ 12096, 1 ] ]
\end{verbatim}

For the remaining statements,
we take a primitive permutation representation on $120$ points,
and assume that the permutation character is {\tt 1a+35a+84a}.
(See~\cite[p.~85]{CCN85}, note that the three classes of maximal subgroups
of index $120$ in $S$ are conjugate under triality.)

\begin{verbatim}
    gap> matgroup:= DerivedSubgroup( GeneralOrthogonalGroup( 1, 8, 2 ) );;
    gap> points:= NormedRowVectors( GF(2)^8 );;
    gap> orbs:= Orbits( matgroup, points );;
    gap> List( orbs, Length );
    [ 135, 120 ]
    gap> g:= Action( matgroup, orbs[2] );;
    gap> Size( g );
    174182400
    gap> pi:= Sum( Irr( t ){ [ 1, 3, 7 ] } );
    Character( CharacterTable( "O8+(2)" ), [ 120, 24, 32, 0, 0, 8, 36, 0, 0, 3, 
      6, 12, 4, 8, 0, 0, 0, 10, 0, 0, 12, 0, 0, 8, 0, 0, 3, 6, 0, 0, 2, 0, 0, 2, 
      1, 2, 2, 3, 0, 0, 2, 0, 0, 0, 0, 0, 3, 2, 0, 0, 1, 0, 0 ] )
\end{verbatim}

In order to show statement~(c),
we first observe that for $s$ in the class {\tt 15A}
and $g$ \emph{not} in one of the classes {\tt 2A}, {\tt 2B}, {\tt 3A},
$\total(g,s) < 1/3$ holds,
and for the exceptional three classes,
we have $\total(g,s) > 1/2$.

\begin{verbatim}
    gap> approx:= ApproxP( prim, spos );;
    gap> testpos:= PositionsProperty( approx, x -> x >= 1/3 );
    [ 2, 3, 7 ]
    gap> AtlasClassNames( t ){ testpos };
    [ "2A", "2B", "3A" ]
    gap> approx{ testpos };
    [ 254/315, 334/315, 1093/1120 ]
    gap> ForAll( approx{ testpos }, x -> x > 1/2 );
    true
\end{verbatim}

Now we compute the values $\prop(g,s)$, for $s$ in the class {\tt 15A}
and $g$ in one of the classes {\tt 2A}, {\tt 2B}, {\tt 3A}.

By our choice of the character of the permutation representation we use,
the class {\tt 15A} is determined as the unique class of element order $15$
with one fixed point.
(Note that the three classes of element order $15$ in $S$ are conjugate
under triality.)
A {\tt 2A} element can be found as the fourth power of any element of order
$8$ in $S$,
a {\tt 3A} element can be found as the fifth power of a {\tt 15A} element,
and a {\tt 2B} element can be found as the sixth power of an element
of order $12$, with $32$ fixed points.

\begin{verbatim}
    gap> ResetGlobalRandomNumberGenerators();
    gap> repeat s:= Random( g );
    >    until Order( s ) = 15 and NrMovedPoints( g ) = 1 + NrMovedPoints( s );
    gap> 3A:= s^5;;
    gap> repeat x:= Random( g ); until Order( x ) = 8;
    gap> 2A:= x^4;;
    gap> repeat x:= Random( g ); until Order( x ) = 12 and
    >      NrMovedPoints( g ) = 32 + NrMovedPoints( x^6 );
    gap> 2B:= x^6;;
    gap> prop15A:= List( [ 2A, 2B, 3A ],
    >                    x -> RatioOfNongenerationTransPermGroup( g, x, s ) );
    [ 23/35, 29/42, 149/224 ]
    gap> Maximum( prop15A );
    29/42
\end{verbatim}

This means that for $s$ in the class {\tt 15A},
we have $\prop( S, s ) = 29/42$,
and the same holds for all $s$ of order $15$
since the three classes of element order $15$ are conjugate under triality.
Now we show that for $s$ of order different from $15$,
the value $\prop(g,s)$ is larger than $29/42$,
for $g$ in one of the classes {\tt 2A}, {\tt 2B}, {\tt 3A},
or their images under triality.
This implies statement~(c).

\begin{verbatim}
    gap> test:= List( [ 2A, 2B, 3A ], x -> ConjugacyClass( g, x ) );;
    gap> ccl:= ConjugacyClasses( g );;
    gap> consider:= Filtered( ccl, c -> Size( c ) in List( test, Size ) );;
    gap> Length( consider );
    7
    gap> filt:= Filtered( ccl, c -> ForAll( consider, cc ->
    >       RatioOfNongenerationTransPermGroup( g, Representative( cc ),
    >           Representative( c ) ) <= 29/42 ) );;
    gap> Length( filt );
    3
    gap> List( filt, c -> Order( Representative( c ) ) );
    [ 15, 15, 15 ]
\end{verbatim}


Now we show statement~(d).
First we observe that all those Klein four groups in $S$ whose involutions
lie in the class \verb|2A| are conjugate in $S$.
Note that this is the unique class of length $1\,575$ in $S$,
and also the unique class whose elements have $24$ fixed points
in the degree $120$ permutation representation.

For that, we use the character table of $S$ to read off that $S$ contains
exactly $14\,175$ such subgroups,
and we use the group to compute one such subgroup and its normalizer
of index $14\,175$.

\begin{verbatim}
    gap> SizesConjugacyClasses( t );
    [ 1, 1575, 3780, 3780, 3780, 56700, 2240, 2240, 2240, 89600, 268800, 37800,
      340200, 907200, 907200, 907200, 2721600, 580608, 580608, 580608, 100800,
      100800, 100800, 604800, 604800, 604800, 806400, 806400, 806400, 806400,
      2419200, 2419200, 2419200, 7257600, 24883200, 5443200, 5443200, 6451200,
      6451200, 6451200, 8709120, 8709120, 8709120, 1209600, 1209600, 1209600,
      4838400, 7257600, 7257600, 7257600, 11612160, 11612160, 11612160 ]
    gap> NrPolyhedralSubgroups( t, 2, 2, 2 );
    rec( number := 14175, type := "V4" )
    gap> repeat x:= Random( g );
    >    until     Order( x ) mod 2 = 0
    >          and NrMovedPoints( x^( Order(x)/2 ) ) = 120 - 24;
    gap> x:= x^( Order(x)/2 );;
    gap> repeat y:= x^Random( g );
    >    until NrMovedPoints( x*y ) = 120 - 24;
    gap> v4:= SubgroupNC( g, [ x, y ] );;
    gap> n:= Normalizer( g, v4 );;
    gap> Index( g, n );
    14175
\end{verbatim}

We verify that the triple has the required property.

\begin{verbatim}
    gap> maxorder:= RepresentativesMaximallyCyclicSubgroups( t );;
    gap> maxorderreps:= List( ClassesPerhapsCorrespondingToTableColumns( g, t,
    >        maxorder ), Representative );;
    gap> Length( maxorderreps );
    28
    gap> CommonGeneratorWithGivenElements( g, maxorderreps, [ x, y, x*y ] );
    fail
\end{verbatim}

For the simple group $S$, it remains to show statement~(e).
We want to show that for any choice of two nonidentity elements $x$, $y$
in $S$, there is an element $s$ in the class {\tt 15A} such that
$\langle s, x \rangle = \langle s, y \rangle = S$ holds.
Only $x$, $y$ in the classes given by the list \verb|testpos| must be considered,
by the estimates $\total(g,s)$.

We replace the values $\total(g,s)$ by the exact values $\prop(g,s)$,
for $g$ in one of these three classes.
Each of the three classes is determined by its element order and its number
of fixed points.

\begin{verbatim}
    gap> reps:= List( ccl, Representative );;
    gap> bading:= List( testpos, i -> Filtered( reps,
    >        r -> Order( r ) = OrdersClassRepresentatives( t )[i] and
    >             NrMovedPoints( r ) = 120 - pi[i] ) );;
    gap> List( bading, Length );
    [ 1, 1, 1 ]
    gap> bading:= List( bading, x -> x[1] );;
\end{verbatim}

For each pair $(C_1, C_2)$ of classes represented by this list,
we have to show that for any choice of elements $x \in C_1$, $y \in C_2$
there is $s$ in the class {\tt 15A} such that
$\langle s, x \rangle = \langle s, y \rangle = S$ holds.
This is done with the random approach that is described in
Section~\ref{groups}.

\begin{verbatim}
    gap> for pair in UnorderedTuples( bading, 2 ) do
    >      test:= RandomCheckUniformSpread( g, pair, s, 80 );
    >      if test <> true then
    >        Error( test );
    >      fi;
    >    od;
\end{verbatim}

We get no error message, so statement~(e) holds.

Now we turn to the automorphic extensions of $S$.
First we compute a permutation representation of $\SO^+(8,2) \cong S.2$
and an element $g$ in the class {\tt 2F},
which is the unique conjugacy class of size $120$ in $S.2$.

\begin{verbatim}
    gap> matgrp:= SO(1,8,2);;
    gap> g2:= Image( IsomorphismPermGroup( matgrp ) );;
    gap> IsTransitive( g2, MovedPoints( g2 ) );
    true
    gap> repeat x:= Random( g2 ); until Order( x ) = 14;
    gap> 2F:= x^7;;
    gap> Size( ConjugacyClass( g2, 2F ) );
    120
\end{verbatim}

Only for $s$ in six conjugacy classes of $S$,
there is a nonzero probability to have $S.2 = \langle g, s \rangle$.

\begin{verbatim}
    gap> der:= DerivedSubgroup( g2 );;
    gap> cclreps:= List( ConjugacyClasses( der ), Representative );;
    gap> nongen:= List( cclreps,
    >               x -> RatioOfNongenerationTransPermGroup( g2, 2F, x ) );;
    gap> goodpos:= PositionsProperty( nongen, x -> x < 1 );;
    gap> invariants:= List( goodpos, i -> [ Order( cclreps[i] ),
    >      Size( Centralizer( g2, cclreps[i] ) ), nongen[i] ] );;
    gap> SortedList( invariants );
    [ [ 10, 20, 1/3 ], [ 10, 20, 1/3 ], [ 12, 24, 2/5 ], [ 12, 24, 2/5 ],
      [ 15, 15, 0 ], [ 15, 15, 0 ] ]
\end{verbatim}

$S$ contains three classes of element order $10$, which are conjugate
in $S.3$.
For a fixed extension of the type $S.2$,
the element $s$ can be chosen only in two of these three classes,
which means that there is another group of the type $S.2$
(more precisely, another subgroup of index three in $S.S_3$)
in which this choice of $s$ is not suitable
-- note that the general aim is to find $s \in S$ uniformly for all
automorphic extensions of $S$.
Analogous statements hold for the other possibilities for $s$,
so statement~(f) follows.

Statement~(g) follows from the list of maximal subgroups
in~\cite[p.~85]{CCN85}.

Statement~(h) follows from the fact that $S$ is the only maximal
subgroup of $S.3$ that contains elements of order $15$,
according to the list of maximal subgroups in~\cite[p.~85]{CCN85}.
Alternatively, if we do not want to assume this information,
we can use explicit computations, as follows.
All we have to check is that any element in the classes {\tt 3F} and {\tt 3G}
generates $S.3$ together with a fixed element of order $15$ in $S$.

We compute a permutation representation of $S.3$ as the derived subgroup
of a subgroup of the type $S.S_3$ inside the sporadic simple Fischer group
$Fi_{22}$;
these subgroups lie in the fourth class of maximal subgroups of $Fi_{22}$,
see~\cite[p.~163]{CCN85}.
An element in the class {\tt 3F} of $S.3$ can be found as a power of
an order $21$ element,
and an element in the class {\tt 3G} can be found as the fourth power of
a {\tt 12P} element.

\begin{verbatim}
    gap> aut:= Group( AtlasGenerators( "Fi22", 1, 4 ).generators );;
    gap> Size( aut ) = 6 * Size( t );
    true
    gap> g3:= DerivedSubgroup( aut );;
    gap> orbs:= Orbits( g3, MovedPoints( g3 ) );;
    gap> List( orbs, Length );
    [ 3150, 360 ]
    gap> g3:= Action( g3, orbs[2] );;
    gap> repeat s:= Random( g3 ); until Order( s ) = 15;
    gap> repeat x:= Random( g3 ); until Order( x ) = 21;
    gap> 3F:= x^7;;
    gap> RatioOfNongenerationTransPermGroup( g3, 3F, s );
    0
    gap> repeat x:= Random( g3 );
    >    until Order( x ) = 12 and Size( Centralizer( g3, x^4 ) ) = 648;
    gap> 3G:= x^4;;
    gap> RatioOfNongenerationTransPermGroup( g3, 3G, s );
    0
\end{verbatim}

Finally, consider statement~(i).
It implies that~\cite[Corollary~1.5]{BGK} holds
for $\Omega^+(8,2)$, with $s$ of order $15$.
Note that by part~(f),
$s$ \emph{cannot be chosen in a prescribed conjugacy class} of $S$
that is independent of the elements $x$, $y$.

If $x$ and $y$ lie in $S$ then statement~(i) follows from part~(e),
and by part~(g), the case that $x$ or $y$ lie in $S.3 \setminus S$
is also not a problem.
We now show that also $x$ or $y$ in $S.2 \setminus S$ is not a problem.
Here we have to deal with the cases that $x$ and $y$ lie in the same
subgroup of index $3$ in $\Aut(S)$ or in different such subgroups.
Actually we show that for each index $3$ subgroup $H = S.2 < \Aut(S)$,
we can choose $s$ from two of the three classes of element order $15$ in $S$
such that $S$ is the only maximal subgroup of $H$ that contains $s$,
and thus $\langle x, s \rangle$ contains $H$,
for any choice of $x \in H \setminus S$.

For that, we note that no novelty in $S.2$ contains elements of order $15$,
so all maximal subgroups of $S.2$ that contain such elements --besides $S$--
have one of the indices $120, 135, 960, 1120$, or $12096$,
and point stabilizers of the types $S_6(2) \times 2$, $2^6:S_8$, $S_9$,
$S_3 \times U_4(2):2$, or
\tthdump{$S_5 \wr 2$.}
We compute the corresponding permutation characters.

\begin{verbatim}
    gap> t2:= CharacterTable( "O8+(2).2" );;
    gap> s:= CharacterTable( "S6(2)" ) * CharacterTable( "Cyclic", 2 );;
    gap> pi:= PossiblePermutationCharacters( s, t2 );;
    gap> prim:= pi;;
    gap> pi:= PermChars( t2, rec( torso:= [ 135 ] ) );;
    gap> Append( prim, pi );
    gap> pi:= PossiblePermutationCharacters( CharacterTable( "A9.2" ), t2 );;
    gap> Append( prim, pi );
    gap> s:= CharacterTable( "Dihedral(6)" ) * CharacterTable( "U4(2).2" );;
    gap> pi:= PossiblePermutationCharacters( s, t2 );;
    gap> Append( prim, pi );
    gap> s:= CharacterTableWreathSymmetric( CharacterTable( "S5" ), 2 );;
    gap> pi:= PossiblePermutationCharacters( s, t2 );;
    gap> Append( prim, pi );
    gap> Length( prim );
    5
    gap> ord15:= PositionsProperty( OrdersClassRepresentatives( t2 ),
    >                               x -> x = 15 );
    [ 39, 40 ]
    gap> List( prim, pi -> pi{ ord15 } );
    [ [ 1, 0 ], [ 2, 0 ], [ 2, 0 ], [ 1, 0 ], [ 1, 0 ] ]
    gap> List( ord15, i -> Maximum( ApproxP( prim, i ) ) );
    [ 307/120, 0 ]
\end{verbatim}

Here it is appropriate to clean the workspace again.

\begin{verbatim}
    gap> CleanWorkspace();
\end{verbatim}

\subsection{$O_8^+(3)$}\label{O8p3}

We show that $S = O_8^+(3)$ satisfies the following.
\begin{enumerate}
\item[(a)]
    $\total(S) = 863/1820$,
    and this value is attained exactly for $\total(S,s)$
    with $s$ of order $20$.
\item[(b)]
    For $s \in S$ of order $20$,
    $\M(S,s)$ consists of two nonconjugate groups of the type
    $O_7(3) = \Omega(7,3)$,
    two conjugate subgroups of the type $3^6:L_4(3)$,
    two nonconjugate subgroups of the type $(A_4 \times U_4(2)):2$,
    and one subgroup of each of the types $2.U_4(3).(2^2)_{122}$
    and $(A_6 \times A_6):2^2$.
\item[(c)]
    $\prop(S) = 194/455$,
    and this value is attained exactly for $\prop(S,s)$
    with $s$ of order $20$.
\item[(d)]
    The uniform spread of $S$ is at least three,
    with $s$ of order $20$.
\item[(e)]
    The preimage of $s$ in the matrix group $2.S = \Omega^+(8,3)$
    can be chosen of order $40$,
    and then the maximal subgroups of $2.S$ containing $s$
    have the structures $2.O_7(3)$, $3^6:2.L_4(3)$,
    $4.U_4(3).2^2 = \SU(4,3).2^2$,
    $2.(A_4 \times U_4(2)).2 = 2.(\PSp(2,3) \otimes \PSp(4,3)).2$, and
    $2.(A_6 \times A_6):2^2 = 2.(\Omega^-(4,3) \times \Omega^-(4,3)):2^2$,
    respectively.
\item[(f)]
    For $s \in S$ of order $20$, we have
    $\prop^{\prime}(S.2_1, s) \in \{ 83/567, 574/1215 \}$,
    $\prop^{\prime}(S.2_2, s) \in \{ 0, 1 \}$
    (depending on the choice of $s$), and
    $\total^{\prime}(S.3, s) = 0$.

    Furthermore, for any choice of $s^{\prime} \in S$,
    we have $\total^{\prime}(S.2_2, s^{\prime}) = 1$
    for some group $S.2_2$.
    However, if it is allowed to choose $s$ from an $\Aut(S)$-class
    of elements of order $20$ (and not from a fixed $S$-class)
    then we can achieve $\total(g,s) = 0$ for any given
    $g \in S.2_2 \setminus S$.
\item[(g)]
    The maximal subgroups of $S.2_1$ that contain an element of order $20$
    are either
    $S$ and the extensions of the subgroups listed in statement~(b)
    or they are
    $S$ and $L_4(3).2^2$, $3^6:L_4(3).2$ (twice),
    $2.U_4(3).(2^2)_{122}.2$, and
    $(A_6 \times A_6):2^2.2$.

    In the former case, the groups have the structures
    $O_7(3):2$ (twice),
    $3^6:(L_4(3) \times 2)$ (twice),
    $S_4 \times U_4(2).2$ (twice),
    $2.U_4(3).(2^2)_{122}.2$, and
    $(A_6 \times A_6):2^2 \times 2$.
\end{enumerate}

Statement~(a) follows from inspection of the primitive permutation
characters.

\begin{verbatim}
    gap> t:= CharacterTable( "O8+(3)" );;
    gap> ProbGenInfoSimple( t );
    [ "O8+(3)", 863/1820, 2, [ "20A", "20B", "20C" ], [ 8, 8, 8 ] ]
\end{verbatim}

Also statement~(b) follows from the information provided by the
character table of $S$ (cf.~\cite[p.~140]{CCN85}).

\begin{verbatim}
    gap> prim:= PrimitivePermutationCharacters( t );;
    gap> ord:= OrdersClassRepresentatives( t );;
    gap> spos:= Position( ord, 20 );;
    gap> filt:= PositionsProperty( prim, x -> x[ spos ] <> 0 );
    [ 1, 2, 7, 15, 18, 19, 24 ]
    gap> Maxes( t ){ filt };
    [ "O7(3)", "O8+(3)M2", "3^6:L4(3)", "2.U4(3).(2^2)_{122}", "(A4xU4(2)):2", 
      "O8+(3)M19", "(A6xA6):2^2" ]
    gap> prim{ filt }{ [ 1, spos ] };
    [ [ 1080, 1 ], [ 1080, 1 ], [ 1120, 2 ], [ 189540, 1 ], [ 7960680, 1 ], 
      [ 7960680, 1 ], [ 9552816, 1 ] ]
\end{verbatim}

For statement~(c), we first show that $\prop(S,s) = 194/455$ holds.
Since this value is larger than $1/3$,
we have to inspect only those classes $g^S$ for which
$\total(g,s) \geq 1/3$ holds,

\begin{verbatim}
    gap> ord:= OrdersClassRepresentatives( t );;
    gap> ord20:= PositionsProperty( ord, x -> x = 20 );;
    gap> cand:= [];;
    gap> for i in ord20 do
    >      approx:= ApproxP( prim, i );
    >      Add( cand, PositionsProperty( approx, x -> x >= 1/3 ) );
    >    od;
    gap> cand;
    [ [ 2, 6, 7, 10 ], [ 3, 6, 8, 11 ], [ 4, 6, 9, 12 ] ]
    gap> AtlasClassNames( t ){ cand[1] };
    [ "2A", "3A", "3B", "3E" ]
\end{verbatim}

The three possibilities form one orbit under the outer automorphism group
of $S$.

\begin{verbatim}
    gap> t3:= CharacterTable( "O8+(3).3" );;
    gap> tfust3:= GetFusionMap( t, t3 );;
    gap> List( cand, x -> tfust3{ x } );
    [ [ 2, 4, 5, 6 ], [ 2, 4, 5, 6 ], [ 2, 4, 5, 6 ] ]
\end{verbatim}

By symmetry, we may consider only the first possibility,
and assume that $s$ is in the class {\tt 20A}.

We work with a permutation representation of degree $1\,080$,
and assume that the permutation character is {\tt 1a+260a+819a}.
(Note that all permutation characters of $S$ of degree $1\,080$
are conjugate under $\Aut(S)$.)

\begin{verbatim}
    gap> g:= Action( SO(1,8,3), NormedRowVectors( GF(3)^8 ), OnLines );;
    gap> Size( g );
    9904359628800
    gap> g:= DerivedSubgroup( g );;  Size( g );
    4952179814400
    gap> orbs:= Orbits( g, MovedPoints( g ) );;
    gap> List( orbs, Length );
    [ 1080, 1080, 1120 ]
    gap> g:= Action( g, orbs[1] );;
    gap> PositionProperty( Irr( t ), chi -> chi[1] = 819 );
    9
    gap> permchar:= Sum( Irr( t ){ [ 1, 2, 9 ] } );
    Character( CharacterTable( "O8+(3)" ), [ 1080, 128, 0, 0, 24, 108, 135, 0, 0, 
      108, 0, 0, 27, 27, 0, 0, 18, 9, 12, 16, 0, 0, 4, 15, 0, 0, 20, 0, 0, 12, 
      11, 0, 0, 20, 0, 0, 15, 0, 0, 12, 0, 0, 2, 0, 0, 3, 3, 0, 0, 6, 6, 0, 0, 3, 
      2, 2, 2, 18, 0, 0, 9, 0, 0, 0, 0, 0, 0, 3, 3, 0, 0, 3, 0, 0, 12, 0, 0, 3, 
      0, 0, 0, 0, 0, 4, 3, 3, 0, 0, 1, 0, 0, 4, 0, 0, 1, 1, 2, 0, 0, 0, 0, 0, 3, 
      0, 0, 2, 0, 0, 5, 0, 0, 1, 0, 0 ] )
\end{verbatim}

Now we show that for $s$ in the class {\tt 20A} (which fixes one point),
the proportion of nongenerating elements $g$ in one of the classes
{\tt 2A}, {\tt 3A}, {\tt 3B}, {\tt 3E} has the maximum $194/455$,
which is attained exactly for {\tt 3A}.
(We find a {\tt 2A} element as a power of $s$,
a {\tt 3A} element as a power of any element of order $18$,
a {\tt 3B} and a {\tt 3E} element as elements with $135$ and $108$
fixed points, respectively, which occur as powers of suitable
elements of order $15$.)

\begin{verbatim}
    gap> permchar{ ord20 };
    [ 1, 0, 0 ]
    gap> AtlasClassNames( t )[ PowerMap( t, 10 )[ ord20[1] ] ];
    "2A"
    gap> ord18:= PositionsProperty( ord, x -> x = 18 );;
    gap> Set( AtlasClassNames( t ){ PowerMap( t, 6 ){ ord18 } } );
    [ "3A" ]
    gap> ord15:= PositionsProperty( ord, x -> x = 15 );;
    gap> PowerMap( t, 5 ){ ord15 };
    [ 7, 8, 9, 10, 11, 12 ]
    gap> AtlasClassNames( t ){ [ 7 .. 12 ] };
    [ "3B", "3C", "3D", "3E", "3F", "3G" ]
    gap> permchar{ [ 7 .. 12 ] };
    [ 135, 0, 0, 108, 0, 0 ]
    gap> mp:= NrMovedPoints( g );;
    gap> ResetGlobalRandomNumberGenerators();
    gap> repeat 20A:= Random( g );
    >    until Order( 20A ) = 20 and mp - NrMovedPoints( 20A ) = 1;
    gap> 2A:= 20A^10;;
    gap> repeat x:= Random( g ); until Order( x ) = 18;
    gap> 3A:= x^6;;
    gap> repeat x:= Random( g );
    >    until Order( x ) = 15 and mp - NrMovedPoints( x^5 ) = 135;
    gap> 3B:= x^5;;
    gap> repeat x:= Random( g );
    >    until Order( x ) = 15 and mp - NrMovedPoints( x^5 ) = 108;
    gap> 3E:= x^5;;
    gap> nongen:= List( [ 2A, 3A, 3B, 3E ],
    >                   c -> RatioOfNongenerationTransPermGroup( g, c, 20A ) );
    [ 3901/9477, 194/455, 451/1092, 451/1092 ]
    gap> Maximum( nongen );
    194/455
\end{verbatim}

Next we compute the values $\prop(g,s)$, for $g$ is in the class {\tt 3A}
and certain elements $s$.
It is enough to consider representatives $s$ of maximally cyclic subgroups
in $S$, but here we can do better, as follows.
Since {\tt 3A} is the unique class of length $72\,800$,
it is fixed under $\Aut(S)$,
so it is enough to consider one element $s$ from each $\Aut(S)$-orbit
on the classes of $S$.
We use the class fusion between the character tables of $S$ and $\Aut(S)$
for computing orbit representatives.

\begin{verbatim}
    gap> maxorder:= RepresentativesMaximallyCyclicSubgroups( t );;
    gap> Length( maxorder );
    57
    gap> autt:= CharacterTable( "O8+(3).S4" );;
    gap> fus:= PossibleClassFusions( t, autt );;
    gap> orbreps:= Set( List( fus, map -> Set( ProjectionMap( map ) ) ) );
    [ [ 1, 2, 5, 6, 7, 13, 17, 18, 19, 20, 23, 24, 27, 30, 31, 37, 43, 46, 50, 
          54, 55, 56, 57, 58, 64, 68, 72, 75, 78, 84, 85, 89, 95, 96, 97, 100, 
          106, 112 ] ]
    gap> totest:= Intersection( maxorder, orbreps[1] );
    [ 43, 50, 54, 56, 57, 64, 68, 75, 78, 84, 85, 89, 95, 97, 100, 106, 112 ]
    gap> Length( totest );
    17
    gap> AtlasClassNames( t ){ totest };
    [ "6Q", "6X", "6B1", "8A", "8B", "9G", "9K", "12A", "12D", "12J", "12K", 
      "12O", "13A", "14A", "15A", "18A", "20A" ]
\end{verbatim}

This means that we have to test one element of each of the element orders
$13$, $14$, $15$, and $18$
(note that we know already a bound for elements of order $20$),
plus certain elements of the orders $6$, $8$, $9$, and $12$
which can be identified by their centralizer orders and (for elements of
order $6$ and $8$) perhaps the centralizer orders of some powers.

\begin{verbatim}
    gap> elementstotest:= [];;
    gap> for elord in [ 13, 14, 15, 18 ] do
    >      repeat s:= Random( g );
    >      until Order( s ) = elord;
    >      Add( elementstotest, s );
    >    od;
\end{verbatim}

The next elements to be tested are
in the classes \verb|6B1| (centralizer order $162$),
in one of \verb|9G|--\verb|9J| (centralizer order $729$),
in one of \verb|9K|--\verb|9N| (centralizer order $81$),
in one of \verb|12A|--\verb|12C| (centralizer order $1\,728$),
in one of \verb|12D|--\verb|12I| (centralizer order $432$),
in \verb|12J| (centralizer order $192$),
in one of \verb|12K|--\verb|12N| (centralizer order $108$),
and in one of \verb|12O|--\verb|12T| (centralizer order $72$).

\begin{verbatim}
    gap> ordcent:= [ [ 6, 162 ], [ 9, 729 ], [ 9, 81 ], [ 12, 1728 ],
    >                [ 12, 432 ], [ 12, 192 ], [ 12, 108 ], [ 12, 72 ] ];;
    gap> cents:= SizesCentralizers( t );;
    gap> for pair in ordcent do
    >      Print( pair, ": ", AtlasClassNames( t ){
    >          Filtered( [ 1 .. Length( ord ) ],
    >                    i -> ord[i] = pair[1] and cents[i] = pair[2] ) }, "\n" );
    >      repeat s:= Random( g );
    >      until Order( s ) = pair[1] and Size( Centralizer( g, s ) ) = pair[2];
    >      Add( elementstotest, s );
    >    od;
    [ 6, 162 ]: [ "6B1" ]
    [ 9, 729 ]: [ "9G", "9H", "9I", "9J" ]
    [ 9, 81 ]: [ "9K", "9L", "9M", "9N" ]
    [ 12, 1728 ]: [ "12A", "12B", "12C" ]
    [ 12, 432 ]: [ "12D", "12E", "12F", "12G", "12H", "12I" ]
    [ 12, 192 ]: [ "12J" ]
    [ 12, 108 ]: [ "12K", "12L", "12M", "12N" ]
    [ 12, 72 ]: [ "12O", "12P", "12Q", "12R", "12S", "12T" ]
\end{verbatim}

The next elements to be tested are
in one of the classes \verb|6Q|--\verb|6S| (centralizer order $648$).

\begin{verbatim}
    gap> AtlasClassNames( t ){ Filtered( [ 1 .. Length( ord ) ],
    >        i -> cents[i] = 648 and cents[ PowerMap( t, 2 )[i] ] = 52488
    >                            and cents[ PowerMap( t, 3 )[i] ] = 26127360 ) };
    [ "6Q", "6R", "6S" ]
    gap> repeat s:= Random( g );
    >    until Order( s ) = 6 and Size( Centralizer( g, s ) ) = 648
    >      and Size( Centralizer( g, s^2 ) ) = 52488
    >      and Size( Centralizer( g, s^3 ) ) = 26127360;
    gap> Add( elementstotest, s );
\end{verbatim}

The next elements to be tested are
in the class \verb|6X|--\verb|6A1| (centralizer order $648$).

\begin{verbatim}
    gap> AtlasClassNames( t ){ Filtered( [ 1 .. Length( ord ) ],
    >        i -> cents[i] = 648 and cents[ PowerMap( t, 2 )[i] ] = 52488
    >                            and cents[ PowerMap( t, 3 )[i] ] = 331776 ) };
    [ "6X", "6Y", "6Z", "6A1" ]
    gap> repeat s:= Random( g );
    >    until Order( s ) = 6 and Size( Centralizer( g, s ) ) = 648
    >      and Size( Centralizer( g, s^2 ) ) = 52488
    >      and Size( Centralizer( g, s^3 ) ) = 331776;
    gap> Add( elementstotest, s );
\end{verbatim}

Finally, we add elements from the classes \verb|8A| and \verb|8B|.

\begin{verbatim}
    gap> AtlasClassNames( t ){ Filtered( [ 1 .. Length( ord ) ],
    >        i -> ord[i] = 8 and cents[ PowerMap( t, 2 )[i] ] = 13824 ) };
    [ "8A" ]
    gap> repeat s:= Random( g );
    >    until Order( s ) = 8 and Size( Centralizer( g, s^2 ) ) = 13824;
    gap> Add( elementstotest, s );
    gap> AtlasClassNames( t ){ Filtered( [ 1 .. Length( ord ) ],
    >        i -> ord[i] = 8 and cents[ PowerMap( t, 2 )[i] ] = 1536 ) };
    [ "8B" ]
    gap> repeat s:= Random( g );
    >    until Order( s ) = 8 and Size( Centralizer( g, s^2 ) ) = 1536;
    gap> Add( elementstotest, s );
\end{verbatim}

Now we compute the ratios.
It turns out that from these candidates,
only elements $s$ of the orders $14$ and $15$ satisfy $\prop(g,s) < 194/455$.

\begin{verbatim}
    gap> nongen:= List( elementstotest,
    >                   s -> RatioOfNongenerationTransPermGroup( g, 3A, s ) );;
    gap> smaller:= PositionsProperty( nongen, x -> x < 194/455 );
    [ 2, 3 ]
    gap> nongen{ smaller };
    [ 127/325, 1453/3640 ]
\end{verbatim}

So the only candidates for $s$ that may be better than order $20$ elements
are elements of order $14$ or $15$.
In order to exclude these two possibilities,
we compute $\prop(g,s)$ for $s$ in the class {\tt 14A}
and $g = s^7$ in the class {\tt 2A},
and for $s$ in the class {\tt 15A} and $g$ in the class {\tt 2A},
which yields values that are larger than $194/455$.

\begin{verbatim}
    gap> repeat s:= Random( g );
    >    until Order( s ) = 14 and NrMovedPoints( s ) = 1078;
    gap> 2A:= s^7;;
    gap> nongen:= RatioOfNongenerationTransPermGroup( g, 2A, s );
    1573/3645
    gap> nongen > 194/455;
    true
    gap> repeat s:= Random( g );
    >    until Order( s ) = 15 and NrMovedPoints( s ) = 1080 - 3;
    gap> nongen:= RatioOfNongenerationTransPermGroup( g, 2A, s );
    490/1053
    gap> nongen > 194/455;
    true
\end{verbatim}

For statement~(d), we show that for each triple of elements
in the union of the classes {\tt 2A}, {\tt 3A}, {\tt 3B}, {\tt 3E}
there is an element in the class {\tt 20A} that generates $S$
together with each element of the triple.

\begin{verbatim}
    gap> for tup in UnorderedTuples( [ 2A, 3A, 3B, 3E ], 3 ) do
    >      cl:= ShallowCopy( tup );
    >      test:= RandomCheckUniformSpread( g, cl, 20A, 100 );
    >      if test <> true then
    >        Error( test );
    >      fi;
    >    od;
\end{verbatim}

We get no error message, so statement~(d) is true.

For statement~(e), first we show that $2.S = \Omega^+(8,3)$ contains
elements of order $40$ but $S$ does not.

\begin{verbatim}
    gap> der:= DerivedSubgroup( SO(1,8,3) );;
    gap> repeat x:= PseudoRandom( der ); until Order( x ) = 40;
    gap> 40 in ord;
    false
\end{verbatim}

Thus elements of order $40$ must arise as preimages of order $20$ elements
under the natural epimorphism from $2.S$ to $S$,
which means that we may choose an order $40$ preimage $\hat{s}$ of $s$.
Then $\M(2.S, \hat{s})$ consists of central extensions of the subgroups
listed in statement~(b).
The perfect subgroups $O_7(3)$, $L_4(3)$, $2.U_4(3)$, and $U_4(2)$
of these groups must lift to their Schur double covers in $2.S$ because
otherwise the preimages would not contain elements of order $40$.

Next we consider the preimage of the subgroup $U = (A_4 \times U_4(2)).2$
of $S$.
We show that the preimages of the two direct factors $A_4$ and $U_4(2)$
in $U^{\prime} = A_4 \times U_4(2)$ are Schur covers.
For $A_4$, this follows from the fact that the preimage of $U^{\prime}$
must contain elements of order $20$,
and that $U_4(2)$ does not contain elements of order $10$.

\begin{verbatim}
    gap> u42:= CharacterTable( "U4(2)" );;
    gap> Filtered( OrdersClassRepresentatives( u42 ), x -> x mod 5 = 0 );
    [ 5 ]
\end{verbatim}

In order to show that the $U_4(2)$ type subgroup of $U^{\prime}$ lifts to its
double cover in $2.S$, we note that the class {\tt 2B} of $U_4(2)$ lifts
to a class of elements of order four in the double cover $2.U_4(2)$,
and that the corresponding class of elements in $U$
is $S$-conjugate to the class of involutions in the direct factor $A_4$
(which is the unique class of length three in $U$).

\begin{verbatim}
    gap> u:= CharacterTable( Maxes( t )[18] );
    CharacterTable( "(A4xU4(2)):2" )
    gap> 2u42:= CharacterTable( "2.U4(2)" );;
    gap> OrdersClassRepresentatives( 2u42 )[4];
    4
    gap> GetFusionMap( 2u42, u42 )[4];
    3
    gap> OrdersClassRepresentatives( u42 )[3];
    2
    gap> List( PossibleClassFusions( u42, u ), x -> x[3] );
    [ 8 ]
    gap> PositionsProperty( SizesConjugacyClasses( u ), x -> x = 3 );
    [ 2 ]
    gap> ForAll( PossibleClassFusions( u, t ), x -> x[2] = x[8] );
    true
\end{verbatim}


The last subgroup for which the structure of the preimage has to be shown
is $U = (A_6 \times A_6):2^2$.
We claim that each of the $A_6$ type subgroups in the derived subgroup
$U^{\prime} = A_6 \times A_6$ lifts to its double cover in $2.S$.
Since all elements of order $20$ in $U$ lie in $U^{\prime}$,
at least one of the two direct factors must lift to its double cover,
in order to give rise to an order $40$ element in $U$.
In fact both factors lift to the double cover
since the two direct factors are interchanged by conjugation in $U$;
the latter follows form tha fact that $U$ has no normal subgroup
of type $A_6$.

\begin{verbatim}
    gap> u:= CharacterTable( Maxes( t )[24] );
    CharacterTable( "(A6xA6):2^2" )
    gap> ClassPositionsOfDerivedSubgroup( u );
    [ 1 .. 22 ]
    gap> PositionsProperty( OrdersClassRepresentatives( u ), x -> x = 20 );
    [ 8 ]
    gap> List( ClassPositionsOfNormalSubgroups( u ),
    >          x -> Sum( SizesConjugacyClasses( u ){ x } ) );
    [ 1, 129600, 259200, 259200, 259200, 518400 ]
\end{verbatim}

So statement~(e) holds.

For statement~(f), we have to consider the upward extensions $S.2_1$,
$S.2_2$, and $S.3$.

First we look at $S.2_1$, an extension by an outer automorphism that
acts as a double transposition in the outer automorphism group $S_4$.
Note that the symmetry between the three classes of element oder $20$ in $S$
is broken in $S.2_1$,
two of these classes have square roots in $S.2_1$, the third has not.

\begin{verbatim}
    gap> t2:= CharacterTable( "O8+(3).2_1" );;
    gap> ord20:= PositionsProperty( OrdersClassRepresentatives( t2 ),
    >                x -> x = 20 );;
    gap> ord20:= Intersection( ord20, ClassPositionsOfDerivedSubgroup( t2 ) );
    [ 84, 85, 86 ]
    gap> List( ord20, x -> x in PowerMap( t2, 2 ) );
    [ false, true, true ]
\end{verbatim}

Changing the viewpoint, we see that for each class of element order $20$
in $S$,
there is a group of the type $S.2_1$ in which the elements in this class
do not have square roots,
and there are groups of this type in which these elements have square roots.
So we have to deal with two different cases,
and we do this by first collecting the permutation characters induced from
{\bf all} maximal subgroups of $S.2_1$ (other than $S$) that contain
elements of order $20$ in $S$,
and then considering $s$ in each of these classes of $S$.

We fix an embedding of $S$ into $S.2_1$ in which the elements in the class
{\tt 20A} do not have square roots.
This situation is given for the stored class fusion between the tables
in the {\GAP} Character Table Library.

\begin{verbatim}
    gap> tfust2:= GetFusionMap( t, t2 );;
    gap> tfust2{ PositionsProperty( OrdersClassRepresentatives( t ),
    >                x -> x = 20 ) };
    [ 84, 85, 86 ]
\end{verbatim}

The six different actions of $S$ on the cosets of $O_7(3)$ type subgroups
induce pairwise different permutation characters that form an orbit under
the action of $\Aut(S)$.
Four of these characters cannot extend to $S.2_1$,
the other two extend to permutation characters of $S.2_1$ on the cosets of
$O_7(3).2$ type subgroups;
these subgroups contain {\tt 20A} elements.

\begin{verbatim}
    gap> primt2:= [];;
    gap> poss:= PossiblePermutationCharacters( CharacterTable( "O7(3)" ), t );;
    gap> invfus:= InverseMap( tfust2 );;
    gap> List( poss, pi -> ForAll( CompositionMaps( pi, invfus ), IsInt ) );
    [ false, false, false, false, true, true ]
    gap> PossiblePermutationCharacters(
    >        CharacterTable( "O7(3)" ) * CharacterTable( "Cyclic", 2 ), t2 );
    [  ]
    gap> ext:= PossiblePermutationCharacters( CharacterTable( "O7(3).2" ), t2 );;
    gap> List( ext, pi -> pi{ ord20 } );
    [ [ 1, 0, 0 ], [ 1, 0, 0 ] ]
    gap> Append( primt2, ext );
\end{verbatim}

The novelties in $S.2_1$ that arise from $O_7(3)$ type subgroups of $S$
have the structure $L_4(3).2^2$.
These subgroups contain elements in the classes {\tt 20B} and {\tt 20C}
of $S$.

\begin{verbatim}
    gap> ext:= PossiblePermutationCharacters( CharacterTable( "L4(3).2^2" ), t2 );;
    gap> List( ext, pi -> pi{ ord20 } );
    [ [ 0, 0, 1 ], [ 0, 1, 0 ] ]
    gap> Append( primt2, ext );
\end{verbatim}

Note that from the possible permutation characters of $S.2_1$ on the cosets
of $L_4(3):2 \times 2$ type subgroups,
we see that such subgroups must contain {\tt 20A} elements,
i.~e., all such subgroups of $S.2_1$ lie inside $O_7(3).2$ type subgroups.
This means that the structure description of these novelties
in~\cite[p.~140]{CCN85} is not correct.
The correct structure is $L_4(3).2^2$.)

\begin{verbatim}
    gap> List( PossiblePermutationCharacters( CharacterTable( "L4(3).2_2" ) *
    >              CharacterTable( "Cyclic", 2 ), t2 ), pi -> pi{ ord20 } );
    [ [ 1, 0, 0 ] ]
\end{verbatim}

All $3^6:L_4(3)$ type subgroups of $S$ extend to $S.2_1$.
We compute these permutation characters as the possible permutation characters
of the right degree.

\begin{verbatim}
    gap> ext:= PermChars( t2, rec( torso:= [ 1120 ] ) );;
    gap> List( ext, pi -> pi{ ord20 } );
    [ [ 2, 0, 0 ], [ 0, 0, 2 ], [ 0, 2, 0 ] ]
    gap> Append( primt2, ext );
\end{verbatim}

Also all $2.U_4(3).2^2$ type subgroups of $S$ extend to $S.2_1$.
We compute the permutation characters as the extensions of the corresponding
permutation characters of $S$.

\begin{verbatim}
    gap> filt:= Filtered( prim, x -> x[1] = 189540 );;
    gap> cand:= List( filt, x -> CompositionMaps( x, invfus ) );;
    gap> ext:= Concatenation( List( cand,
    >              pi -> PermChars( t2, rec( torso:= pi ) ) ) );;
    gap> List( ext, x -> x{ ord20 } );
    [ [ 1, 0, 0 ], [ 0, 1, 0 ], [ 0, 0, 1 ] ]
    gap> Append( primt2, ext );
\end{verbatim}

The extensions of $(A_4 \times U_4(2)):2$ type subgroups of $S$ to $S.2_1$
have the type $S_4 \times U_4(2):2$, they contain {\tt 20A} elements.

\begin{verbatim}
    gap> ext:= PossiblePermutationCharacters( CharacterTable( "Symmetric", 4 ) *
    >              CharacterTable( "U4(2).2" ), t2 );;
    gap> List( ext, x -> x{ ord20 } );
    [ [ 1, 0, 0 ], [ 1, 0, 0 ] ]
    gap> Append( primt2, ext );
\end{verbatim}

All $(A_6 \times A_6):2^2$ type subgroups of $S$ extend to $S.2_1$.
We compute the permutation characters as the extensions of the corresponding
permutation characters of $S$.

\begin{verbatim}
    gap> filt:= Filtered( prim, x -> x[1] = 9552816 );;
    gap> cand:= List( filt, x -> CompositionMaps( x, InverseMap( tfust2 ) ));;
    gap> ext:= Concatenation( List( cand,
    >              pi -> PermChars( t2, rec( torso:= pi ) ) ) );;
    gap> List( ext, x -> x{ ord20 } );
    [ [ 1, 0, 0 ], [ 0, 1, 0 ], [ 0, 0, 1 ] ]
    gap> Append( primt2, ext );
\end{verbatim}

We have found all relevant permutation characters of $S.2_1$.
This together with the list in~\cite[p.~140]{CCN85} implies statement~(g).

Now we compute the bounds $\total^{\prime}(S.2_1, s)$.

\begin{verbatim}
    gap> Length( primt2 );
    15
    gap> approx:= List( ord20, x -> ApproxP( primt2, x ) );;
    gap> outer:= Difference(
    >      PositionsProperty( OrdersClassRepresentatives( t2 ), IsPrimeInt ),
    >      ClassPositionsOfDerivedSubgroup( t2 ) );;
    gap> List( approx, l -> Maximum( l{ outer } ) );
    [ 574/1215, 83/567, 83/567 ]
\end{verbatim}

Next we look at $S.2_2$, an extension by an outer automorphism that
acts as a transposition in the outer automorphism group $S_4$.
Similar to the above situation,
the symmetry between the three classes of element oder $20$ in $S$
is broken also in $S.2_2$:
The first is a conjugacy class of $S.2_2$,
the other two classes are fused in $S.2_2$,

\begin{verbatim}
    gap> t2:= CharacterTable( "O8+(3).2_2" );;
    gap> ord20:= PositionsProperty( OrdersClassRepresentatives( t2 ),
    >                x -> x = 20 );;
    gap> ord20:= Intersection( ord20, ClassPositionsOfDerivedSubgroup( t2 ) );
    [ 82, 83 ]
    gap> tfust2:= GetFusionMap( t, t2 );;
    gap> tfust2{ PositionsProperty( OrdersClassRepresentatives( t ),
    >                x -> x = 20 ) };
    [ 82, 83, 83 ]
\end{verbatim}

Like in the case $S.2_1$, we compute the permutation characters induced from
{\bf all} maximal subgroups of $S.2_2$ (other than $S$) that contain
elements of order $20$ in $S$.

We fix the embedding of $S$ into $S.2_2$ in which the class
{\tt 20A} of $S$ is a class of $S.2_2$.
This situation is given for the stored class fusion between the tables
in the {\GAP} Character Table Library.

Exactly two classes of $O_7(3)$ type subgroups in $S$ extend to $S.2_2$,
these groups contain {\tt 20A} elements.

\begin{verbatim}
    gap> primt2:= [];;
    gap> ext:= PermChars( t2, rec( torso:= [ 1080 ] ) );;
    gap> List( ext, pi -> pi{ ord20 } );
    [ [ 1, 0 ], [ 1, 0 ] ]
    gap> Append( primt2, ext );
\end{verbatim}

Only one class of $3^6:L_4(3)$ type subgroups extends to $S.2_2$.
(Note that we need not consider the novelties of the type
$3^{3+6}:(L_3(3) \times 2)$, because the order of these groups is not
divisible by $5$.)

\begin{verbatim}
    gap> ext:= PermChars( t2, rec( torso:= [ 1120 ] ) );;
    gap> List( ext, pi -> pi{ ord20 } );
    [ [ 2, 0 ] ]
    gap> Append( primt2, ext );
\end{verbatim}

Only one class of $2.U_4(3).2^2$ type subgroups of $S$ extends to $S.2_2$.
We compute the permutation character as the extension of the corresponding
permutation characters of $S$.

\begin{verbatim}
    gap> filt:= Filtered( prim, x -> x[1] = 189540 );;
    gap> cand:= List( filt, x -> CompositionMaps( x, InverseMap( tfust2 ) ));;
    gap> ext:= Concatenation( List( cand,
    >              pi -> PermChars( t2, rec( torso:= pi ) ) ) );;
    gap> List( ext, x -> x{ ord20 } );
    [ [ 1, 0 ] ]
    gap> Append( primt2, ext );
\end{verbatim}

Two classes of $(A_4 \times U_4(2)):2$ type subgroups of $S$ extend
to $S.2_2$.

\begin{verbatim}
    gap> filt:= Filtered( prim, x -> x[1] = 7960680 );;
    gap> cand:= List( filt, x -> CompositionMaps( x, InverseMap( tfust2 ) ));;
    gap> ext:= Concatenation( List( cand,
    >              pi -> PermChars( t2, rec( torso:= pi ) ) ) );;
    gap> List( ext, x -> x{ ord20 } );
    [ [ 1, 0 ], [ 1, 0 ] ]
    gap> Append( primt2, ext );
\end{verbatim}

Exactly one class of $(A_6 \times A_6):2^2$ type subgroups in $S$
extends to $S.2_2$, and the extensions have the structure
\tthdump{$S_6 \wr 2$.}

\begin{verbatim}
    gap> ext:= PossiblePermutationCharacters( CharacterTableWreathSymmetric(
    >              CharacterTable( "S6" ), 2 ), t2 );;
    gap> List( ext, x -> x{ ord20 } );
    [ [ 1, 0 ] ]
    gap> Append( primt2, ext );
\end{verbatim}

We have found all relevant permutation characters of $S.2_2$,
and compute the bounds $\total^{\prime}(S.2_2, s)$.

\begin{verbatim}
    gap> Length( primt2 );
    7
    gap> approx:= List( ord20, x -> ApproxP( primt2, x ) );;
    gap> outer:= Difference(
    >      PositionsProperty( OrdersClassRepresentatives( t2 ), IsPrimeInt ),
    >      ClassPositionsOfDerivedSubgroup( t2 ) );;
    gap> List( approx, l -> Maximum( l{ outer } ) );
    [ 14/9, 0 ]
\end{verbatim}

This means that there is an extension of the type $S.2_2$ in which
$s$ cannot be chosen such that the bound is less than $1/2$.
More precisely, we have $\total(g,s) \geq 1/2$ exactly for $g$ in the
unique outer involution class of size $1\,080$.

\begin{verbatim}
    gap> approx:= ApproxP( primt2, ord20[1] );;
    gap> bad:= Filtered( outer, i -> approx[i] >= 1/2 );
    [ 84 ]
    gap> OrdersClassRepresentatives( t2 ){ bad };
    [ 2 ]
    gap> SizesConjugacyClasses( t2 ){ bad };
    [ 1080 ]
    gap> Number( SizesConjugacyClasses( t2 ), x -> x = 1080 );
    1
\end{verbatim}

So we compute the proportion of elements in this class that generate $S.2_2$
together with an element $s$ of order $20$ in $S$.
(As above, we have to consider two conjugacy classes.)
For that, we first compute a permutation representation of $S.2_2$,
using that $S.2_2$ is isomporphic to the two subgroups of index $2$ in
$\PGO^+(8,3) = O_8^+(3).2^2_{122}$ that are different from
$\PSO^+(8,3) = O_8^+(3).2_1$, cf.~\cite[p.~140]{CCN85}.

\begin{verbatim}
    gap> go:= GO(1,8,3);;
    gap> so:= SO(1,8,3);;
    gap> outerelm:= First( GeneratorsOfGroup( go ), x -> not x in so );;
    gap> g2:= ClosureGroup( DerivedSubgroup( so ), outerelm );;
    gap> Size( g2 );
    19808719257600
    gap> dom:= NormedVectors( GF(3)^8 );;
    gap> orbs:= Orbits( g2, dom, OnLines );;
    gap> List( orbs, Length );
    [ 1080, 1080, 1120 ]
    gap> act:= Action( g2, orbs[1], OnLines );;
\end{verbatim}

An involution $g$ can be found as a power of one of the given generators.

\begin{verbatim}
    gap> Order( outerelm );
    26
    gap> g:= Permutation( outerelm^13, orbs[1], OnLines );;
    gap> Size( ConjugacyClass( act, g ) );
    1080
\end{verbatim}

Now we find the candidates for the elements $s$,
and compute their ratios of nongeneration.

\begin{verbatim}
    gap> ord20;
    [ 82, 83 ]
    gap> SizesCentralizers( t2 ){ ord20 };
    [ 40, 20 ]
    gap> der:= DerivedSubgroup( act );;
    gap> repeat 20A:= Random( der );
    >    until Order( 20A ) = 20 and Size( Centralizer( act, 20A ) ) = 40;
    gap> RatioOfNongenerationTransPermGroup( act, g, 20A );
    1
    gap> repeat 20BC:= Random( der );
    >    until Order( 20BC ) = 20 and Size( Centralizer( act, 20BC ) ) = 20;
    gap> RatioOfNongenerationTransPermGroup( act, g, 20BC );
    0
\end{verbatim}

This means that for $s$ in one $S$-class of elements of order $20$,
we have $\prop^{\prime}(g, s) = 1$,
and $s$ in the other two $S$-classes of elements of order $20$ generates
with any conjugate of $g$.

Concerning $S.2_2$, it remains to show that we cannot find a better element
than $s$.
For that, we first compute class representatives $s^{\prime}$
in $S$, w.r.t.~conjugacy in $S.2_2$,
and then compute $\prop^{\prime}( s^{\prime}, g )$.
(It would be enough to check representatives of classes of maximal element
order, but computing all classes is easy enough.)

\begin{verbatim}
    gap> ccl:= ConjugacyClasses( act );;
    gap> der:= DerivedSubgroup( act );;
    gap> reps:= Filtered( List( ccl, Representative ), x -> x in der );;
    gap> Length( reps );
    83
    gap> ratios:= List( reps,
    >                   s -> RatioOfNongenerationTransPermGroup( act, g, s ) );;
    gap> cand:= PositionsProperty( ratios, x -> x < 1 );;
    gap> ratios:= ratios{ cand };;
    gap> SortParallel( ratios, cand );
    gap> ratios;
    [ 0, 1/10, 1/10, 16/135, 1/3, 1/3, 11/27, 7/15, 7/15 ]
\end{verbatim}

For $S.2_2$, it remains to show that there is no element $s^{\prime} \in S$
such that $\prop^{\prime}( {s^{\prime}}^x, g ) < 1$ holds for any
$x \in \Aut(S)$ and $g \in S.2_2$.
So we are done when we can show that each class given by \verb|cand|
is conjugate in $S.3$ to a class outside \verb|cand|.
The classes can be identified by element orders and centralizer orders.

\begin{verbatim}
    gap> invs:= List( cand,
    >       x -> [ Order( reps[x] ), Size( Centralizer( der, reps[x] ) ) ] );
    [ [ 20, 20 ], [ 18, 108 ], [ 18, 108 ], [ 14, 28 ], [ 15, 45 ], [ 15, 45 ], 
      [ 10, 40 ], [ 12, 72 ], [ 12, 72 ] ]
\end{verbatim}

Namely, \verb|cand| contains no full $S.3$-orbit of classes of the element orders
$20$, $18$, $14$, $15$, and $10$; also, \verb|cand| does not contain full
$S.3$-orbits on the classes {\tt 12O}--{\tt 12T}.

Finally, we deal with $S.3$.
The fact that no maximal subgroup of $S$ containing an element of order $20$
extends to $S.3$ follows either from the list of maximal subgroups of $S$
in~\cite[p.~140]{CCN85}
or directly from the permutation characters.

\begin{verbatim}
    gap> t3:= CharacterTable( "O8+(3).3" );;
    gap> tfust3:= GetFusionMap( t, t3 );;
    gap> inv:= InverseMap( tfust3 );;
    gap> filt:= PositionsProperty( prim, x -> x[ spos ] <> 0 );;
    gap> ForAll( prim{ filt },
    >            pi -> ForAny( CompositionMaps( pi, inv ), IsList ) );
    true
\end{verbatim}

So we have to consider only the classes of novelties in $S.3$,
but the order of none of these groups is divisible by $20$
--again see~\cite[p.~140]{CCN85}).
This means that {\bf any} element in $S.3 \setminus S$ together with an element
of order $20$ in $S$ generates $S.3$.
This is in fact stronger than statement~(f),
which claims this property only for elements of prime order in
$S.3 \setminus S$ (and their roots);
note that $S.3 \setminus S$ contains elements of the orders $9$ and $27$.

\begin{verbatim}
    gap> outer:= Difference( [ 1 .. NrConjugacyClasses( t3 ) ],
    >                ClassPositionsOfDerivedSubgroup( t3 ) );
    [ 53, 54, 55, 56, 57, 58, 59, 60, 61, 62, 63, 64, 65, 66, 67, 68, 69, 70, 71, 
      72, 73, 74, 75, 76, 77, 78, 79, 80, 81, 82, 83, 84, 85, 86, 87, 88, 89, 90, 
      91, 92, 93, 94 ]
    gap> Set( OrdersClassRepresentatives( t3 ){ outer } );      
    [ 3, 6, 9, 12, 18, 21, 24, 27, 36, 39 ]
\end{verbatim}

Before we turn to the next computations, we clean the workspace.

\begin{verbatim}
    gap> CleanWorkspace();
\end{verbatim}

\subsection{$O^+_8(4)$}\label{O8p4}

We show that $S = O^+_8(4) = \Omega^+(8,4)$ satisfies the following.
\begin{enumerate}
\item[(a)]
    For suitable $s \in S$ of the type $2^- \perp 6^-$
    (i.~e., $s$ decomposes the natural $8$-dimensional module for
    $S$ into an orthogonal sum of two irreducible modules
    of the dimensions $2$ and $6$, respectively) and of order $65$,
    $\M(S,s)$ consists of exactly three pairwise nonconjugate subgroups
    of the type $(5 \times O^-_6(4)).2 = (5 \times \Omega^-(6,4)).2$.
\item[(b)]
    $\total( S, s ) \leq 34\,817 / 1\,645\,056$.
\item[(c)]
    In the extensions $S.2_1$ and $S.3$ of $S$ by graph automorphisms,
    there is at most one maximal subgroup besides $S$ that contains $s$.
    For the extension $S.2_2$ of $S$ by a field automorphism,
    we have $\total^{\prime}(S.2_2, s) = 0$.
    In the extension $S.2_3$ of $S$ by the product of an involutory graph
    automorphism and a field automorphism,
    there is a unique maximal subgroup besides $S$ that contains $s$.
\end{enumerate}

A safe source for determining $\M(S,s)$ is~\cite{Kle87}.
By inspection of the result matrix in this paper,
we get that the only maximal subgroups of $S$ that contain elements
of order $65$ occur in the rows 9--14 and 23--25;
they have the isomorphism types
$S_6(4) = \Sp(6,4) \cong O_7(4) = \Omega(7,4)$ and
$(5 \times O_6^-(4)).2 = (5 \times \Omega^-(6,4)).2$, respectively,
and for each of these, there are three conjugacy classes of subgroups
in $S$, which are conjugate under the triality graph automorphism of $S$.

We start with the natural matrix representation of $S$.
For convenience,
we compute an isomorphic permutation group on $5\,525$ points.

\begin{verbatim}
    gap> q:= 4;;  n:= 8;;
    gap> G:= DerivedSubgroup( SO( 1, n, q ) );;
    gap> points:= NormedRowVectors( GF(q)^n );;
    gap> orbs:= Orbits( G, points, OnLines );;
    gap> List( orbs, Length );
    [ 5525, 16320 ]
    gap> hom:= ActionHomomorphism( G, orbs[1], OnLines );;
    gap> G:= Image( hom );;
\end{verbatim}

The group $S$ contains exactly six conjugacy classes of (cyclic) subgroups
of order $65$;
this follows from the fact that the centralizer of any Sylow $13$ subgroup
in $S$ has the structure $5 \times 5 \times 13$.

\begin{verbatim}
    gap> Collected( Factors( Size( G ) ) );
    [ [ 2, 24 ], [ 3, 5 ], [ 5, 4 ], [ 7, 1 ], [ 13, 1 ], [ 17, 2 ] ]
    gap> ResetGlobalRandomNumberGenerators();
    gap> repeat x:= Random( G );
    >    until Order( x ) mod 13 = 0;
    gap> x:= x^( Order( x ) / 13 );;
    gap> c:= Centralizer( G, x );;
    gap> IsAbelian( c );  AbelianInvariants( c );
    true
    [ 5, 5, 13 ]
\end{verbatim}

The group $S_6(4)$ contains exactly one class of subgroups of order $65$,
since the conjugacy classes of elements of order $65$ in $S_6(4)$ are
algebraically conjugate.

\begin{verbatim}
    gap> t:= CharacterTable( "S6(4)" );;
    gap> ord65:= PositionsProperty( OrdersClassRepresentatives( t ),
    >                               x -> x = 65 );
    [ 105, 106, 107, 108, 109, 110, 111, 112 ]
    gap> ord65 = ClassOrbit( t, ord65[1] );
    true
\end{verbatim}

Thus there are at least three classes of order $65$ elements in $S$
that are \emph{not} contained in $S_6(4)$ type subgroups of $S$.
So we choose such an element $s$,
and have to consider only overgroups of the type
$(5 \times \Omega^-(6,4)).2$.

The group $\Omega^-(6,4) \cong U_4(4)$ contains exactly one class
of subgroups of order $65$.

\begin{verbatim}
    gap> t:= CharacterTable( "U4(4)" );;
    gap> ords:= OrdersClassRepresentatives( t );;
    gap> ord65:= PositionsProperty( ords, x -> x = 65 );;
    gap> ord65 = ClassOrbit( t, ord65[1] );
    true
\end{verbatim}

So $5 \times \Omega^-(6,4)$ contains exactly six such classes.
Furthermore, subgroups in different classes are not $S$-conjugate.

\begin{verbatim}
    gap> syl5:= SylowSubgroup( c, 5 );;
    gap> elms:= Filtered( Elements( syl5 ), y -> Order( y ) = 5 );;
    gap> reps:= Set( List( elms, SmallestGeneratorPerm ) );;  Length( reps );
    6
    gap> reps65:= List( reps, y -> SubgroupNC( G, [ y * x ] ) );;
    gap> pairs:= Filtered( UnorderedTuples( [ 1 .. 6 ], 2 ),
    >                      p -> p[1] <> p[2] );;
    gap> ForAny( pairs, p -> IsConjugate( G, reps65[ p[1] ], reps65[ p[2] ] ) );
    false
\end{verbatim}

We consider only subgroups $M \leq S$ in the three $S$-classes of the type
$(5 \times \Omega^-(6,4)).2$.

\begin{verbatim}
    gap> cand:= List( reps, y -> Normalizer( G, SubgroupNC( G, [ y ] ) ) );;
    gap> cand:= Filtered( cand, y -> Size( y ) = 10 * Size( t ) );;
    gap> Length( cand );
    3
\end{verbatim}

%
%
%
%

(Note that one of the members in $\M(S,s)$ is the stabilizer in $S$ of the
orthogonal decomposition $2^- \perp 6^-$,
the other two members are not reducible.)

By the above, the classes of subgroups of order $65$ in each such $M$
are in bijection with the corresponding classes in $S$.
Since $N_S(\langle g \rangle) \subseteq M$ holds
for any $g \in M$ of order $65$, also the conjugacy classes of
\emph{elements} of order $65$ in $M$ are in bijection with those in $S$.

\begin{verbatim}
    gap> norms:= List( reps65, y -> Normalizer( G, y ) );;
    gap> ForAll( norms, y -> ForAll( cand, M -> IsSubset( M, y ) ) );
    true
\end{verbatim}

As a consequence, we have $g^S \cap M = g^M$ and thus $1_M^S(g) = 1$.
This implies statement~(a).

In order to show statement~(b),
we want to use the function \verb|UpperBoundFixedPointRatios| introduced in
Section~\ref{groups}.
For that, we first compute the conjugacy classes of
the three class representatives $M$.
(Since the groups have elementary abelian Sylow $5$ subgroups of the order
$5^4$, computing all conjugacy classes appears to be faster than using
\verb|ClassesOfPrimeOrder|.)
Then we compute an upper bounds for $\total(S,s)$.

\begin{verbatim}
    gap> syl5:= SylowSubgroup( cand[1], 5 );;
    gap> Size( syl5 );  IsElementaryAbelian( syl5 );
    625
    true
    gap> UpperBoundFixedPointRatios( G, List( cand, ConjugacyClasses ), false );
    [ 34817/1645056, false ]
\end{verbatim}

\begin{remark}
{\rm
Computing the exact value $\total(S,s)$ in the above setup would require
to test the $S$-conjugacy of certain order $5$ elements in $M$.
With the current {\GAP} implementation,
some of the relevant tests need several hours of CPU time.

An alternative approach would be to compute the permutation action of
$S$ on the cosets of $M$, of degree $6\,580\,224$,
and to count the fixed points of conjugacy class representatives
of prime order.
The currently available {\GAP} library methods are not sufficient
for computing this in reasonable time.
``Ad-hoc code'' for this special case works,
but it seemed to be not appropriate to include it here.}
\end{remark}

In the proof of statement~(c),
again we consult the result matrix in~\cite{Kle87}.
For $S.3$,
the maximal subgroups are in the rows $4$, $15$, $22$, $26$, and $61$.
Only row $26$ yields subgroups that contain elements $s$ of order $65$,
they have the isomorphism type
$(5 \times \GU(3,4)).2 \cong (5^2 \times U_3(4)).2$.
Note that the conjugacy classes of the members in $\M(S,s)$ are permuted
by the outer automorphism of order $3$,
so none of the subgroups in $\M(S,s)$ extends to $S.3$.
By~\cite[Lemma~2.4~(2)]{BGK}, if there is a maximal subgroup of $S.3$
besides $S$ that contains $s$ then this subgroup is the normalizer in $S.3$
of the intersection of the three members of $\M(S,s)$,
i.~e., $s$ is contained in at most one such subgroup.

For $S.2_1$,
only the rows $9$ and $23$ yield maximal subgroups containing elements of
order $65$, and since we had chosen $s$ in such a way
that row $9$ was excluded already for the simple group,
only extensions of the elements in $\M(S,s)$ can appear.
Exactly one of these three subgroups of $S$ extends to $S.2_1$,
so again we get just one maximal subgroup of $S.2_1$, besides $S$,
that contains $s$.

All subgroups in $\M(S,s)$ extend to $S.2_2$, see~\cite{Kle87}.
We compute the extensions of the above subgroups $M$ of $S$ to $S.2_2$,
by constructing the action of the field automorphism
in the permutation representation we used for $S$.
In other words, we compute the projective action of the Frobenius map.

\begin{verbatim}
    gap> frob:= PermList( List( orbs[1], v -> Position( orbs[1],
    >              List( v, x -> x^2 ) ) ) );;
    gap> G2:= ClosureGroupDefault( G, frob );;
    gap> cand2:= List( cand, M -> Normalizer( G2, M ) );;
    gap> ccl:= List( cand2,
    >                M2 -> PcConjugacyClassReps( SylowSubgroup( M2, 2 ) ) );;
    gap> List( ccl, l -> Number( l, x -> Order( x ) = 2 and not x in G ) );
    [ 0, 0, 0 ]
\end{verbatim}

So in each case, the extension of $M$ to its normalizer in $S.2_2$
is non-split.
This implies $\total^{\prime}(S.2_2,s) = 0$.

Finally, in the extension of $S$ by the product of a graph automorphism
and the field automorphism, exactly that member of $\M(S,s)$ is invariant
that is invariant under the graph automorphism,
hence statement~(c) holds.

It is again time to clean the workspace.

\begin{verbatim}
    gap> CleanWorkspace();
\end{verbatim}

\subsection{$\ast$~$O_9(3)$}\label{O93}

The group $S = O_9(3) = \Omega_9(3)$ is the first member in the series
dealt with in~\cite[Proposition~5.7]{BGK},
and serves as an example to illustrate this statement.

\begin{enumerate}
\item[(a)] 
    For $s \in S$ of the type $1 \perp 8^-$
    (i.~e., $s$ decomposes the natural $9$-dimensional module for $S$
    into an orthogonal sum of two irreducible modules of the dimensions
    $1$ and $8$, respectively) and of order $(3^4 + 1)/2 = 41$,
    $\M(S,s)$ consists of one group of the type $O_8^-(3).2_1 = \PGO^-(8,3)$.
\item[(b)]
    $\total(S,s) = 1/3$.
\item[(c)]
    The uniform spread of $S$ is at least three,
    with $s$ of order $41$.
\end{enumerate}

By~\cite{MSW94}, the only maximal subgroup of $S$ that contains $s$
is the stabilizer $M$ of the orthogonal decomposition.
The group $2 \times O_8^-(3).2_1 = \GO^-(8,3)$ embeds naturally into
$\SO(9,3)$, its intersection with $S$ is $\PGO^-(8,3)$.
This proves statement~(a).

The group $M$ is the stabilizer of a $1$-space, it has index $3\,240$ in $S$.

\begin{verbatim}
    gap> g:= SO( 9, 3 );;
    gap> g:= DerivedSubgroup( g );;
    gap> Size( g );
    65784756654489600
    gap> orbs:= Orbits( g, NormedRowVectors( GF(3)^9 ), OnLines );;
    gap> List( orbs, Length ) / 41;
    [ 3240/41, 81, 80 ]
    gap> Size( SO( 9, 3 ) ) / Size( GO( -1, 8, 3 ) );
    3240
\end{verbatim}

So we compute the unique transitive permutation character of $S$ that has
degree $3\,240$.

\begin{verbatim}
    gap> t:= CharacterTable( "O9(3)" );;
    gap> pi:= PermChars( t, rec( torso:= [ 3240 ] ) );
    [ Character( CharacterTable( "O9(3)" ), [ 3240, 1080, 380, 132, 48, 324, 378,
          351, 0, 0, 54, 27, 54, 27, 0, 118, 0, 36, 46, 18, 12, 2, 8, 45, 0, 108,
          108, 135, 126, 0, 0, 56, 0, 0, 36, 47, 38, 27, 39, 36, 24, 12, 18, 18,
          15, 24, 2, 18, 15, 9, 0, 0, 0, 2, 0, 18, 11, 3, 9, 6, 6, 9, 6, 3, 6, 3,
          0, 6, 16, 0, 4, 6, 2, 45, 36, 0, 0, 0, 0, 0, 0, 0, 9, 9, 6, 3, 0, 0,
          15, 13, 0, 5, 7, 36, 0, 10, 0, 10, 19, 6, 15, 0, 0, 0, 0, 12, 3, 10, 0,
          3, 3, 7, 0, 6, 6, 2, 8, 0, 4, 0, 2, 0, 1, 3, 0, 0, 3, 0, 3, 2, 2, 3, 3,
          6, 2, 2, 9, 6, 3, 0, 0, 18, 9, 0, 0, 12, 0, 0, 8, 0, 6, 9, 5, 0, 0, 0,
          0, 0, 0, 0, 0, 3, 3, 3, 2, 1, 3, 3, 1, 0, 0, 4, 1, 0, 0, 1, 0, 3, 3, 1,
          1, 2, 2, 0, 0, 1, 3, 4, 0, 1, 2, 0, 0, 1, 0, 4, 1, 0, 0, 0, 0, 1, 0, 0,
          1, 0, 0, 1, 1, 1, 1, 1, 0, 0, 1, 1, 1, 0 ] ) ]
    gap> spos:= Position( OrdersClassRepresentatives( t ), 41 );
    208
    gap> approx:= ApproxP( pi, spos );;
    gap> Maximum( approx );
    1/3
    gap> PositionsProperty( approx, x -> x = 1/3 );
    [ 2 ]
    gap> SizesConjugacyClasses( t )[2];
    3321
    gap> OrdersClassRepresentatives( t )[2];
    2
\end{verbatim}

We see that $\prop( S, s ) = \total( S, s ) = 1/3$ holds,
and that $\total( g, s )$ attains this maximum only for $g$ in one
class of involutions in $S$;
let us call this class \verb|2A|.
(This class consists of the negatives of a class of \emph{reflections}
in $\GO(9,3)$.)
This shows statement~(b).

In order to show that the uniform spread of $S$ is at least three,
it suffices to show that for each triple of \verb|2A| elements,
there is an element $s$ of order $41$ in $S$ that generates $S$ with each
element of the triple.

We work with the primitive permutation representation of $S$ on $3\,240$
points.
In this representation, $s$ fixes exactly one point,
and by statement~(a),
$s$ generates $S$ with $x \in S$ if and only if $x$ moves this point.
Since the number of fixed points of each \verb|2A| involution in $S$ is
exactly one third of the moved points of $S$,
it suffices to show that we cannot choose three such involutions with
mutually disjoint fixed point sets.
And this is shown particularly easily because it will turn out that
already for any two different \verb|2A| involutions,
the sets of fixed points of are never disjoint.

First we compute a \verb|2A| element, which is determined as an involution
with exactly $1\,080$ fixed points.

\begin{verbatim}
    gap> g:= Action( g, orbs[1], OnLines );;
    gap> repeat
    >      repeat x:= Random( g ); ord:= Order( x ); until ord mod 2 = 0;
    >      y:= x^(ord/2);
    > until NrMovedPoints( y ) = 3240 - 1080;
\end{verbatim}

Next we compute the sets of fixed points of the elements in the class \verb|2A|,
by forming the $S$-orbit of the set of fixed points of the chosen
\verb|2A| element.

\begin{verbatim}
    gap> fp:= Difference( MovedPoints( g ), MovedPoints( y ) );;
    gap> orb:= Orbit( g, fp, OnSets );;
\end{verbatim}

Finally, we show that for any pair of \verb|2A| elements, their sets of
fixed points intersect nontrivially.
(Of course we can fix one of the two elements.)
This proves statement~(c).

\begin{verbatim}
    gap> ForAny( orb, l -> IsEmpty( Intersection( l, fp ) ) );
    false
\end{verbatim}


\subsection{$O_{10}^-(3)$}\label{O10m3}

We show that the group $S = O_{10}^-(3) = \POmega^-(10,3)$ satisfies the
following.
\begin{enumerate}
\item[(a)]
    For $s \in S$ irreducible of order $(3^5 + 1)/4 = 61$,
    $\M(S,s)$ consists of one subgroup of the type $\SU(5,3) \cong U_5(3)$.
\item[(b)]
    $\total(S,s) = 1/1\,066$.
\end{enumerate}

By~\cite{Be00}, the maximal subgroups of $S$ containing $s$ are of
extension field type,
and by~\cite[Prop.~4.3.18 and~4.3.20]{KlL90}, these groups have the structure
$\SU(5,3) = U_5(3)$
(which lift to $2 \times U_5(3) < \GU(5,3)$ in $\Omega^-(10,3) = 2.S$)
or $\Omega(5,9).2$,
but the order of the latter group is not divisible by $|s|$.
Furthermore, by~\cite[Lemma~2.12~(b)]{BGK}, $s$ is contained in only one
member of the former class.

\begin{verbatim}
    gap> Size( GO(5,9) ) / 61;
    6886425600/61
\end{verbatim}

\emph{(When the first version of these computations was written,
the character tables of both $S$ and $U_5(3)$ were not
contained in the {\GAP} Character Table Library,
so we worked with the groups.
Meanwhile the character tables are available,
thus we can show also a character theoretic solution.)}

\begin{verbatim}
    gap> t:= CharacterTable( "O10-(3)" );  s:= CharacterTable( "U5(3)" );
    CharacterTable( "O10-(3)" )
    CharacterTable( "U5(3)" )
    gap> SigmaFromMaxes( t, "61A", [ s ], [ 1 ] );
    1/1066
\end{verbatim}

\emph{(Now follow the computations with groups.)}

The first step is the construction of the embedding of $M = \SU(5,3)$
into the matrix group $2.S$,
that is, we write the matrix generators of $M$ as linear mappings on
the natural module for $2.S$, and then conjugate them such that the
result matrices respect the bilinear form of $2.S$.
For convenience, we choose a basis for the field extension
$\F_9/\F_3$ such that the $\F_3$-linear mapping given by the invariant
form of $M$ is invariant under the $\F_3$-linear mappings given by
the generators of $M$.


\begin{verbatim}
    gap> m:= SU(5,3);;
    gap> so:= SO(-1,10,3);;
    gap> omega:= DerivedSubgroup( so );;
    gap> om:= InvariantBilinearForm( so ).matrix;;
    gap> Display( om );
     . 1 . . . . . . . .
     1 . . . . . . . . .
     . . 1 . . . . . . .
     . . . 2 . . . . . .
     . . . . 2 . . . . .
     . . . . . 2 . . . .
     . . . . . . 2 . . .
     . . . . . . . 2 . .
     . . . . . . . . 2 .
     . . . . . . . . . 2
    gap> b:= Basis( GF(9), [ Z(3)^0, Z(3^2)^2 ] );
    Basis( GF(3^2), [ Z(3)^0, Z(3^2)^2 ] )
    gap> blow:= List( GeneratorsOfGroup( m ), x -> BlownUpMat( b, x ) );;
    gap> form:= BlownUpMat( b, InvariantSesquilinearForm( m ).matrix );;
    gap> ForAll( blow, x -> x * form * TransposedMat( x ) = form );
    true
    gap> Display( form );
     . . . . . . . . 1 .
     . . . . . . . . . 1
     . . . . . . 1 . . .
     . . . . . . . 1 . .
     . . . . 1 . . . . .
     . . . . . 1 . . . .
     . . 1 . . . . . . .
     . . . 1 . . . . . .
     1 . . . . . . . . .
     . 1 . . . . . . . .
\end{verbatim}

The matrix {\tt om} of the invariant bilinear form of $2.S$
is equivalent to the identity matrix $I$.
So we compute matrices {\tt T1} and {\tt T2}
that transform {\tt om} and {\tt form}, respectively, to $\pm I$.

\begin{verbatim}
    gap> T1:= IdentityMat( 10, GF(3) );;
    gap> T1{[1..3]}{[1..3]}:= [[1,1,0],[1,-1,1],[1,-1,-1]]*Z(3)^0;;
    gap> pi:= PermutationMat( (1,10)(3,8), 10, GF(3) );;
    gap> tr:= NullMat( 10,10,GF(3) );;
    gap> tr{[1, 2]}{[1, 2]}:= [[1,1],[1,-1]]*Z(3)^0;;
    gap> tr{[3, 4]}{[3, 4]}:= [[1,1],[1,-1]]*Z(3)^0;;
    gap> tr{[7, 8]}{[7, 8]}:= [[1,1],[1,-1]]*Z(3)^0;;
    gap> tr{[9,10]}{[9,10]}:= [[1,1],[1,-1]]*Z(3)^0;;
    gap> tr{[5, 6]}{[5, 6]}:= [[1,0],[0,1]]*Z(3)^0;;
    gap> tr2:= IdentityMat( 10,GF(3) );;
    gap> tr2{[1,3]}{[1,3]}:= [[-1,1],[1,1]]*Z(3)^0;;
    gap> tr2{[7,9]}{[7,9]}:= [[-1,1],[1,1]]*Z(3)^0;;
    gap> T2:= tr2 * tr * pi;;
    gap> D:= T1^-1 * T2;;
    gap> tblow:= List( blow, x -> D * x * D^-1 );;
    gap> IsSubset( omega, tblow );
    true
\end{verbatim}

Now we switch to a permutation representation of $S$,
and use the embedding of $M$ into $2.S$ to obtain the corresponding
subgroup of type $M$ in $S$.
Then we compute an upper bound for $\max\{ \fpr(g,S/M); g \in S^{\times} \}$.

\begin{verbatim}
    gap> orbs:= Orbits( omega, NormedRowVectors( GF(3)^10 ), OnLines );;
    gap> List( orbs, Length );
    [ 9882, 9882, 9760 ]
    gap> permgrp:= Action( omega, orbs[3], OnLines );;
    gap> M:= SubgroupNC( permgrp,
    >            List( tblow, x -> Permutation( x, orbs[3], OnLines ) ) );;
    gap> ccl:= ClassesOfPrimeOrder( M, Set( Factors( Size( M ) ) ),
    >                               TrivialSubgroup( M ) );;
    gap> UpperBoundFixedPointRatios( permgrp, [ ccl ], false );
    [ 1/1066, true ]
\end{verbatim}

The entry \verb|true| in the second position of the result indicates
that in fact the \emph{exact} value for the maximum of $\fpr(g,S/M)$
has been computed.
This implies statement~(b).

We clean the workspace.

\begin{verbatim}
    gap> CleanWorkspace();
\end{verbatim}

\subsection{$O_{14}^-(2)$}\label{O14m2}

We show that the group $S = O_{14}^-(2) = \Omega^-(14,2)$ satisfies the
following.
\begin{enumerate}
\item[(a)]
    For $s \in S$ irreducible of order $2^7+1 = 129$,
    $\M(S,s)$ consists of one subgroup $M$ of the type
    $\GU(7,2) \cong 3 \times U_7(2)$.
\item[(b)]
    $\total(S,s) = 1/2\,015$.
\end{enumerate}

By~\cite{Be00}, any maximal subgroup of $S$ containing $s$ is of extension
field type, and by~\cite[Table~3.5.F, Prop.~4.3.18]{KlL90},
these groups have the type $\GU(7,2)$,
and there is exactly one class of subgroups of this type.
Furthermore, by~\cite[Lemma~2.12~(a)]{BGK}, $s$ is contained in only one
member of this class.

We embed $U_7(2)$ into $S$,
by first replacing each element in $\F_4$ by the $2 \times 2$ matrix
of the induced $\F_2$-linear mapping w.r.t.~a suitable basis,
and then conjugating the images of the generators such that the invariant
quadratic form of $S$ is respected.

\begin{verbatim}
    gap> o:= SO(-1,14,2);;
    gap> g:= SU(7,2);;
    gap> b:= Basis( GF(4) );;
    gap> blow:= List( GeneratorsOfGroup( g ), x -> BlownUpMat( b, x ) );;
    gap> form:= NullMat( 14, 14, GF(2) );;
    gap> for i in [ 1 .. 14 ] do form[i][ 15-i ]:= Z(2); od;
    gap> ForAll( blow, x -> x * form * TransposedMat( x ) = form );
    true
    gap> pi:= PermutationMat( (1,13)(3,11)(5,9), 14, GF(2) );;
    gap> pi * form * TransposedMat( pi ) = InvariantBilinearForm( o ).matrix;
    true
    gap> pi2:= PermutationMat( (7,3)(8,4), 14, GF(2) );;
    gap> D:= pi2 * pi;;
    gap> tblow:= List( blow, x -> D * x * D^-1 );;
    gap> IsSubset( o, tblow );
    true
\end{verbatim}

Note that the central subgroup of order three in $\GU(7,2)$ consists of
scalar matrices.

\begin{verbatim}
    gap> omega:= DerivedSubgroup( o );;
    gap> IsSubset( omega, tblow );
    true
    gap> z:= Z(4) * One( g );;
    gap> tz:= D * BlownUpMat( b, z ) * D^-1;;
    gap> tz in omega;
    true
\end{verbatim}

Now we switch to a permutation representation of $S$,
and compute the conjugacy classes of prime element order in the subgroup $M$.
The latter is done in two steps, first class representatives of the
simple subgroup $U_7(2)$ of $M$ are computed, and then they are multiplied
with the scalars in $M$.

\begin{verbatim}
    gap> orbs:= Orbits( omega, NormedVectors( GF(2)^14 ), OnLines );;
    gap> List( orbs, Length );
    [ 8127, 8256 ]
    gap> omega:= Action( omega, orbs[1], OnLines );;
    gap> gens:= List( GeneratorsOfGroup( g ),
    >             x -> Permutation( D * BlownUpMat( b, x ) * D^-1, orbs[1] ) );;
    gap> g:= Group( gens );;
    gap> ccl:= ClassesOfPrimeOrder( g, Set( Factors( Size( g ) ) ),
    >                               TrivialSubgroup( g ) );;
    gap> tz:= Permutation( tz, orbs[1] );;
    gap> primereps:= List( ccl, Representative );;
    gap> Add( primereps, () );
    gap> reps:= Concatenation( List( primereps,
    >               x -> List( [ 0 .. 2 ], i -> x * tz^i ) ) );;
    gap> primereps:= Filtered( reps, x -> IsPrimeInt( Order( x ) ) );;
    gap> Length( primereps );
    48
\end{verbatim}

Finally, we apply \verb|UpperBoundFixedPointRatios| (see Section~\ref{groups})
to compute an upper bound for $\fpr(g,S/M)$, for $g \in S^{\times}$.

\begin{verbatim}
    gap> M:= ClosureGroup( g, tz );;
    gap> bccl:= List( primereps, x -> ConjugacyClass( M, x ) );;
    gap> UpperBoundFixedPointRatios( omega, [ bccl ], false );
    [ 1/2015, true ]
\end{verbatim}

Although some of the classes of $M$ in the list \verb|bccl| may be $S$-conjugate,
the entry \verb|true| in the second position of the result indicates
that in fact the \emph{exact} value for the maximum of $\fpr(g,S/M)$,
for $g \in S^{\times}$, has been computed.
This implies statement~(b).

We clean the workspace.

\begin{verbatim}
    gap> CleanWorkspace();
\end{verbatim}

\subsection{$O_{12}^+(3)$}\label{O12p3}

We show that the group $S = O_{12}^+(3) = \POmega^+(12,3)$ satisfies the
following.
\begin{enumerate}
\item[(a)]
    $S$ has a maximal subgroup $M$ of the type $N_S(\POmega^+(6,9))$,
    which has the structure $\POmega^+(6,9).[4]$.
\item[(b)]
    $\fpr(g,S/M) \leq 2/88\,209$ holds for all $g \in S^{\times}$.
\end{enumerate}

(This result is used in the proof of~\cite[Proposition~5.14]{BGK},
where it is shown that for $s \in S$ of order $205$,
$\M(S,s)$ consists of one reducible subgroup $G_8$
and at most two extension field type subgroups
of the type $N_S(\POmega^+(6,9))$.
By~\cite[Proposition~3.16]{GK}, $\fpr(g,S/G_8) \leq 19/3^5$ holds
for all $g \in S^{\times}$.
This implies
$\prop(g,s) \leq 19/3^5 + 2 \cdot 2/88\,209 = 6\,901/88\,209 < 1/3$.)

Statement~(a) follows from~\cite[Prop.~4.3.14]{KlL90}.

For statement~(b), we embed $\GO^+(6,9) \cong \Omega^+(6,9).2^2$
into $\SO^+(12,3) = 2.S.2$,
by replacing each element in $\F_9$ by the $2 \times 2$ matrix of the
induced $\F_3$-linear mapping w.r.t.~a suitable basis $(b_1, b_2)$.
We choose a basis with the property $b_1 = 1$ and $b_2^2 = 1 + b_2$,
because then the image of a symmetric matrix is again symmetric
(so the image of the invariant form is an invariant form for the image
of the group),
and apply an appropriate transformation to the images of the generators.

\begin{verbatim}
    gap> so:= SO(+1,12,3);;
    gap> Display( InvariantBilinearForm( so ).matrix );
     . 1 . . . . . . . . . .
     1 . . . . . . . . . . .
     . . 1 . . . . . . . . .
     . . . 2 . . . . . . . .
     . . . . 2 . . . . . . .
     . . . . . 2 . . . . . .
     . . . . . . 2 . . . . .
     . . . . . . . 2 . . . .
     . . . . . . . . 2 . . .
     . . . . . . . . . 2 . .
     . . . . . . . . . . 2 .
     . . . . . . . . . . . 2
    gap> g:= GO(+1,6,9);;
    gap> Z(9)^2 = Z(3)^0 + Z(9);
    true
    gap> b:= Basis( GF(9), [ Z(3)^0, Z(9) ] );
    Basis( GF(3^2), [ Z(3)^0, Z(3^2) ] )
    gap> blow:= List( GeneratorsOfGroup( g ), x -> BlownUpMat( b, x ) );;
    gap> m:= BlownUpMat( b, InvariantBilinearForm( g ).matrix );;
    gap> Display( m );
     . . 1 . . . . . . . . .
     . . . 1 . . . . . . . .
     1 . . . . . . . . . . .
     . 1 . . . . . . . . . .
     . . . . 2 . . . . . . .
     . . . . . 2 . . . . . .
     . . . . . . 2 . . . . .
     . . . . . . . 2 . . . .
     . . . . . . . . 2 . . .
     . . . . . . . . . 2 . .
     . . . . . . . . . . 2 .
     . . . . . . . . . . . 2
    gap> pi:= PermutationMat( (2,3), 12, GF(3) );;
    gap> tr:= IdentityMat( 12, GF(3) );;
    gap> tr{[3,4]}{[3,4]}:= [[1,-1],[1,1]]*Z(3)^0;;
    gap> D:= tr * pi;;
    gap> D * m * TransposedMat( D ) = InvariantBilinearForm( so ).matrix;
    true
    gap> tblow:= List( blow, x -> D * x * D^-1 );;
    gap> IsSubset( so, tblow );
    true
\end{verbatim}

The image of $\GO^+(6,9)$ under the embedding into $\SO^+(12,3)$
does not lie in $\Omega^+(12,3) = 2.S$,
so a factor of two is missing in $\GO^+(6,9) \cap 2.S$ for getting
(the preimage $2.M$ of) the required maximal subgroup $M$ of $S$.
Because of this, and also because currently it is time consuming
to compute the derived subgroup of $\SO^+(12,3)$,
we work with the upward extension $\PSO^+(12,3) = S.2$.
Note that $M$ extends to a maximal subgroup of $S.2$.

First we factor out the centre of $\SO^+(12,3)$,
and switch to a permutation representation of $S.2$.

\begin{verbatim}
    gap> orbs:= Orbits( so, NormedVectors( GF(3)^12 ), OnLines );;
    gap> List( orbs, Length );
    [ 88452, 88452, 88816 ]
    gap> act:= Action( so, orbs[1], OnLines );;
    gap> SetSize( act, Size( so ) / 2 );
\end{verbatim}

%
%

Next we rewrite the matrix generators for $\GO^+(6,9)$ accordingly,
and compute the normalizer in $S.2$ of the subgroup they generate;
this is the maximal subgroup $M.2$ we need.

\begin{verbatim}
    gap> u:= SubgroupNC( act,
    >            List( tblow, x -> Permutation( x, orbs[1], OnLines ) ) );;
    gap> n:= Normalizer( act, u );;
    gap> Size( n ) / Size( u );
    2
\end{verbatim}


Now we compute class representatives of prime order in $M.2$,
in a smaller faithful permutation representation,
and then the desired upper bound for $\fpr(g, S/M)$.


\begin{verbatim}
    gap> norbs:= Orbits( n, MovedPoints( n ) );;
    gap> List( norbs, Length );
    [ 58968, 29484 ]
    gap> hom:= ActionHomomorphism( n, norbs[2] );;
    gap> nact:= Image( hom );;
    gap> Size( nact ) = Size( n );
    true
    gap> ccl:= ClassesOfPrimeOrder( nact, Set( Factors( Size( nact ) ) ),
    >                               TrivialSubgroup( nact ) );;
    gap> Length( ccl );
    26
    gap> preim:= List( ccl,
    >        x -> PreImagesRepresentative( hom, Representative( x ) ) );;
    gap> pccl:= List( preim, x -> ConjugacyClass( n, x ) );;
    gap> for i in [ 1 .. Length( pccl ) ] do
    >      SetSize( pccl[i], Size( ccl[i] ) );
    >    od;
    gap> UpperBoundFixedPointRatios( act, [ pccl ], false );
    [ 2/88209, true ]
\end{verbatim}


Note that we have computed
$\max\{ \fpr(g,S.2/M.2), g \in S.2^{\times} \} \geq
 \max\{ \fpr(g,S.2/M.2), g \in S^{\times} \} =
 \max\{ \fpr(g,S/M), g \in S^{\times} \}$.

\subsection{$\ast$~$S_4(8)$}\label{S48}

We show that the group $S = S_4(8) = \Sp(4,8)$ satisfies the following.
\begin{enumerate}
\item[(a)]
    For $s \in S$ irreducible of order $65$,
    $\M(S,s)$ consists of two nonconjugate subgroups of the type
    $S_2(64).2 = \Sp(2,64).2 \cong L_2(64).2
    \cong O_4^-(8).2 = \Omega^-(4,8).2$.
\item[(b)]
    $\total(S,s) = 8/63$.
\end{enumerate}

By~\cite{Be00}, the only maximal subgroups of $S$ that contain $s$
are $O_4^-(8).2 = \SO^-(4,8)$ or of extension field type.
By~\cite[Prop.~4.3.10, 4.8.6]{KlL90}, there is one class of each of these
subgroups (which happen to be isomorphic).

These classes of subgroups induce different permutation characters.
One argument to see this is that the involutions in the outer half
of extension field type subgroup $S_2(64).2 < S_4(8)$ have
a two-dimensional fixed space,
whereas the outer involutions in $\SO^-(4,8)$ have a three-dimensional
fixed space.

The former statement can be seen by using a normal basis of the field
extension $\F_{64}/\F_8$, such that the action of the Frobenius
automorphism (which yields a suitable outer involution) is just a
double transposition on the basis vectors of the natural module for $S$.

\begin{verbatim}
    gap> sp:= SP(4,8);;
    gap> Display( InvariantBilinearForm( sp ).matrix );
     . . . 1
     . . 1 .
     . 1 . .
     1 . . .
    gap> z:= Z(64);;
    gap> f:= AsField( GF(8), GF(64) );;
    gap> repeat
    >      b:= Basis( f, [ z, z^8 ] );
    >      z:= z * Z(64);
    > until b <> fail;
    gap> sub:= SP(2,64);;
    gap> Display( InvariantBilinearForm( sub ).matrix );
     . 1
     1 .
    gap> ext:= Group( List( GeneratorsOfGroup( sub ),
    >                       x -> BlownUpMat( b, x ) ) );;
    gap> tr:= PermutationMat( (3,4), 4, GF(2) );;
    gap> conj:= ConjugateGroup( ext, tr );;
    gap> IsSubset( sp, conj );
    true
    gap> inv:= [[0,1,0,0],[1,0,0,0],[0,0,0,1],[0,0,1,0]] * Z(2);;
    gap> inv in sp;
    true
    gap> inv in conj;
    false
    gap> Length( NullspaceMat( inv - inv^0 ) );
    2
\end{verbatim}

The latter statement can be shown by looking at an outer involution in
$\SO^-(4,8)$.

\begin{verbatim}
    gap> so:= SO(-1,4,8);;
    gap> der:= DerivedSubgroup( so );;
    gap> x:= First( GeneratorsOfGroup( so ), x -> not x in der );;
    gap> x:= x^( Order(x)/2 );;
    gap> Length( NullspaceMat( x - x^0 ) );
    3
\end{verbatim}

%

The character table of $L_2(64).2$ is currently not available in the
{\GAP} Character Table Library, so we compute the possible permutation
characters with a combinatorial approach,
and show statement~(a).

\begin{verbatim}
    gap> CharacterTable( "L2(64).2" );
    fail
    gap> t:= CharacterTable( "S4(8)" );;
    gap> degree:= Size( t ) / ( 2 * Size( SL(2,64) ) );;
    gap> pi:= PermChars( t, rec( torso:= [ degree ] ) );
    [ Character( CharacterTable( "S4(8)" ), [ 2016, 0, 256, 32, 0, 36, 0, 8, 1, 
          0, 4, 0, 0, 0, 28, 28, 28, 0, 0, 0, 0, 0, 0, 36, 36, 36, 0, 0, 0, 0, 0, 
          0, 1, 1, 1, 0, 0, 0, 4, 4, 4, 0, 0, 0, 4, 4, 4, 0, 0, 0, 1, 1, 1, 0, 0, 
          0, 0, 0, 0, 0, 0, 0, 1, 1, 1, 1, 1, 1, 1, 1, 1, 1, 1, 1, 1, 1, 1, 1, 1, 
          1, 1, 1, 1 ] ), Character( CharacterTable( "S4(8)" ), 
        [ 2016, 256, 0, 32, 36, 0, 0, 8, 1, 4, 0, 28, 28, 28, 0, 0, 0, 0, 0, 0, 
          36, 36, 36, 0, 0, 0, 0, 0, 0, 0, 0, 0, 1, 1, 1, 4, 4, 4, 0, 0, 0, 4, 4, 
          4, 0, 0, 0, 1, 1, 1, 0, 0, 0, 1, 1, 1, 1, 1, 1, 1, 1, 1, 0, 0, 0, 0, 0, 
          0, 0, 0, 0, 1, 1, 1, 1, 1, 1, 1, 1, 1, 1, 1, 1 ] ) ]
    gap> spos:= Position( OrdersClassRepresentatives( t ), 65 );;
    gap> List( pi, x -> x[ spos ] );
    [ 1, 1 ]
\end{verbatim}

Now we compute $\total(S,s)$, which yields statement~(b).

\begin{verbatim}
    gap> Maximum( ApproxP( pi, spos ) );
    8/63
\end{verbatim}

We clean the workspace.

\begin{verbatim}
    gap> CleanWorkspace();
\end{verbatim}

\subsection{$S_6(2)$}\label{S62}

We show that the group $S = S_6(2) = \Sp(6,2)$ satisfies the following.
\begin{enumerate}
\item[(a)]
    $\total(S) = 4/7$,
    and this value is attained exactly for $\total(S,s)$
    with $s$ of order $9$.
\item[(b)]
    For $s \in S$ of order $9$,
    $\M(S,s)$ consists of one subgroup of the type
    $U_4(2).2 = \Omega^-(6,2).2$ and three conjugate subgroups
    of the type $L_2(8).3 = \Sp(2,8).3$.
\item[(c)]
    For $s \in S$ of order $9$,
    and $g \in S^{\times}$,
    we have $\prop(g,s) < 1/3$, except if $g$ is in one of the classes
    {\tt 2A} (the transvection class) or {\tt 3A}.
\item[(d)]
    For $s \in S$ of order $15$,
    and $g \in S^{\times}$,
    we have $\prop(g,s) < 1/3$, except if $g$ is in one of the classes
    {\tt 2A} or {\tt 2B}.
\item[(e)]
    $\prop(S) = 11/21$,
    and this value is attained exactly for $\prop(S,s)$
    with $s$ of order $15$.
\item[(f)]
    For all $s^{\prime} \in S$,
    we have $\prop(g,s^{\prime}) > 1/3$ for $g$ in at least two classes.
\item[(g)]
    The uniform spread of $S$ is at least two,
    with $s$ of order $9$.
\end{enumerate}

(Note that in this example, the optimal choice of $s$ w.r.t. $\total(S,s)$
is not optimal w.r.t. $\prop(S,s)$.)

Statement~(a) follows from the inspection of the primitive permutation
characters, cf.~Section~\ref{easyloop}.

\begin{verbatim}
    gap> t:= CharacterTable( "S6(2)" );;
    gap> ProbGenInfoSimple( t );
    [ "S6(2)", 4/7, 1, [ "9A" ], [ 4 ] ]
\end{verbatim}

Also statement~(b) follows from the information provided by the
character table of $S$ (cf.~\cite[p.~46]{CCN85}).

\begin{verbatim}
    gap> prim:= PrimitivePermutationCharacters( t );;
    gap> ord:= OrdersClassRepresentatives( t );;
    gap> spos:= Position( ord, 9 );;
    gap> filt:= PositionsProperty( prim, x -> x[ spos ] <> 0 );
    [ 1, 8 ]
    gap> Maxes( t ){ filt };
    [ "U4(2).2", "L2(8).3" ]
    gap> List( prim{ filt }, x -> x[ spos ] );
    [ 1, 3 ]
\end{verbatim}

Now we consider statement~(c).
For $s$ of order $9$ and $g$ in one of the classes {\tt 2A}, {\tt 3A},
we observe that $\prop(g,s) = \total(g,s)$ holds.
This is because exactly one maximal subgroup of $S$
contains both $s$ and $g$.
For all other elements $g$, we have even $\total(g,s) < 1/3$.

\begin{verbatim}
    gap> prim:= PrimitivePermutationCharacters( t );;
    gap> spos9:= Position( ord, 9 );;
    gap> approx9:= ApproxP( prim, spos9 );;
    gap> filt9:= PositionsProperty( approx9, x -> x >= 1/3 );
    [ 2, 6 ]
    gap> AtlasClassNames( t ){ filt9 };
    [ "2A", "3A" ]
    gap> approx9{ filt9 };
    [ 4/7, 5/14 ]
    gap> List( Filtered( prim, x -> x[ spos9 ] <> 0 ), x -> x{ filt9 } );
    [ [ 16, 10 ], [ 0, 0 ] ]
\end{verbatim}

Similarly, statement~(d) follows.
For $s$ of order $15$ and $g$ in one of the classes {\tt 2A}, {\tt 2B},
already the degree $36$ permutation character yields $\prop(g,s) \geq 1/3$.
And for all other elements $g$, again we have $\total(g,s) < 1/3$.

\begin{verbatim}
    gap> spos15:= Position( ord, 15 );;
    gap> approx15:= ApproxP( prim, spos15 );;
    gap> filt15:= PositionsProperty( approx15, x -> x >= 1/3 );
    [ 2, 3 ]
    gap> PositionsProperty( ApproxP( prim{ [ 2 ] }, spos15 ), x -> x >= 1/3 );
    [ 2, 3 ]
    gap> AtlasClassNames( t ){ filt15 };
    [ "2A", "2B" ]
    gap> approx15{ filt15 };
    [ 46/63, 8/21 ]
\end{verbatim}

For the remaining statements, we use explicit computations with $S$,
in the transitive degree $63$ permutation representation.
We start with a function that computes a transvection in $S_d(2)$;
note that the invariant bilinear form used for symplectic groups in {\GAP}
is described by a matrix with nonzero entries exactly in the positions
$(i,d+1-i)$, for $1 \leq i \leq d$.

\begin{verbatim}
    gap> transvection:= function( d )
    >     local mat;
    >     mat:= IdentityMat( d, Z(2) );
    >     mat{ [ 1, d ] }{ [ 1, d ] }:= [ [ 0, 1 ], [ 1, 0 ] ] * Z(2);
    >     return mat;
    > end;;
\end{verbatim}

First we compute, for statement~(d), the exact values $\prop(g,s)$ for $g$
in one of the classes {\tt 2A} or {\tt 2B}, and $s$ of order $15$.
Note that the classes {\tt 2A}, {\tt 2B} are the unique classes of the
lengths $63$ and $315$, respectively.

\begin{verbatim}
    gap> PositionsProperty( SizesConjugacyClasses( t ), x -> x in [ 63, 315 ] );
    [ 2, 3 ]
    gap> d:= 6;;
    gap> matgrp:= Sp(d,2);;
    gap> hom:= ActionHomomorphism( matgrp, NormedRowVectors( GF(2)^d ) );;
    gap> g:= Image( hom, matgrp );;
    gap> ResetGlobalRandomNumberGenerators();
    gap> repeat s15:= Random( g );
    >    until Order( s15 ) = 15;
    gap> 2A:= Image( hom, transvection( d ) );;
    gap> Size( ConjugacyClass( g, 2A ) );
    63
    gap> IsTransitive( g, MovedPoints( g ) );
    true
    gap> RatioOfNongenerationTransPermGroup( g, 2A, s15 );
    11/21
    gap> repeat 12C:= Random( g );
    >    until Order( 12C ) = 12 and Size( Centralizer( g, 12C ) ) = 12;
    gap> 2B:= 12C^6;;
    gap> Size( ConjugacyClass( g, 2B ) );
    315
    gap> RatioOfNongenerationTransPermGroup( g, 2B, s15 );
    8/21
\end{verbatim}

For statement~(e), we compute $\prop(g, s^{\prime})$,
for a transvection $g$ and class representatives $s^{\prime}$ of $S$.
It turns out that the minimum is $11/21$,
and it is attained for exactly one $s^{\prime}$;
by the above, this element has order $15$.

\begin{verbatim}
    gap> ccl:= ConjugacyClasses( g );;
    gap> reps:= List( ccl, Representative );;
    gap> nongen2A:= List( reps,
    >        x -> RatioOfNongenerationTransPermGroup( g, 2A, x ) );;
    gap> min:= Minimum( nongen2A );
    11/21
    gap> Number( nongen2A, x -> x = min );
    1
\end{verbatim}

For statement~(f), we show that for any choice of $s^{\prime}$,
at least two of the values $\prop(g,s^{\prime})$,
with $g$ in the classes {\tt 2A}, {\tt 2B}, or {\tt 3A},
are larger than $1/3$.

\begin{verbatim}
    gap> nongen2B:= List( reps,
    >        x -> RatioOfNongenerationTransPermGroup( g, 2B, x ) );;
    gap> 3A:= s15^5;;
    gap> nongen3A:= List( reps,
    >        x -> RatioOfNongenerationTransPermGroup( g, 3A, x ) );;
    gap> bad:= List( [ 1 .. NrConjugacyClasses( t ) ],
    >                i -> Number( [ nongen2A, nongen2B, nongen3A ],
    >                             x -> x[i] > 1/3 ) );;
    gap> Minimum( bad );
    2
\end{verbatim}

Finally, for statement~(g), we have to consider only the case that the
two elements $x$, $y$ are transvections.

\begin{verbatim}
    gap> PositionsProperty( approx9, x -> x + approx9[2] >= 1 );
    [ 2 ]
\end{verbatim}

We use the random approach described in Section~\ref{groups}.

\begin{verbatim}
    gap> repeat s9:= Random( g );
    >    until Order( s9 ) = 9;
    gap> RandomCheckUniformSpread( g, [ 2A, 2A ], s9, 20 );
    true
\end{verbatim}

\subsection{$S_8(2)$}\label{S82}

We show that the group $S = S_8(2)$ satisfies the following.
\begin{enumerate}
\item[(a)]
    For $s \in S$ of order $17$,
    $\M(S,s)$ consists of one subgroup of each of the types
    $O_8^-(2).2 = \Omega^-(8,2).2$, $S_4(4).2 = \Sp(4,4).2$,
    and $L_2(17) = \PSL(2,17)$.
\item[(b)]
    For $s \in S$ of order $17$,
    and $g \in S^{\times}$,
    we have $\prop(g,s) < 1/3$, except if $g$ is a transvection.
\item[(c)]
    The uniform spread of $S$ is at least two,
    with $s$ of order $17$.
\end{enumerate}

Statement~(a) follows from the list of maximal subgroups of $S$
in~\cite[p.~123]{CCN85},
and the fact that $1_H^S(s) = 1$ holds for each $H \in \M(S,s)$.
Note that $17$ divides the indices of the maximal subgroups of the types
$O_8^+(2).2$ and $2^7 : S_6(2)$ in $S$,
and obviously $17$ does not divide the orders of the remaining maximal
subgroups.

The permutation characters induced from the first two subgroups
are uniquely determined by the ordinary character tables.
The permutation character induced from the last subgroup is
uniquely determined if one considers also the corresponding Brauer tables;
the correct class fusion is stored in the {\GAP} Character Table Library,
see~\cite{AmbigFus}.

\begin{verbatim}
    gap> t:= CharacterTable( "S8(2)" );;
    gap> pi1:= PossiblePermutationCharacters( CharacterTable( "O8-(2).2" ), t );;
    gap> pi2:= PossiblePermutationCharacters( CharacterTable( "S4(4).2" ), t );;
    gap> pi3:= [ TrivialCharacter( CharacterTable( "L2(17)" ) )^t ];;
    gap> prim:= Concatenation( pi1, pi2, pi3 );;
    gap> Length( prim );
    3
    gap> spos:= Position( OrdersClassRepresentatives( t ), 17 );;
    gap> List( prim, x -> x[ spos ] );
    [ 1, 1, 1 ]
\end{verbatim}

For statement~(b),
we observe that $\total(g,s) < 1/3$ if $g$ is not a transvection,
and that $\prop(g,s) = \total(g,s)$ for transvections $g$
because exactly one of the three permutation characters is nonzero
on both $s$ and the class of transvections.

\begin{verbatim}
    gap> approx:= ApproxP( prim, spos );;
    gap> PositionsProperty( approx, x -> x >= 1/3 );
    [ 2 ]
    gap> Number( prim, pi -> pi[2] <> 0 and pi[ spos ] <> 0 );
    1
    gap> approx[2];
    8/15
\end{verbatim}

In statement~(c), we have to consider only the case that the
two elements $x$, $y$ are transvections.

\begin{verbatim}
    gap> PositionsProperty( approx, x -> x + approx[2] >= 1 );
    [ 2 ]
\end{verbatim}

We use the random approach described in Section~\ref{groups}.

\begin{verbatim}
    gap> d:= 8;;
    gap> matgrp:= Sp(d,2);;
    gap> hom:= ActionHomomorphism( matgrp, NormedRowVectors( GF(2)^d ) );;
    gap> x:= Image( hom, transvection( d ) );;
    gap> g:= Image( hom, matgrp );;
    gap> C:= ConjugacyClass( g, x );;  Size( C );
    255
    gap> ResetGlobalRandomNumberGenerators();
    gap> repeat s:= Random( g );
    >    until Order( s ) = 17;
    gap> RandomCheckUniformSpread( g, [ x, x ], s, 20 );
    true
\end{verbatim}

\subsection{$\ast$~$S_{10}(2)$}\label{S102}

We show that the group $S = S_{10}(2)$ satisfies the following.
\begin{enumerate}
\item[(a)]
    For $s \in S$ of order $33$,
    $\M(S,s)$ consists of one subgroup of each of the types
    $\Omega^-(10,2).2$ and $L_2(32).5 = \Sp(2,32).5$.
\item[(b)]
    For $s \in S$ of order $33$,
    and $g \in S^{\times}$,
    we have $\prop(g,s) < 1/3$, except if $g$ is a transvection.
\item[(c)]
    The uniform spread of $S$ is at least two,
    with $s$ of order $33$.
\end{enumerate}

By~\cite{Be00}, the only maximal subgroups of $S$ that contain $s$
have the types stated in~(a),
and by~\cite[Prop.~4.3.10 and 4.8.6]{KlL90}, there is exactly one class
of each of these subgroups.

We compute the values $\total( g, s )$, for all $g \in S^{\times}$.

\begin{verbatim}
    gap> t:= CharacterTable( "S10(2)" );;
    gap> pi1:= PossiblePermutationCharacters( CharacterTable( "O10-(2).2" ), t );;
    gap> pi2:= PossiblePermutationCharacters( CharacterTable( "L2(32).5" ), t );;
    gap> prim:= Concatenation( pi1, pi2 );;  Length( prim );
    2
    gap> spos:= Position( OrdersClassRepresentatives( t ), 33 );;
    gap> approx:= ApproxP( prim, spos );;
\end{verbatim}

For statement~(b),
we observe that $\total(g,s) < 1/3$ if $g$ is not a transvection,
and that $\prop(g,s) = \total(g,s)$ for transvections $g$
because exactly one of the two permutation characters is nonzero
on both $s$ and the class of transvections.

\begin{verbatim}
    gap> PositionsProperty( approx, x -> x >= 1/3 );
    [ 2 ]
    gap> Number( prim, pi -> pi[2] <> 0 and pi[ spos ] <> 0 );
    1
    gap> approx[2];
    16/31
\end{verbatim}

In statement~(c), we have to consider only the case that the
two elements $x$, $y$ are transvections.
We use the random approach described in Section~\ref{groups}.

\begin{verbatim}
    gap> d:= 10;;
    gap> matgrp:= Sp(d,2);;
    gap> hom:= ActionHomomorphism( matgrp, NormedRowVectors( GF(2)^d ) );;
    gap> x:= Image( hom, transvection( d ) );;
    gap> g:= Image( hom, matgrp );;
    gap> C:= ConjugacyClass( g, x );;  Size( C );
    1023
    gap> ResetGlobalRandomNumberGenerators();
    gap> repeat s:= Random( g );
    >    until Order( s ) = 33;
    gap> RandomCheckUniformSpread( g, [ x, x ], s, 20 );
    true
\end{verbatim}

\subsection{$U_4(2)$}\label{U42}

We show that $S = U_4(2) = \SU(4,2) \cong S_4(3) = \PSp(4,3)$
satisfies the following.
\begin{enumerate}
\item[(a)]
    $\total(S) = 21/40$,
    and this value is attained exactly for $\total(S,s)$
    with $s$ of order $12$.
\item[(b)]
    For $s \in S$ of order $9$,
    $\M(S,s)$ consists of two groups,
    of the types $3^{1+2}_+ \colon 2A_4 = \GU(3,2)$ and $3^3 \colon S_4$,
    respectively.
\item[(c)]
    $\prop(S) = 2/5$,
    and this value is attained exactly for $\prop(S,s)$
    with $s$ of order $9$.
\item[(d)]
    The uniform spread of $S$ is at least three,
    with $s$ of order $9$.
\item[(e)]
    $\total^{\prime}(\Aut(S),s) = 7/20$.
\end{enumerate}


(Note that in this example, the optimal choice of $s$ w.r.t. $\total(S,s)$
is not optimal w.r.t. $\prop(S,s)$.)

Statement~(a) follows from inspection of the primitive permutation
characters, cf.~Section~\ref{easyloop}.

\begin{verbatim}
    gap> t:= CharacterTable( "U4(2)" );;
    gap> ProbGenInfoSimple( t );
    [ "U4(2)", 21/40, 1, [ "12A" ], [ 2 ] ]
\end{verbatim}

Statement~(b) can be read off from the permutation characters,
and the fact that the only classes of maximal subgroups that contain
elements of order $9$ consist of groups of the structures
$3^{1+2}_+:2A_4$ and $3^3:S_4$,
see~\cite[p.~26]{CCN85}.

\begin{verbatim}
    gap> OrdersClassRepresentatives( t );
    [ 1, 2, 2, 3, 3, 3, 3, 4, 4, 5, 6, 6, 6, 6, 6, 6, 9, 9, 12, 12 ]
    gap> prim:= PrimitivePermutationCharacters( t );
    [ Character( CharacterTable( "U4(2)" ), [ 27, 3, 7, 0, 0, 9, 0, 3, 1, 2, 0, 
          0, 3, 3, 0, 1, 0, 0, 0, 0 ] ), Character( CharacterTable( "U4(2)" ), 
        [ 36, 12, 8, 0, 0, 6, 3, 0, 2, 1, 0, 0, 0, 0, 3, 2, 0, 0, 0, 0 ] ), 
      Character( CharacterTable( "U4(2)" ), [ 40, 8, 0, 13, 13, 4, 4, 4, 0, 0, 5, 
          5, 2, 2, 2, 0, 1, 1, 1, 1 ] ), Character( CharacterTable( "U4(2)" ), 
        [ 40, 16, 4, 4, 4, 1, 7, 0, 2, 0, 4, 4, 1, 1, 1, 1, 1, 1, 0, 0 ] ), 
      Character( CharacterTable( "U4(2)" ), [ 45, 13, 5, 9, 9, 6, 3, 1, 1, 0, 1, 
          1, 4, 4, 1, 2, 0, 0, 1, 1 ] ) ]
\end{verbatim}

For statement~(c),
we use a primitive permutation representation on $40$ points
that occurs in the natural action of $\SU(4,2)$.

\begin{verbatim}
    gap> g:= SU(4,2);;
    gap> orbs:= Orbits( g, NormedRowVectors( GF(4)^4 ), OnLines );;
    gap> List( orbs, Length );
    [ 45, 40 ]
    gap> g:= Action( g, orbs[2], OnLines );;
\end{verbatim}

First we show that for $s$ of order $9$, $\prop(S,s) = 2/5$ holds.
For that, we have to consider only $\prop(g,s)$,
with $g$ in one of the classes {\tt 2A} (of length $45$)
and {\tt 3A} (of length $40$);
since the class {\tt 3B} contains the inverses of the elements in
the class {\tt 3A}, we need not test it.

\begin{verbatim}
    gap> spos:= Position( OrdersClassRepresentatives( t ), 9 );
    17
    gap> approx:= ApproxP( prim, spos );
    [ 0, 3/5, 1/10, 17/40, 17/40, 1/8, 11/40, 1/10, 1/20, 0, 9/40, 9/40, 3/40, 
      3/40, 3/40, 1/40, 1/20, 1/20, 1/40, 1/40 ]
    gap> badpos:= PositionsProperty( approx, x -> x >= 2/5 );
    [ 2, 4, 5 ]
    gap> PowerMap( t, 2 )[4];
    5
    gap> OrdersClassRepresentatives( t );
    [ 1, 2, 2, 3, 3, 3, 3, 4, 4, 5, 6, 6, 6, 6, 6, 6, 9, 9, 12, 12 ]
    gap> SizesConjugacyClasses( t );
    [ 1, 45, 270, 40, 40, 240, 480, 540, 3240, 5184, 360, 360, 720, 720, 1440, 
      2160, 2880, 2880, 2160, 2160 ]
\end{verbatim}

A representative $g$ of a class of length $40$ can be found as the third
power of any order $9$ element.

\begin{verbatim}
    gap> PowerMap( t, 3 )[ spos ];
    4
    gap> ResetGlobalRandomNumberGenerators();
    gap> repeat s:= Random( g );
    >    until Order( s ) = 9;
    gap> Size( ConjugacyClass( g, s^3 ) );
    40
    gap> prop:= RatioOfNongenerationTransPermGroup( g, s^3, s );
    13/40
\end{verbatim}

Next we examine $g$ in the class {\tt 2A}.

\begin{verbatim}
    gap> repeat x:= Random( g ); until Order( x ) = 12;
    gap> Size( ConjugacyClass( g, x^6 ) );
    45
    gap> prop:= RatioOfNongenerationTransPermGroup( g, x^6, s );
    2/5
\end{verbatim}

Finally, we compute that for $s$ of order different from $9$
and $g$ in the class {\tt 2A}, $\prop(g,s)$ is larger than $2/5$.

\begin{verbatim}
    gap> ccl:= List( ConjugacyClasses( g ), Representative );;
    gap> SortParallel( List( ccl, Order ), ccl );
    gap> List( ccl, Order );
    [ 1, 2, 2, 3, 3, 3, 3, 4, 4, 5, 6, 6, 6, 6, 6, 6, 9, 9, 12, 12 ]
    gap> prop:= List( ccl, r -> RatioOfNongenerationTransPermGroup( g, x^6, r ) );
    [ 1, 1, 1, 1, 1, 1, 1, 1, 1, 5/9, 1, 1, 1, 1, 1, 1, 2/5, 2/5, 7/15, 7/15 ]
    gap> Minimum( prop );
    2/5
\end{verbatim}

In order to show statement~(d),
we have to consider triples $(x_1, x_2, x_3)$
with $x_i$ of prime order and $\sum_{i=1}^3 \prop(x_i,s) \geq 1$.
This means that it suffices to check $x$ in the class {\tt 2A},
$y$ in ${\tt 2A} \cup {\tt 3A}$,
and $z$ in ${\tt 2A} \cup {\tt 3A} \cup {\tt 3D}$.

\begin{verbatim}
    gap> approx[2]:= 2/5;;
    gap> approx[4]:= 13/40;;
    gap> primeord:= PositionsProperty( OrdersClassRepresentatives( t ),
    >                                  IsPrimeInt );
    [ 2, 3, 4, 5, 6, 7, 10 ]
    gap> RemoveSet( primeord, 5 );
    gap> primeord;
    [ 2, 3, 4, 6, 7, 10 ]
    gap> approx{ primeord };
    [ 2/5, 1/10, 13/40, 1/8, 11/40, 0 ]
    gap> AtlasClassNames( t ){ primeord };
    [ "2A", "2B", "3A", "3C", "3D", "5A" ]
    gap> triples:= Filtered( UnorderedTuples( primeord, 3 ),
    >                  t -> Sum( approx{ t } ) >= 1 );
    [ [ 2, 2, 2 ], [ 2, 2, 4 ], [ 2, 2, 7 ], [ 2, 4, 4 ], [ 2, 4, 7 ] ]
\end{verbatim}

We use the random approach described in Section~\ref{groups}.

\begin{verbatim}
    gap> repeat 6E:= Random( g );
    >    until Order( 6E ) = 6 and Size( Centralizer( g, 6E ) ) = 18;
    gap> 2A:= 6E^3;;
    gap> 3A:= s^3;;
    gap> 3D:= 6E^2;;
    gap> RandomCheckUniformSpread( g, [ 2A, 2A, 2A ], s, 50 );
    true
    gap> RandomCheckUniformSpread( g, [ 2A, 2A, 3A ], s, 50 );
    true
    gap> RandomCheckUniformSpread( g, [ 3D, 2A, 2A ], s, 50 );
    true
    gap> RandomCheckUniformSpread( g, [ 2A, 3A, 3A ], s, 50 );
    true
    gap> RandomCheckUniformSpread( g, [ 3D, 3A, 2A ], s, 50 );
    true
\end{verbatim}

Statement~(e) can be proved using \verb|ProbGenInfoAlmostSimple|,
cf. Section~\ref{easyloopaut}.

\begin{verbatim}
    gap> t:= CharacterTable( "U4(2)" );;
    gap> t2:= CharacterTable( "U4(2).2" );;
    gap> spos:= PositionsProperty( OrdersClassRepresentatives( t ), x -> x = 9 );;
    gap> ProbGenInfoAlmostSimple( t, t2, spos );
    [ "U4(2).2", 7/20, [ "9AB" ], [ 2 ] ]
\end{verbatim}

\subsection{$U_4(3)$}\label{U43}

We show that $S = U_4(3) = \PSU(4,3)$ satisfies the following.
\begin{enumerate}
\item[(a)]
    $\total(S) = 53/153$,
    and this value is attained exactly for $\total(S,s)$
    with $s$ of order $7$.
\item[(b)]
    For $s \in S$ of order $7$,
    $\M(S,s)$ consists of two nonconjugate groups of the type $L_3(4)$,
    one group of the type $U_3(3)$, and four pairwise nonconjugate
    groups of the type $A_7$.
\item[(c)]
    $\prop(S) = 43/135$,
    and this value is attained exactly for $\prop(S,s)$
    with $s$ of order $7$.
\item[(d)]
    The uniform spread of $S$ is at least three,
    with $s$ of order $7$.
\item[(e)]
    The preimage of $s$ in the matrix group $\SU(4,3) \cong 4.U_4(3)$
    has order $28$,
    the preimages of the groups in $\M(S,s)$ have the structures
    $4_2.L_3(4)$, $4 \times U_3(3) \cong \GU(3,3)$, and $4.A_7$
    (the latter being a central product of a cyclic group of order four
    and $2.A_7$).
\item[(f)]
    $\prop^{\prime}(S.2_1,s) = 13/27$,
    $\total^{\prime}(S.2_2) = 1/3$,
    and $\total^{\prime}(S.2_3) = 31/162$,
    with $s$ of order $7$ in each case.
\end{enumerate}

Statement~(a) follows from inspection of the primitive permutation
characters, cf.~Section~\ref{easyloop}.

\begin{verbatim}
    gap> t:= CharacterTable( "U4(3)" );;
    gap> ProbGenInfoSimple( t );
    [ "U4(3)", 53/135, 2, [ "7A" ], [ 7 ] ]
\end{verbatim}

Statement~(b) can be read off from the permutation characters,
and the fact that the only classes of maximal subgroups that contain
elements of order $7$ consist of groups of the structures
as claimed,
see~\cite[p.~52]{CCN85}.

\begin{verbatim}
    gap> prim:= PrimitivePermutationCharacters( t );;
    gap> spos:= Position( OrdersClassRepresentatives( t ), 7 );
    13
    gap> List( Filtered( prim, x -> x[ spos ] <> 0 ), l -> l{ [ 1, spos ] } );
    [ [ 162, 1 ], [ 162, 1 ], [ 540, 1 ], [ 1296, 1 ], [ 1296, 1 ], [ 1296, 1 ], 
      [ 1296, 1 ] ]
\end{verbatim}

In order to show statement~(c) (which then implies statement~(d)),
we use a permutation representation on $112$ points.
It corresponds to an orbit of one-dimensional subspaces in the
natural module of $\Omega^-(6,3) \cong S$.

\begin{verbatim}
    gap> matgrp:= DerivedSubgroup( SO( -1, 6, 3 ) );;
    gap> orbs:= Orbits( matgrp, NormedRowVectors( GF(3)^6 ), OnLines );;
    gap> List( orbs, Length );
    [ 126, 126, 112 ]
    gap> G:= Action( matgrp, orbs[3], OnLines );;
\end{verbatim}

It is sufficient to compute $\prop(g,s)$, for involutions $g \in S$.

\begin{verbatim}
    gap> approx:= ApproxP( prim, spos );
    [ 0, 53/135, 1/10, 1/24, 1/24, 7/45, 4/45, 1/27, 1/36, 1/90, 1/216, 1/216, 
      7/405, 7/405, 1/270, 0, 0, 0, 0, 1/270 ]
    gap> Filtered( approx, x -> x >= 43/135 );
    [ 53/135 ]
    gap> OrdersClassRepresentatives( t );
    [ 1, 2, 3, 3, 3, 3, 4, 4, 5, 6, 6, 6, 7, 7, 8, 9, 9, 9, 9, 12 ]
    gap> ResetGlobalRandomNumberGenerators();
    gap> repeat g:= Random( G ); until Order(g) = 2;
    gap> repeat s:= Random( G );
    >    until Order(s) = 7;
    gap> bad:= RatioOfNongenerationTransPermGroup( G, g, s );
    43/135
    gap> bad < 1/3;
    true
\end{verbatim}

Statement~(e) can be shown easily with character-theoretic methods,
as follows.
Since $\SU(4,3)$ is a Schur cover of $S$ and the groups in $\M(S,s)$
are simple, only very few possibilities have to be checked.
The Schur multiplier of $U_3(3)$ is trivial (see, e.~g., \cite[p.~14]{CCN85}),
so the preimage in $\SU(4,3)$ is a direct product of $U_3(3)$ and the
centre of $\SU(4,3)$.
Neither $L_3(4)$ nor its double cover $2.L_3(4)$ can be a subgroup of
$\SU(4,3)$, so the preimage of $L_3(4)$ must be a Schur cover of $L_3(4)$,
i.~e., it must have either the type $4_1.L_3(4)$ or $4_2.L_3(4)$
(see~\cite[p.~23]{CCN85});
only the type $4_2.L_3(4)$ turns out to be possible.

\begin{verbatim}
    gap> 4t:= CharacterTable( "4.U4(3)" );;
    gap> Length( PossibleClassFusions( CharacterTable( "L3(4)" ), 4t ) );
    0
    gap> Length( PossibleClassFusions( CharacterTable( "2.L3(4)" ), 4t ) );
    0
    gap> Length( PossibleClassFusions( CharacterTable( "4_1.L3(4)" ), 4t ) );
    0
    gap> Length( PossibleClassFusions( CharacterTable( "4_2.L3(4)" ), 4t ) );
    4
\end{verbatim}

As for the preimage of the $A_7$ type subgroups,
we first observe that the double cover of $A_7$ cannot be a subgroup of the
double cover of $S$,
so the preimage of $A_7$ in the double cover of $U_4(3)$ is a direct product
$2 \times A_7$.
The group $\SU(4,3)$ does not contain $A_7$ type subgroups,
thus the $A_7$ type subgroups in $2.U_4(3)$ lift to double covers of $A_7$
in $\SU(4,3)$.
This proves the claimed structure.

\begin{verbatim}
    gap> 2t:= CharacterTable( "2.U4(3)" );;
    gap> Length( PossibleClassFusions( CharacterTable( "2.A7" ), 2t ) );
    0
    gap> Length( PossibleClassFusions( CharacterTable( "A7" ), 4t ) );
    0
\end{verbatim}

For statement~(f), we consider automorphic extensions of $S$.
The bound for $S.2_3$ has been computed in Section~\ref{easyloopaut}.
That for $S.2_2$ can be computed form the fact that the classes of
maximal subgroups of $S.2_2$ containing $s$ of order $7$ are
$S$, one class of $U_3(3).2$ type subgroups, and two classes of $S_7$ type
subgroups which induce the same permutation character
(see~\cite[p.~52]{CCN85}).

\begin{verbatim}
    gap> t2:= CharacterTable( "U4(3).2_2" );;
    gap> pi1:= PossiblePermutationCharacters( CharacterTable( "U3(3).2" ), t2 );
    [ Character( CharacterTable( "U4(3).2_2" ), [ 540, 12, 54, 0, 0, 9, 8, 0, 0, 
          6, 0, 0, 1, 2, 0, 0, 0, 2, 0, 24, 4, 0, 0, 0, 0, 0, 0, 3, 2, 0, 4, 0, 
          0, 0 ] ) ]
    gap> pi2:= PossiblePermutationCharacters( CharacterTable( "A7.2" ), t2 );
    [ Character( CharacterTable( "U4(3).2_2" ), [ 1296, 48, 0, 27, 0, 9, 0, 4, 1, 
          0, 3, 0, 1, 0, 0, 0, 0, 0, 216, 24, 0, 4, 0, 0, 0, 9, 0, 3, 0, 1, 0, 1, 
          0, 0 ] ) ]
    gap> prim:= Concatenation( pi1, pi2, pi2 );;
    gap> outer:= Difference(
    >      PositionsProperty( OrdersClassRepresentatives( t2 ), IsPrimeInt ),
    >      ClassPositionsOfDerivedSubgroup( t2 ) );;
    gap> spos:= Position( OrdersClassRepresentatives( t2 ), 7 );;
    gap> Maximum( ApproxP( prim, spos ){ outer } );
    1/3
\end{verbatim}

Finally, Section~\ref{easyloopaut} shows that the character tables are
not sufficient for what we need, so we compute the exact proportion of
nongeneration for $U_4(3).2_1 \cong \SO^-(6,3)$.

\begin{verbatim}
    gap> matgrp:= SO( -1, 6, 3 );
    SO(-1,6,3)
    gap> orbs:= Orbits( matgrp, NormedRowVectors( GF(3)^6 ), OnLines );;
    gap> List( orbs, Length );
    [ 126, 126, 112 ]
    gap> G:= Action( matgrp, orbs[3], OnLines );;
    gap> repeat s:= Random( G );
    >    until Order( s ) = 7;
    gap> repeat
    >      repeat 2B:= Random( G ); until Order( 2B ) mod 2 = 0;
    >      2B:= 2B^( Order( 2B ) / 2 );
    >      c:= Centralizer( G, 2B );
    >    until Size( c ) = 12096;
    gap> RatioOfNongenerationTransPermGroup( G, 2B, s );
    13/27
    gap> repeat
    >      repeat 2C:= Random( G ); until Order( 2C ) mod 2 = 0;
    >      2C:= 2C^( Order( 2C ) / 2 );
    >      c:= Centralizer( G, 2C );
    >    until Size( c ) = 1440;
    gap> RatioOfNongenerationTransPermGroup( G, 2C, s );
    0
\end{verbatim}

\subsection{$U_6(3)$}\label{U63}

We show that $S = U_6(3) = \PSU(6,3)$ satisfies the following.
\begin{enumerate}
\item[(a)]
    For $s \in S$ of the type $1 \perp 5$
    (i.~e., the preimage of $s$ in $2.S = \SU(6,3)$ decomposes the natural
    $6$-dimensional module for $2.S$ into an orthogonal sum of two
    irreducible modules of the dimensions $1$ and $5$, respectively)
    and of order $(3^5 + 1)/2 = 122$,
    $\M(S,s)$ consists of one group of the type $2 \times U_5(3)$,
    which lifts to a subgroup of the type $4 \times U_5(3) = \GU(5,3)$
    in $2.S$.
    (The preimage of $s$ in $2.S$ has order $3^5 + 1 = 244$.)
\item[(b)]
    $\total(S,s) = 353/3\,159$.
\end{enumerate}

By~\cite{MSW94}, the only maximal subgroup of $S$ that contains $s$
is the stabilizer $H \cong 2 \times U_5(3)$ of the orthogonal decomposition.
This proves statement~(a).

The character table of $S$ is currently not available in the {\GAP}
Character Table Library.
We consider the permutation action of $S$ on the orbit of the stabilized
$1$-space.
So $M$ can be taken as a point stabilizer in this action.

\begin{verbatim}
    gap> CharacterTable( "U6(3)" );
    fail
    gap> g:= SU(6,3);;
    gap> orbs:= Orbits( g, NormedRowVectors( GF(9)^6 ), OnLines );;
    gap> List( orbs, Length );
    [ 22204, 44226 ]
    gap> repeat x:= PseudoRandom( g ); until Order( x ) = 244;
    gap> List( orbs, o -> Number( o, v -> OnLines( v, x ) = v ) );
    [ 0, 1 ]
    gap> g:= Action( g, orbs[2], OnLines );;
    gap> M:= Stabilizer( g, 1 );;
\end{verbatim}

Then we compute a list of elements in $M$ that covers the conjugacy classes
of prime element order, from which the numbers of fixed points
and thus $\max\{ \fpr( S/M, g ); g \in M^{\times} \} = \total( S, s )$
can be derived.
This way we avoid completely to check the $S$-conjugacy
of elements (class representatives of Sylow subgroups in $M$).

\begin{verbatim}
    gap> elms:= [];;
    gap> for p in Set( Factors( Size( M ) ) ) do
    >      syl:= SylowSubgroup( M, p );
    >      Append( elms, Filtered( PcConjugacyClassReps( syl ),
    >                              r -> Order( r ) = p ) );
    >    od;
    gap> 1 - Minimum( List( elms, NrMovedPoints ) ) / Length( orbs[2] );
    353/3159
\end{verbatim}


%

\subsection{$U_8(2)$}\label{U82}

We show that $S = U_8(2) = \SU(8,2)$ satisfies the following.
\begin{enumerate}
\item[(a)]
    For $s \in S$ of the type $1 \perp 7$
    (i.~e., $s$ decomposes the natural $8$-dimensional module for
    $S$ into an orthogonal sum of two irreducible modules
    of the dimensions $1$ and $7$, respectively) and of order
    $2^7 + 1 = 129$,
    $\M(S,s)$ consists of one group of the type $3 \times U_7(2) = \GU(7,2)$.
\item[(b)]
    $\total(S,s) = 2\,753/10\,880$.
\end{enumerate}

By~\cite{MSW94}, the only maximal subgroup of $S$ that contains $s$
is the stabilizer $M \cong \GU(7,2)$ of the orthogonal decomposition.
This proves statement~(a).

The character table of $S$ is currently not available in the {\GAP}
Character Table Library.
We proceed exactly as in Section~\ref{U63} in order to prove statement~(b).

\begin{verbatim}
    gap> CharacterTable( "U8(2)" );
    fail
    gap> g:= SU(8,2);;
    gap> orbs:= Orbits( g, NormedRowVectors( GF(4)^8 ), OnLines );;
    gap> List( orbs, Length );
    [ 10965, 10880 ]
    gap> repeat x:= PseudoRandom( g ); until Order( x ) = 129;
    gap> List( orbs, o -> Number( o, v -> OnLines( v, x ) = v ) );
    [ 0, 1 ]
    gap> g:= Action( g, orbs[2], OnLines );;
    gap> M:= Stabilizer( g, 1 );;
    gap> elms:= [];;
    gap> for p in Set( Factors( Size( M ) ) ) do
    >      syl:= SylowSubgroup( M, p );
    >      Append( elms, Filtered( PcConjugacyClassReps( syl ),
    >                              r -> Order( r ) = p ) );
    >    od;
    gap> Length( elms );
    611
    gap> 1 - Minimum( List( elms, NrMovedPoints ) ) / Length( orbs[2] );
    2753/10880
\end{verbatim}


%

\tthdump{\addcontentsline{toc}{section}{References}}

\bibliographystyle{amsalpha}
\newcommand{\etalchar}[1]{$^{#1}$}
\providecommand{\bysame}{\leavevmode\hbox to3em{\hrulefill}\thinspace}
\providecommand{\MR}{\relax\ifhmode\unskip\space\fi MR }
\providecommand{\MRhref}[2]{%
  \href{http://www.ams.org/mathscinet-getitem?mr=#1}{#2}
}
\providecommand{\href}[2]{#2}


\end{document}